\documentclass[11pt, reqno]{amsart}

\usepackage[usenames,dvipsnames]{color}
\usepackage[OT2, T1]{fontenc}
\usepackage{url}
\usepackage{amsmath}
\usepackage{array}
\usepackage{graphicx}
\usepackage{amssymb}
\usepackage{amsthm}
\definecolor{link}{RGB}{11,0,128}
\usepackage[colorlinks=true,citecolor=NavyBlue,linkcolor=OliveGreen,urlcolor=link]{hyperref}
\usepackage{hyperref}
\usepackage{paralist}
\usepackage{colonequals}			% for nice :=
\usepackage{color}
\usepackage{mathabx}			% for \widec
\usepackage{amscd}
\usepackage[all,cmtip]{xy}			% for commutative diagrams
\usepackage{sseq}				% for spectral sequences
\usepackage{verbatim}
\usepackage{parskip}
\usepackage{microtype}
\usepackage[in]{fullpage}
\usepackage{mathrsfs} 			% for \mathscr (script letters)
\usepackage{cleveref}
\usepackage{enumitem}
\usepackage{moreenum}			% for enumeration via Greek letters
\usepackage[alphabetic,lite]{amsrefs} 	% for bibliography
\usepackage{stmaryrd}			% this is for \llbracket and \rrbracket
\usepackage[foot]{amsaddr}
\usepackage{subfiles}
\usepackage{ragged2e}
\usepackage{stackengine}

\newcommand{\bG}{\mathbb{G}}

\newcommand{\bP}{\mathbb{P}}

\newcommand{\cO}{\mathcal{O}}

\DeclareSymbolFont{cyrletters}{OT2}{wncyr}{m}{n}
\DeclareMathSymbol{\Sha}{\mathalpha}{cyrletters}{"58}	% Sha
\newcommand{\ce}{\colonequals}
\renewcommand{\i}{^{-1}}
\newcommand{\isomto}{\overset{\sim}{\longrightarrow}}

\newcommand{\csub[1]}{\centering\subsection{\text{#1} \nopunct} \hfill \justify}

\providecommand{\up}[1]{{\upshape(}#1{\upshape)}}
\providecommand{\uref}[1]{{\upshape\ref{#1}}}
\providecommand{\uS}{{\upshape\S}}

\renewcommand{\b}{\textbf}
\providecommand{\ucolon}{{\upshape:} }
\providecommand{\uscolon}{{\upshape;} }

\newcommand{\brems}{\begin{rems} \hfill \begin{enumerate}[label=\b{\thenumberingbase.},ref=\thenumberingbase]}
\newcommand{\remi}{\addtocounter{numberingbase}{1} \item}
\newcommand{\erems}{\end{enumerate} \end{rems}}
\newcommand{\begs}{\begin{egs} \hfill \begin{enumerate}[label=\b{\thenumberingbase.},ref=\thenumberingbase]}

\newcommand{\eegs}{\end{enumerate} \end{egs}}
\newcommand{\m}{\item}
\newcommand{\bsm}{\begin{smallmatrix}}
\newcommand{\esm}{\end{smallmatrix}}
\newcommand{\blem}{\begin{lemma}}
\newcommand{\elem}{\end{lemma}}

\newcommand{\bconj}{\begin{conj}}
\newcommand{\econj}{\end{conj}}
\newcommand{\bprob}{\begin{Problem}}
\newcommand{\eprob}{\end{Problem}}

\newcommand{\bq}{\begin{Q}}
\newcommand{\eq}{\end{Q}}
\newcommand{\benum}{\begin{enumerate}[label={{\upshape(\alph*)}}]}
\newcommand{\benuma}{\begin{enumerate}[label={{\upshape(\arabic*)}}]}
\newcommand{\benumb}{\begin{enumerate}[label={{\upshape\b{\arabic*.}}}]}
\newcommand{\benumr}{\begin{enumerate}[label={{\upshape(\roman*)}}]}
\newcommand{\eenum}{\end{enumerate}}
\newcommand{\bitem}{\begin{itemize}}
\newcommand{\eitem}{\end{itemize}}
\newcommand{\bc}{\begin{comment}}
\newcommand{\ec}{\end{comment}}
\newcommand{\bd}{\begin{defn}}
\newcommand{\ed}{\end{defn}}

\newcommand{\beg}{\begin{eg}}
\newcommand{\eeg}{\end{eg}}

\newcommand{\bcl}{\begin{claim}}
\newcommand{\ecl}{\end{claim}}

\newcommand{\x}{\text}

\newcommand{\q}{\quad}
\providecommand{\qxq}[1]{\quad\text{#1}\quad}

\newcommand{\ba}{\begin{aligned}}
\newcommand{\ea}{\end{aligned}}
\newcommand{\be}{\begin{equation}}
\newcommand{\ee}{\end{equation}}
\newcommand{\bpf}{\begin{proof}}
\newcommand{\epf}{\end{proof}}
\newcommand{\bthm}{\begin{thm}}
\newcommand{\ethm}{\end{thm}}
\newcommand{\bprop}{\begin{prop}}
\newcommand{\eprop}{\end{prop}}
\newcommand{\bcor}{\begin{cor}}
\newcommand{\ecor}{\end{cor}}
\newcommand{\brem}{\begin{rem}}
\newcommand{\erem}{\end{rem}}

\usepackage{aliascnt}
\newaliascnt{numberingbase}{subsection}
\numberwithin{equation}{subsection}

\newtheoremstyle{thms}{0.5em}{0.5em}{\itshape}{}{\bfseries}{.}{ }{}
\theoremstyle{thms}
\newtheorem{conj}[numberingbase]{Conjecture}

\newtheorem{cor}[numberingbase]{Corollary}

\newtheorem{lemma}[numberingbase]{Lemma}

\newtheorem{prop}[numberingbase]{Proposition}

\newtheorem{proposition}[numberingbase]{Proposition}
\newtheorem{Q}[numberingbase]{Question}
\newtheorem{thm}[numberingbase]{Theorem}
\newtheorem{theorem}[numberingbase]{Theorem}

\newtheoremstyle{claims}{0.5em}{0.5em}{}{}{\itshape}{.}{ }{}
\theoremstyle{claims}
\newtheorem{claim}[equation]{Claim}

\newtheoremstyle{defs}{0.5em}{0.5em}{}{}{\bfseries}{.}{ }{}
\theoremstyle{defs}
\newtheorem{defn}[numberingbase]{Definition}

\newtheorem{caution}[numberingbase]{Caution}
\newtheorem{eg}[numberingbase]{Example}
\newtheorem*{egs}{Examples}
\newtheorem{rem}[numberingbase]{Remark}

\newtheorem{remark}[numberingbase]{Remark}
\newtheorem*{rems}{Remarks}

\Crefname{claim}{Claim}{Claims}
\Crefname{bclaim}{Claim}{Claims}
\Crefname{sublemma}{Lemma}{Lemmas}
\Crefname{conj}{Conjecture}{Conjectures}
\Crefname{cor}{Corollary}{Corollaries}
\Crefname{defn}{Definition}{Definitions}
\Crefname{eg}{Example}{Examples}
\Crefname{prop}{Proposition}{Propositions} 
\Crefname{Q}{Question}{Questions}
\Crefname{rem}{Remark}{Remarks}
\Crefname{thm}{Theorem}{Theorems}
\Crefname{Theorem}{Theorem}{Theorems}
\Crefname{variant}{Variant}{Variants}
\Crefname{caution}{Caution}{Cautions}

\theoremstyle{thms}
\newtheorem{thm-tweak}[subsection]{Theorem}
\Crefname{thm-tweak}{Theorem}{Theorems}
\newtheorem{lemma-tweak}[subsection]{Lemma}
\Crefname{lemma-tweak}{Lemma}{Lemmas}
\newtheorem{cor-tweak}[subsection]{Corollary}
\Crefname{cor-tweak}{Corollary}{Corollaries}
\newtheorem{prop-tweak}[subsection]{Proposition}
\Crefname{prop-tweak}{Proposition}{Propositions} 
\newtheorem{conj-tweak}[subsection]{Conjecture}
\Crefname{conj-tweak}{Conjecture}{Conjectures} 
\newtheorem{q-tweak}[subsection]{Question}
\Crefname{q-tweak}{Question}{Questions} 

\theoremstyle{defs}
\newtheorem{defn-tweak}[subsection]{Definition}
\Crefname{defn-tweak}{Definition}{Definitions}
\newtheorem{eg-tweak}[subsection]{Example}
\Crefname{eg-tweak}{Example}{Examples}
\newtheorem*{rems-tweak}{Remarks}
\newtheorem{rem-tweak}[subsection]{Remark}
\Crefname{rem-tweak}{Remark}{Remarks}

\newtheoremstyle{subsection-tweak}
   {2pt}
   {3pt}%
   {}
   {}%
   {\bfseries}
   {}%
   {.5em}
   {\thmnumber{\@{#1}{}\@{#2}.}%
    \thmnote{~{\bfseries#3.}}}    
    
\theoremstyle{subsection-tweak}
\newtheorem{pp}[numberingbase]{}
\newcommand{\bpp}{\begin{pp}}
\newcommand{\epp}{\end{pp}}

\theoremstyle{subsection-tweak}
\newtheorem{pp-tweak}[subsection]{}

\makeatletter
\def\@tocline#1#2#3#4#5#6#7{
    \begingroup 
    \@ifempty{#4}{}{}

    \parindent\z@ \leftskip#3\relax \advance\leftskip\@tempdima\relax
    #5\hskip-\@tempdima
      \ifcase #1
       \or\or \hskip 2em \or \hskip 1em \else \hskip 3em \fi%
      #6\nobreak\relax
    \dotfill\hbox to\@pnumwidth{\@tocpagenum{#7}}\par
    \nobreak
    \endgroup
 }
 \def\l@section{\@tocline{1}{0pt}{1pc}{}{}}

\renewcommand{\tocsection}[3]{%
  \indentlabel{\@ifnotempty{#2}{\makebox[1.3em][l]{%
    \ignorespaces#1 \bfseries{#2}.\hfill}}}\bfseries{#3}
    \vspace{-5pt}}

\renewcommand{\tocsubsection}[3]{%
  \indentlabel{\@ifnotempty{#2}{\hspace*{-0.5em}\makebox[2.1em][l]{%
    \ignorespaces#1#2.\hfill}}}#3
    \vspace{-5pt}}
   
\makeatother 

\makeatletter % This is for proper formatting of appendices in the table of contents
\newcommand\appendix@section[1]{%
  \refstepcounter{section}%
  \orig@section*{Appendix \@Alph\c@section. #1}%
}
\let\orig@section\section
\g@addto@macro\appendix{\let\section\appendix@section}
\makeatother

\makeatletter
\AtBeginDocument{%
  \ifdefined\Hy@deferred@thm@head % newer hyperref wrapped the theorem head
    \patchcmd\Hy@deferred@thm@head{\addvspace{-\parskip}}{}{}%
      {\typeout{Couldn't patch \string\Hy@deferred@thm@head}}%
  \fi}
\@ifpackagelater{microtype}{2022/01/01}{\microtypesetup{nopatch=all}}{}
\makeatother
\UseMicrotypeSet[expansion]{basictext}

\makeatletter
\@namedef{subjclassname@2020}{%
  \textup{2020} Mathematics Subject Classification}
\makeatother

\author{K\k{e}stutis \v{C}esnavi\v{c}ius}
\address{Sorbonne Universit\'e, Universit\'e Paris Cit\'e, CNRS, IMJ-PRG, F-75005 Paris, France}
\email{kestutis@imj-prg.fr}

\date{\today}

\usepackage{tikz}
\usepackage{tikz-cd}

\makeatletter
\newcommand\xleftrightarrow[2][]{%
  \ext@arrow 9999{\longleftrightarrowfill@}{#1}{#2}}
\newcommand\longleftrightarrowfill@{%
  \arrowfill@\leftarrow\relbar\rightarrow}
\makeatother

\makeatletter
\newcommand{\primeitem}{%
  \item[\upshape{(\roman{enumi}$'$)}]%
  \csname phantomsection\endcsname % hyperref anchor; does nothing without hyperref
  \protected@edef\@currentlabel{(\roman{enumi}$'$)}%
}
\makeatother

\makeatletter
\newcommand{\primeprimeitem}{%
  \item[\upshape{(\roman{enumi}$''$)}]%
  \csname phantomsection\endcsname % hyperref anchor; does nothing without hyperref
  \protected@edef\@currentlabel{(\roman{enumi}$''$)}%
}
\makeatother

\begin{document}

%\subjclass[2020]{Primary ; Secondary .}
%\keywords{}

\title{Kneser--Tits over fields of cohomological dimension $\le 2$ \\ \vspace{5pt} \scriptsize{with an appendix by Alexis Bouthier and Federico Scavia}}

\maketitle

\begin{abstract} 
We settle the Kneser--Tits problem over every field $k$ of (separable) cohomological dimension $\le 2$, thus resolving a conjecture of Gille. More precisely, for such a $k$ and a semisimple, simply connected $k$-group $G$ whose indecomposable factors over $k$ are isotropic, we show that $G(k)$ is generated by “elementary matrices” (the $k$-points of the unipotent radicals of the parabolic $k$-subgroups of $G$), in other words, that $K_1^G(k) \cong 1$. In fact, we show this vanishing already when $k$ has separable $p$-dimension $\le 2$ in the sense of Gille with respect to a subset of the torsion primes $p$ of $G$, and we reinterpret the separable $p$-dimension itself via syntomic cohomology. To approach those $G$ that escape the net of previously known cases, for an inclusion of reductive $k$-groups $M \subset L$, we develop a new method based on the Rost invariant for transferring $R$-equivalence on $L/M$ to $R$-equivalence on $L$. To find homogeneous spaces to which our method fruitfully applies, we use the theories of prehomogeneous vector spaces and of special symplectic representations: concretely, for a $32$-dimensional half-spin representation $V$, we argue that the vectors in its unique open orbit $\mathcal O \subset V$ admit unique up to interchange decompositions into sums of two pure spinors, and show that various double covers of $\mathcal O$ built from such decompositions are rational varieties over any base field. %To illustrate our result, we give a precise description of the derived subgroup of the group of $k$-points of suitably isotropic reductive $k$-groups.
\end{abstract}

\hypersetup{
    linktoc=page,     %set to all if you want both sections and subsections linked
}
\renewcommand*\contentsname{}
\q\\
\tableofcontents

\section{Introduction}

For a connected reductive group $G$ over a field $k$, its elementary subgroup $G(k)^+ \lhd G(k)$ is generated by the $k$-points of the unipotent radicals of the parabolic $k$-subgroups of $G$. Over infinite $k$, the unipotent radicals of a fixed pair of suitably transversal opposite $k$-parabolics suffice already, for instance, $\mathrm{GL}_n(k)^+$ is the group generated by the $n \times n$ elementary matrices. The quotient $K_1^G(k) \colonequals G(k)/G(k)^+$ is the unstable $K_1$-group of $G$, and the Kneser--Tits problem, surveyed in \cite{Gil09}, asks when this group is trivial. Examples of nontrivial $K_1^G(k)$ exist for general $k$, however, a conjecture of Gille \cite{Gil09}*{question 8.5} (sharpened in \cite{Gil19}*{conjecture 9.5.2}) predicts this triviality when $k$ is of dimension $\le 2$. We settle the remaining cases of this prediction.

The Kneser--Tits problem is one of several interrelated, long-standing arithmetic questions about algebraic groups: it is connected to the Serre Conjecture II, to the classification of retract rational groups, to the projective simplicity and the commutator problem for $G(k)$, see \cite{Gil09} for an overview. As an example, the groups $G(k)$ describe glueing data in various patching diagrams for torsors, and the elements of $G(k)^+$ give torsors that are easier to manage because, for instance, elementary matrices lift along surjections and can be approximated---this makes $K_1^G(-)$ important for studying $G$-torsors, much like the usual $K_1(-)$ is important for analyzing the Grothendieck group $K_0(-)$.

\bthm[\Cref{thm:M}]\label{thm:MA} For a field $k$ of separable dimension $\le 2$ and a semisimple, simply connected $k$-group $G$ whose indecomposable $k$-factors are all isotropic,
\[
G(k) = G(k)^+, \quad \text{equivalently,} \quad K_1^G(k)\cong 1.
\]
\ethm

As we will see when reviewing prior work in \S\ref{pp:prior-work}, \Cref{thm:MA} covers new cases already in characteristic $0$. In positive characteristic, the separable dimension used in its formulation is as in Gille’s weakening \cite{Gil00}*{section 1.1 (c)} of the prior definition of Kato \cite{Kat82}*{p.~220, Definition}. Concretely, our dimension assumption means that $H^i_{\mathrm{\text{ét}}}(k, \mathscr F(j)) \cong 0$ for all $i > 2$, all $i > j \ge 0$, and all torsion abelian sheaves $\mathscr F$ on $k_{\mathrm{\text{ét}}}$---we show in \Cref{sec:appendix} that this reinterpretation in terms of the syntomic cohomology of $k$ is equivalent to the more hands-on definition of Gille. Such syntomic cohomology groups are relevant because they host the values of cohomological invariants of algebraic groups, for instance, the Rost invariant is valued in $H^3_{\mathrm{\text{ét}}}(k, (\mathbb Q/\mathbb Z)(2))$.

To simplify, we have understated what \Cref{thm:M} actually shows: in fact, $K_1^G(k) \cong 1$ already when $k$ has separable $p$-dimension $\mathrm{dim}_p^{\mathrm{sep}}(k)\le 2$ for a short list of primes $p$ that is a subset of the torsion primes of $G$, and away from type $A_n$ only $p = 2$ is relevant. This resolves the isotropic case of \cite{Gil19}*{conjecture 9.5.2}, even in a sharper form than there because the isotropicity and general rationality results force many torsion primes of $G$ to not matter.

\bpp[Prior work and the gist of the problem] \label{pp:prior-work} \Cref{thm:MA} was known in the following cases, the vast majority of them for perfect $k$ and without considering the finer version where separable $p$-dimension $\le 2$ is only assumed for torsion primes $p$.
\benum
\m
\label{pw:A} The case when $k$ is perfect and $G$ has no absolutely indecomposable constituent of Tits index $E_{7,1}^{66}$ was settled by Gille \cite{Gil09}*{théorème 8.6} (with \cite{Gil10}*{end of Section 7}), building on prior work of Yanchevskij \cite{Yan74}, \cite{Yan75}, Tits \cite{Tit78}, Chernousov–Platonov \cite{CP98}, Merkurjev \cite{Mer98}, and Garibaldi \cite{Gar07}. 

\m
\label{pw:B} The case when $k$ is a global field was settled by Gille \cite{Gil09}*{théorème 1.3}, building on prior work of Tits \cite{Tit78}, Prasad--Raghunathan \cite{PR85}, \cite{Pra08}, Chernousov--Timoshenko \cite{CT99}, and Garibaldi \cite{Gar07}.

\m
\label{pw:C} The case when $k$ is a nonarchimedean local field was settled by Platonov \cite{Pla69}, \cite{Pla70}, see also \cite{Tit78}*{section 1.2}, \cite{PR85}*{Section 2}, \cite{PR94}*{Theorem 7.6}.

\m
\label{pw:D} The case when $k$ is of separable dimension $\le 1$ follows from Steinberg’s theorem \cite{Ste65}*{Theorem 1.9} (see also \cite{BS68}*{Section 8.6}).

\m
\label{pw:E} The case when $G$ is of type $A_n$ was settled in the fine form of an arbitrary $k$ and separable $p$-dimension $\le 2$ for torsion primes $p$ by Grenier-Boley in \cite{GB07}*{Theorem 1.1}, building on prior work of Yanchevskij \cite{Yan74}, \cite{Yan75} (see also \cite{NM43}, \cite{Wan50}), and Gille \cite{Gil00}.

\m
\label{pw:F} Various cases in which $k$ is either of geometric nature or discretely valued were settled in \cite{CTGP04}*{Theorem 4.5, Corollaries 4.6 and 4.11} (see also \cite{CTGP01}*{théorème 2.5}), \cite{Gil09}*{corollaire 8.10}, and \cite{GILA25}*{Section 3}.
\eenum
All in all, ignoring the somewhat thorny bad characteristic aspects, the crux of the matter for \Cref{thm:MA} is to settle “just” the case of semisimple, simply connected groups of Tits index $E_{7,1}^{66}$. The specificity of this is somewhat illusory: in fact, all the other exceptional, semisimple, simply connected, absolutely indecomposable, isotropic  $k$-groups $G$ are either “small enough” to be retract rational (or even rational) as varieties over $k$, so that, as we will recall, automatically $K_1^G(k) \cong 1$, or they are “large enough” so that the cohomological dimension assumption leaves no room for them over $k$ (and often they are both). The $E_{7, 1}^{66}$ happen to fall precisely between the two extremes (see \Cref{prop:no} and \Cref{rem:no}), unless $k$ is so small that it cannot accommodate them either, as happens in \ref{pw:B}--\ref{pw:D} and \ref{pw:F}. In view of all this, instead of trying to nuke the $E_{7,1}^{66}$’s by any means possible, our priority is to simultaneously develop broadly useful, new techniques for analyzing exceptional (and other) groups that would work over arbitrary base fields (see \Cref{caution} \ref{m:caution-2}).\footnote{\label{foot:focused}Of course, without losing sight of \Cref{thm:MA}. For instance, our method ought to adapt to show the rationality of simply connected groups of Tits index ${}^2E_{6,1}^{29}$—these are known to be stably rational by Garibaldi \cite{Gar07}*{Proposition~0.1} (see also \Cref{prop:retract}~\ref{m:R-c}~\ref{m:R-x} below)—but we did not include this (nor some other tangents) for the sake of focus.}
\epp

\bpp[$R$-equivalence and dévissage]
The unipotent radicals that define $G(k)^+$ are affine spaces as $k$-schemes, so the triviality of $K_1^G(k)$ implies that $G(k)$ is exhausted by products of those $k$-points that may be joined to the identity by affine lines $\mathbb A^1_k \subset G$. By a well-known comparison due to Gille (see \Cref{prop:W-R}), we lose no generality by replacing $\mathbb A^1_k$ by its semilocalization at $0$ and $1$: 
\[
K_1^G(k) \cong G(k)/R,
\]
where the right hand side is the group of $R$-equivalence classes in the sense of Manin. Therefore, we need to build many rational curves in $G$ with the help of the isotropy assumption and the cohomological dimension of $k$ (for $E_{7,1}^{66}$ merely $\mathrm{dim}^{\mathrm{sep}}_2(k)\le 2$ will suffice). This is simplest if $G$ itself is rational (or at least retract rational), so known rationality results dispose of quite a few cases, see the general \Cref{thm:RR} about $R$-equivalence for retract rational varieties. 

In general, we develop a dévissage technique for lifting the $R$-triviality of a homogeneous space to the $R$-triviality of a group $L$ that acts transitively on it, granted that the (usually much lower-dimensional and hence simpler) stabilizer group $M$ is $R$-trivial as well, see \Cref{prop:nonabelian-R}. This is a bit similar to the control of $R$-equivalence in short exact sequences \cite{BKG04}*{Appendix, Theorem~1}, or to \cite{GS25}*{Lemma 2.6 (2)} that considered $R$-equivalence on homogeneous spaces (such as $L/M$), however, our criterion applies more broadly. This is because we only require the triviality of those $M$-torsors over $k(t)$ that become trivial when induced to $L$-torsors, in fact, even with $k(t)$ replaced by discretely valued fields such as $k(\!(t)\!)$ thanks to a trick of Gabber that circumvents the noncanonicity of the field of representatives in the Cohen structure theorem over imperfect fields.
\epp

\bpp[Using the Rost invariant to trivialize $M$-torsors] 
The robustness of our dévissage has a price: to apply it, we need to argue the triviality of some very specific $M$-torsors over $k(\!(t)\!)$. For this, we use the Rost invariant and the cohomological dimension assumption on $k$—this is the only place in the overall argument for $E_{7,1}^{66}$ where cohomological dimension plays a role. Assuming that we are working with groups for which the Rost invariant is defined, the point is that the Rost multiplier of the inclusion $M \hookrightarrow L$ is often $1$, so the triviality of the induced $L$-torsor at least tells us that the Rost invariant of the initial $M$-torsor over $k(\!(t)\!)$ vanishes. We can then draw on the vast literature about the triviality of the kernel of the Rost invariant, for instance, away from bad characteristics we could trivialize the $M$-torsors over $k(\!(t)\!)$ relevant for our argument by Preeti's \cite{Pre13}*{Theorem~1.7}.\footnote{\label{foot:Preeti} The statement of \emph{loc.~cit.}~should be restricted to absolutely indecomposable groups, otherwise the additivity of the Rost invariant with respect to products gives counterexamples.} 

We wish to include all characteristics and, as mentioned in \Cref{caution}~\ref{m:caution-1}, we have no hope of controlling $\mathrm{dim}^{\mathrm{sep}}_p(k(\!(t)\!))$ anyway, so under the assumptions that matter for us ($M$ of type $A_n$), we devise a direct argument for the triviality of the kernel of the Rost invariant in \Cref{thm:Rost-A2}. This argument is based on results of Berhuy–Frings–Tignol \cite{BFT07} about the classification of unitary Hermitian forms when $\mathrm{dim}_2^{\mathrm{sep}}(k) \le 2$. We came to it by analyzing \cite{PS24}, where, in a somewhat different way, Parimala–Suresh have shown the triviality of the kernel of the Rost invariant in type $A_n$ over (for example) $k(\!(t)\!)$ when $k$ is either a local field or a global field without real places.
\epp

\bpp[Building the relevant homogeneous space]\label{pp:build}
To apply our homogeneous space dévissage, given a semisimple, simply connected $k$-group $G$ of isotropic $k$-rank $1$, we let $L$ be a Levi of a proper parabolic $k$-subgroup of $G$. Thus, $L$ is the intersection of two opposite proper $k$-parabolics of $G$, and we let $M$ be the (identity component) of the intersection with a third, sufficiently general $k$-parabolic. By the open cell decomposition, the sought $R$-triviality of $G$ amounts to that of $L$, and when $G$ has Tits index $E_{7,1}^{66}$ this $M$ is small enough to be $R$-trivial: it is semisimple, simply connected of type $A_5$, so it is retract rational, as reviewed in \Cref{prop:retract}~\ref{m:R-c} (in fact, this $M$ is of outer type $A_5$ and the division algebra that intervenes in its description \S\ref{pp:conv} \ref{m:conv-i} \ref{m:conv-b} has index $2$ over its center, so $M$ is even rational over $k$ by a result of Chernousov–Merkurjev \cite{CM98}*{Corollary 6.5}). 

Intuitively, we are removing the circled node from the Tits index
\[
%\begin{comment}
\begin{tikzpicture}[
  x=0.8cm,
  y=0.8cm,
  line width=0.5pt,
  font=\Large
]
  \node[anchor=east, inner sep=0pt] at (-1.05,0) {\scalebox{0.8}{$E_{7,1}^{66}$}};

  \draw (0,0) -- (5,0);
  \draw (3,0) -- (3,1);

  \foreach \x/\i in {0/7,1/6,2/5,3/4,4/3,5/1} {
    \fill (\x,0) circle[radius=2.1pt];
    \node[below=6pt] at (\x,0) {\scalebox{0.8}{}};
  }

  \fill (3,1) circle[radius=2.1pt];
  \node[right=6pt] at (3,1) {\scalebox{0.8}{}};

  \draw (5,0) circle[radius=0.2cm];
\end{tikzpicture}
%\end{comment}
\]
to get the derived group $H \colonequals [L, L]$ of type $D_6$, which comes with a distinguished branch node (to which the circled node was attached), and we are then removing that branch node to get to an $M$ of type $A_5$. To apply our dévissage, we need to understand $R$-equivalence on $X \colonequals L/M$, and in \Cref{cor:LM-ratl}~\ref{m:LM-c} we show that this homogeneous space is even a rational variety of dimension $32$ over every $k$. For this, we change the perspective and forget about $E_{7,1}^{66}$: by the theory of Tits algebras and highest weight, the distinguished branch node of $H$ corresponds to a $k$-descent $V$ of a $32$-dimensional half-spin representation of $H_{\overline{k}} \simeq \mathrm{Spin}_{12,\, \overline{k}}$, on which the action extends to $L$, and the key to the rationality of $X$ ends up hidden in very special representation-theoretic properties of $V$.
\epp

\bpp[Rationality of the homogeneous space via spin geometry] \label{pp:spin}
The representation $V$ is exceptional: it appears in the Sato–Kimura classification \cite{SK77} of prehomogeneous vector spaces, meaning that the $L$-action has an open orbit $\mathcal O \subset V$. Ignoring the fact that the classification of \emph{op.~cit.}~is over $\mathbb C$, this suggests a plausible source of the rationality of $X$: the open orbit $\mathcal O \subset V$ is a rational variety, it has a $k$-point $v$ as soon as $k$ is infinite, and $\mathrm{Stab}_L(v)^0$ is the $M$ from \S\ref{pp:build}. The difficulty is that $\mathrm{Stab}_L(v)$ is not connected, which douses the effectiveness of $R$-equivalence. Its component group is $\mathbb Z/2\mathbb Z$ and may be viewed as compensating for the apparent mismatch between the Dynkin diagram of $E_7$ having no nontrivial automorphisms yet $M$ being of \emph{outer} type $A_5$. 

Instead, we use special features of half-spinor geometry in dimension $32$ to show that the double cover $X \cong L/\mathrm{Stab}_L(v)^0$ of $\mathcal O$ is a rational variety for every base field $k$. Our approach is based on the cone $\mathscr S \subset V \setminus \{0\}$ of pure spinors defined as the preimage of the unique closed $L$-orbit $S\subset \mathbb P(V)$ (the orbit of the highest weight line). Two pure spinors are \emph{opposite} if the $H$-stabilizers of their lines are opposite as parabolic subgroups; this defines a dense open $\mathscr O \subset \mathscr S \times_k \mathscr S$. The dimension of $\mathscr S$ is $16$ and  the open $L$-orbit $\mathcal O \subset V$ consists precisely of those vectors that are sums of two opposite pure spinors, necessarily uniquely up to interchanging the summands: in \Cref{thm:dt} we show that the sum map $\mathscr O \rightarrow \mathcal O$ is an $L$-equivariant étale double cover. This is very special to half-spin dimension $32$: in general, a half-spin representation of $D_n$ has dimension $2^{n - 1}$, its cone of pure spinors has dimension $\frac{n(n - 1)}{2} + 1$, and the former is twice the latter only for $n = 6$. 

The unique decomposition into pure spinors in dimension $32$ has been observed by Hitchin over $\mathbb C$ in \cite{Hit03}*{Proposition 2}, by Slupinski–Stanton over fields of characteristic different from $2$ and $3$ in \cite{SS15}*{Theorem 4.7} slightly more generally for special symplectic representations, and by Helenius again away from characteristic $2$ and $3$ in \cite{Hel12}*{Lemma 24} in the context of Freudenthal triple systems. These works produce explicit formulas for the pure spinor summands in terms of the various invariant functions associated to $V$ and they involve division by $2$ and $3$. Our more indirect approach mixes  techniques of arithmetic geometry and representation theory and works over any scheme. 

In the end, our $X \cong L/\mathrm{Stab}_L(v)^0 \rightarrow \mathcal O$ is the twist of the double cover $\mathscr O \rightarrow \mathcal O$ by the quadratic étale $k$-algebra $k_v$ over which the opposite pure spinor summands of $v$ are defined (see \Cref{cor:LM-ratl}). The sought rationality of $X$ then comes from $S_{k_v}$ (via $\mathscr S_{k_v}$ and $\mathscr O_{k_v}$): as any flag variety of a reductive group, $S$ becomes rational over any field over which it acquires a rational point (see \Cref{lem:S-ratl}).
\epp

%The Kneser–Tits problem originated with the Tannaka–Artin problem that aims to determine the commutator subgroup of a central division algebra over a field, see \cite{Tit78}*{page 218} and the references there. \Cref{thm:M} gives the following result on the Tannaka–Artin problem for reductive groups without anisotropic constituents over a field of cohomological dimension $\le 2$. 

%\bcor \label{cor:C}
%...
%\ecor

\bpp[Notation and conventions] \label{pp:conv}
For a field $k$, we let $k^s \subset\overline{k}$ denote its separable and algebraic closures. A \emph{$k$-variety} is a separated $k$-scheme of finite type. A nonempty $k$-variety $X$ is \emph{rational} (resp.,~\emph{retract rational}) if it contains a dense open $U \subset X$ that is isomorphic to (resp.,~that is a retract of) a nonempty open $V \subset \mathbb A^n_k$ (resp.,~meaning that some $k$-morphisms $i \colon U \rightarrow V$ and $r \colon V \rightarrow U$ satisfy $r \circ i = \mathrm{id}_U$). A nonempty $k$-variety $X$ is \emph{stably rational} if $X \times_k \mathbb A^m_k$ is rational for some (equivalently, for every large) $m \ge 0$. A nonempty $k$-variety $X$ is (resp.,~\emph{separably}) \emph{unirational} if it admits a (resp.,~generically smooth) dominant map from a nonempty open $V \subset \mathbb A^n_k$.~All in all,
\[
\text{rational} \Rightarrow \text{stably rational} \Rightarrow \text{retract rational} \Rightarrow \text{unirational},
\]
and each of these notions implies that $X$ has many rational curves. For a central simple algebra $A$ over $k$, its \emph{degree} is $\mathrm{deg}_k(A)\colonequals \sqrt{\mathrm{dim}_k(A)}$, its \emph{period} is the order of the class $[A]$ in the Brauer group of $k$, and its \emph{index} is the degree of a division algebra that represents $[A]$. The period divides the index, and the two have the same prime divisors, see \cite{GS17}*{Theorem 2.8.7}. We follow the convention that a reductive group over a field is assumed to be connected. 

An algebraic group $G$ over $k$ is \emph{indecomposable} (resp., \emph{absolutely indecomposable}) if it cannot be decomposed as a product of two nontrivial $k$-subgroups (resp.,~if $G_{\overline k}$ is indecomposable). For semisimple, simply connected groups, absolute indecomposability amounts to being `absolutely almost simple,' although we find our terminology more intuitive. An algebraic group $G$ over a field $k$ is \emph{isotropic} (resp.,~\emph{anisotropic}) if it contains (resp.,~does not contain) $\mathbb G_{m,\, k}$ as a $k$-subgroup. The \emph{isotropic rank} of a semisimple group $G$ over $k$ is the largest $r\ge 0$ such that $G$ contains $\mathbb G_{m,\, k}^r$ as a $k$-subgroup. We use the geometric classification of semisimple, simply connected, absolutely indecomposable $k$-groups $G$ via Dynkin diagrams \cite{SGA3IIInew}*{exposé XXV, théorème~1.1; exposé~XXII, définition 4.3.3; exposé XXIV, section 5.3, proposition 5.10; exposé XXI, corollaire~7.4.4, remarque 7.4.6} and the following descriptions in the cases $A_n$ and $D_n$ (see \S\ref{pp:build} for their relevance).
\benumr
\m \label{m:conv-i} If $G$ is of type $A_n$, then $G \simeq \mathrm{SU}(A, \sigma)$ for a quadratic étale $k$-algebra $k^\prime$, a central simple algebra $A$ of degree $n + 1$ over $k^\prime$, and an involution $\sigma$ on $A$ that is $k^\prime/k$-unitary in the sense that it restricts to the nontrivial $k$-automorphism of $k^\prime$, where $\mathrm{SU}(A, \sigma)$ is the $k$-group that parametrizes the sections $g$ of $A$ that satisfy $\sigma(g)g = 1$ (equivalently, that satisfy $g \sigma(g) = 1$) and have reduced norm $1$, see \cite{KMRT98}*{Theorem 26.9 (with page 346 for the definition)}. We distinguish the two subcases: either $k^\prime \simeq k \times k$ or $k^\prime$ is a field.
\benum
\m \label{m:conv-a} (Inner type). If $k^\prime \simeq k \times k$, then $G \simeq \mathrm{SL}_r(D)$ for a central division algebra $D$ over $k$ and an $r \ge 1$, see \cite{KMRT98}*{Corollary 1.2}.

\m \label{m:conv-b} (Outer type when $n \ge 2$). If $k^\prime$ is a field, then $G \simeq \mathrm{SU}_r(D, h)$ for a central division algebra $D$ over $k^\prime$ with a $k$-involution $\sigma$ that is nontrivial on $k^\prime$, an $r \ge 1$, and a nondegenerate $D$-valued $\sigma$-Hermitian form $h$ on the right $D$-module $D^r$, that is, $G$ is the $k$-group of the $D$-linear endomorphisms of $D^r$ of reduced norm $1$  preserving $h$, see \cite{KMRT98}*{Corollary~1.2, Theorem 3.1 (2), and Theorem 4.2 (2)} (with \cite{KMRT98}*{Proposition 4.1} to recast via $h$). 
\eenum
The index of the central simple algebra $A$ is equal to that of $D$ and divides $n + 1$.

\m \label{m:conv-ii} 
If $G$ is of nontrialitarian type $D_n$, then $G \simeq \mathrm{Spin}(A, \sigma, f)$ for a central simple algebra $A$ over $k$ of degree $2n$, a $k$-involution $\sigma$ on $A$, and, letting $\mathrm{Trd}_A \colon A \rightarrow k$ denote the reduced trace of $A$, a $k$-linear map $f \colon A^{\sigma = \mathrm{id}} \rightarrow k$ with $f(a + \sigma(a)) = \mathrm{Trd}_A(a)$ such that $\mathrm{dim}_k(A^{\sigma = \mathrm{id}}) = n(2n + 1)$ and $\mathrm{Trd}_A(A^{\sigma = -\mathrm{id}}) = 0$ (which makes the condition on $f$ well defined), see \cite{KMRT98}*{Theorem~26.15 with Definition 5.4 (and Definition 13.30 for the spin group)} in the case $n \ge 5$ and \cite{CF15}*{Propositions 4.4.0.34 and 4.5.2.7, théorème 8.2.0.59} (with \cite{Rue24}*{equation before Lemma 3.2 (with p. 69)}) in general. 
If $\mathrm{char}(k) \neq 2$, then $f$ is determined (and omitted) and $\mathrm{Trd}_A(A^{\sigma = -\mathrm{id}}) = 0$ is automatic. In the case when $A$ splits, by \cite{KMRT98}*{Proposition 5.11 and Definition 13.30}, the data $(\sigma, f)$ amount to a nondegenerate quadratic form $q$ on $k^{2n}$, and $\mathrm{Spin}(A, \sigma, f) \simeq \mathrm{Spin}(q)$. In particular,  $G_{\overline{k}} \simeq \mathrm{Spin}_{2n,\, \overline{k}}$. 

The $k$-involution $\sigma$ forces the period of $A$ to be $\le 2$, see \cite{KMRT98}*{Theorem 3.1 (1)}.
\eenum
For a field $k$, we let $H^i_{\mathrm{\text{\'et}}}(k, (\mathbb Q/\mathbb Z)(j))$ denote the indicated syntomic cohomology group defined using \cite{BM23}*{Notation 1.1}. Concretely, for a prime $p$, its $p$-primary part agrees with the usual definition as étale cohomology of Tate twists for $p \neq \mathrm{char}(k)$ and as étale cohomology of suitable logarithmic de Rham–Witt sheaves for $p = \mathrm{char}(k)$, see \cite{BM23}*{Examples 1.2 and 1.3}. More generally, for $j \ge 0$ and a $p$-primary torsion abelian sheaf $\mathscr F$ on $k_{\mathrm{\text{ét}}}$, we set $\mathscr F(j) \colonequals \mathscr F \otimes^\mathbb L_{\mathbb Z} \mathbb Z_p(j)$ with the complex $\mathbb Z_p(j)$ on $k_{\mathrm{\text{ét}}}$ as in \emph{op.~cit.},~and we define $\mathscr F(j)$ via the primary decomposition  for all torsion abelian sheaves $\mathscr F$ on $k_{\mathrm{\text{ét}}}$. By \emph{loc.~cit.},~$(\mathbb Z/p\mathbb Z)(j) \cong \mu_p^{\otimes j}$ when $p \neq \mathrm{char}(k)$, and $(\mathbb Z/p\mathbb Z)(j) \cong \underline{\Omega}^j_{k,\, \mathrm{log}}[-j]$ when $p = \mathrm{char}(k)$, and, by \cite{BM23}*{Example 1.5}, we have $H^i_{\mathrm{\text{\'et}}}(k, (\mathbb Q/\mathbb Z)(1)) \cong \varinjlim_{n \ge 1} H^i_{\mathrm{fppf}}(k, \mu_n)$. Throughout we identify the fppf and étale cohomologies whenever the coefficient group is smooth. 

\epp

{\bf Acknowledgements.}  I thank Alexis Bouthier for exchanges in November 2025 about the problem addressed in this article, especially, for drawing my attention to the conjecture of Gille on the Kneser–Tits problem. I thank Federico Scavia for showing me work on a different approach at the time when I was preparing the final draft of this article. I thank Nicolas Grenier-Boley for sending me the published version of \cite{GB07}. I thank Alexis Bouthier and Federico Scavia for helpful conversations and for providing Appendix \ref{app:B} that presents a different argument for the $E_{7,1}^{66}$ case of \Cref{thm:M} away from imperfect fields of characteristic $2$.  This project was supported by the Simons Collaboration on Perfection in Algebra, Geometry, and Topology, award ID~MP-SCMPS-00001529-10.

No text in this article originated from, was written or modified in any way by any AI. At various stages of the project I used commercially available versions of large language models (LLMs) for searching the literature and testing my approaches (most of which, as usual, led to nowhere and are not discussed here); this has certainly sped up the process. My goal was to settle the remaining cases of \Cref{thm:MA}, especially, with the constraint of finding a reasonably conceptual method. 

As for the actual contents, AI found the Slupinski–Stanton result \cite{SS15}*{Theorem 4.7, Example~4.8} that supplies unique decompositions into sums of pure spinors. This turned out to be a key missing piece for one of my approaches: modulo issues in small characteristics, it proved the rationality of a relevant homogeneous space. I later evolved this proof into its simpler, equation-free form presented in \S\ref{sec:half-spin}, which made it work in every characteristic at the cost of it now seeming a bit distant from the more formula-oriented \cite{SS15}, and I also found the Helenius reference \cite{Hel12}*{Lemma 24}, which is even closer to my setting. Nevertheless, finding \cite{SS15} was a turning point for the genesis of this article, and it would have been difficult to do without AI-assisted search.

\section{Reducing to $R$-equivalence}

For arguing \Cref{thm:MA}, it is more practical to work with $R$-equivalence. We review this in \Cref{prop:W-R}, show that $R$-equivalence is invariant under specialization with respect to Henselian discrete valuation rings in \Cref{RS} and is trivial for retract rational varieties in \Cref{thm:RR}, and review some classes of retract rational, semisimple, simply connected, isotropic groups in~\Cref{prop:retract}.

\bpp[The group $K_1^G(k)$] \label{pp:K1Gk}
For a reductive group $G$ over a field $k$, we let $G(k)^+ \lhd G(k)$ be the (necessarily normal) subgroup generated by the $k$-points of the unipotent radicals $\mathscr R_{\mathrm{u}}(P)$ of the parabolic $k$-subgroups $P \subset G$. The \emph{unstable $K_1$-group} of $G$ is the quotient
\[
K_1^G(k) \colonequals G(k)/G(k)^+,
\]
also called the \emph{Whitehead group} of $G$ and denoted by $W(k, G)$. The $K$-theoretic notation $K_1^G(k)$ is that of \cite{Sta14}*{Definition 3.1} (with a precursor in \cite{Jar83}*{Section 3}) and is consistent with the case of $\mathrm{GL}_n$, in which $K_1(k) \cong \varinjlim_{n \ge 1} K_1^{\mathrm{GL}_n}(k)$. Both $G(k)^+$ and $K_1^G(k)$ are functorial in $k$ and $G$, and commute with products and restrictions of scalars of groups, see \cite{torsors-algebraic}*{Section 7.1.3}. 

Slightly different definitions for $G(k)^+$ have appeared in the literature (and have sometimes led to small inaccuracies!), for instance, instead of the $k$-points of the unipotent radicals $\mathscr R_{\mathrm u}(P)$, one may prefer the $k$-points of all the $k$-subgroups $\mathbb G_{a,\, k} \subset G$. Any such $\mathbb G_{a,\, k}$ lies in some $\mathscr R_{\mathrm u}(P)$ (see \cite{CGP15}*{Theorem C.3.8} or \cite{torsors-algebraic}*{Corollary 6.2.5}), so the resulting variant of $G(k)^+$ is \emph{a priori} smaller. The two definitions agree if $k \not\cong \mathbb F_2$ or if $G$ does not contain any nonsplit semisimple $k$-subgroup of type $A_2$, see \cite{BT73}*{Proposition 6.2}. Still, they are genuinely different when $k \cong \mathbb F_2$ and $G$ is semisimple, simply connected, absolutely indecomposable, and nonsplit of type $A_2$.
\epp

\bpp[$R$-equivalence] \label{pp:R-equivalence}
For a (contravariant) functor $X$ on the category of affine schemes over a fixed field $k$, we say that $x, x^\prime \in X(k)$ are \emph{directly $R$-equivalent} if there is a $\chi \in X(\mathscr O_{\mathbb A^1_k,\, 0,\, 1})$ with $\chi|_0 = x$ and $\chi|_1 = x'$. We say that $\chi$ \emph{realizes direct $R$-equivalence} between $x, x^\prime \in X(k)$. The equivalence relation $\sim_R$ on $X(k)$ that this generates is the \emph{$R$-equivalence} (introduced in \cite{Man86}*{Definition 14.1}), and we say that $x, x^\prime \in X(k)$ are \emph{$R$-equivalent} if they lie in the same $R$-equivalence class, in other words, if they are related by a finite chain of direct $R$-equivalences. We write $X(k)/R$ for the set of $R$-equivalence classes in $X(k)$. Both direct $R$-equivalence and $R$-equivalence respect products and also restrictions of scalars along any finite extension of fields $k^\prime/k$ (use that $\mathscr O_{\mathbb A^1_k,\, 0,\, 1} \otimes_k k^\prime \cong \mathscr O_{\mathbb A^1_{k^\prime},\, 0,\, 1}$). 

For a (contravariant) group-valued functor $G$ on the category of affine $k$-schemes, $x, y \in G(k)$ are directly $R$-equivalent if and only if so are $xy^{-1}, 1 \in G(k)$. Thus, by multiplying $\mathscr O_{\mathbb A^1_k,\, 0,\, 1}$-points that realize direct $R$-equivalence, we find that $R$-equivalence simply amounts to direct $R$-equivalence in $G(k)$. Similarly, the elements that are $R$-equivalent to $1$ form a normal subgroup $RG(k) \subset G(k)$, the quotient by which is the group of $R$-equivalence classes $G(k)/R$. %For unirational linear algebraic groups over an infinite field, this quotient is a birational invariant, see \cite{CTS77}*{page 197, corollaire}.
\epp

\bprop[Gille] \label{prop:W-R} For a field $k$ and a semisimple, simply connected $k$-group $G$ whose indecomposable $k$-factors are isotropic, we have $K_1^G(k) \xrightarrow{\sim} G(k)/R$.
\eprop

\bpf
The formation of $K_1^G(k)$ and $G(k)/R$ commutes with products and restrictions of scalars, see \S\ref{pp:K1Gk} and \S\ref{pp:R-equivalence}. Thus, by the canonical decomposition of our semisimple, simply connected $k$-group 
\be \label{eqn:can-dec}
G \cong \prod \mathrm{Res}_{k'/k}(G')
\ee
with finite separable $k^\prime/k$ and semisimple, simply connected, absolutely indecomposable $k^\prime$-groups $G^\prime$ (see \cite{SGA3IIInew}*{exposé XXIV, proposition 5.10 (i) (with section 5.3)}) and by the fact that $\mathrm{Res}_{k'/k}(G')$ is isotropic if and only if so is $G^\prime$ (see \cite{torsors-regular}*{last paragraph of Section 1.3.4}), we may assume that $G$ is absolutely indecomposable. It then remains to apply \cite{Gil09}*{théorème 7.2}.
\epf

By \cite{Gil09}*{corollaire 7.3}, for semisimple, simply connected groups $G$ with isotropic indecomposable factors, $K_1^G(-)$ is invariant under specialization along Henselian discrete valuation rings. For $R$-equivalence, as we now show, the same invariance holds for arbitrary reductive groups. This was already known in equal characteristic by \cite{Gil04}*{théorème 0.2 c)} and \cite{GS25}*{Theorem 8.14}.

\bthm \label{RS}
For a reductive group scheme $G$ over a Henselian discrete valuation ring $\mathcal O$ with fraction field $K$ and residue field $k$, the specialization map induces an isomorphism
\[
G(K)/R \xrightarrow{\sim} G(k)/R.
\]
\ethm

The specialization map for $R$-equivalence has been constructed in \cite{Gil04}*{théorème 0.2} and \cite{CTHHKPS22}*{Theorem A.10} and is characterized by being compatible with the map $G(\mathcal O) \twoheadrightarrow G(k)$, see \cite{CTHHKPS22}*{Corollary A.9 (b)} (or \cite{GS25}*{Lemma 8.9}).

\bpf
Hensel lifting gives the surjectivity, so we focus on the injectivity. By \cite{CTHHKPS22}*{Corollary~A.9~(b)} (or \cite{GS25}*{Lemma 8.9}), $G(\mathcal O)$ exhausts $G(K)/R$, so it suffices to show that for a $g \in G(\mathcal O)$, if $g_k$ is directly $R$-equivalent to $1 \in G(k)$, then $g_K$ is directly $R$-equivalent to $1 \in G(K)$. For this, let $G \subset X$ be the proper, smooth, $\mathcal O$-fiberwise dense compactification built in \cite{Nat26}*{Theorem~1.2} (see also \cite{SGA3II}*{exposé X, théorème 5.16}). A direct $R$-equivalence between $g_K$~and $1$ in $X(K)$, equivalently, in $G(K)$, is supplied by \cite{Kol04}*{Theorem 2 (1) and its proof in paragraph~17}.\footnote{A proper, smooth variety $Y$ over a field $\ell$ with $Y_{\ell^\prime}$ separably unirational for a separable extension $\ell^\prime/\ell$ is separably rationally connected in the sense of \cite{Kol04}*{Definition 1}. To see this, by a limit and spreading out argument, we may assume that the separable extension $\ell^\prime/\ell$ is finite, and then let $U \subset \mathbb A^n_{\ell^\prime}$  be a dense open with a dominant, generically smooth $\ell^\prime$-morphism $U \rightarrow Y_{\ell^\prime}$. Since $Y$ is proper, the entire $\mathbb P^1$ that passes through two generic points of $U$ in $\bP^n_{\ell^\prime}$ maps to $Y$. In particular, there are a dense open $V \subset U \times_{\ell^\prime} U$ and an $\ell^\prime$-morphism $\mathbb P^1_\ell \times_\ell V \rightarrow Y$ such that the evaluation at two fixed $\ell$-points of $\mathbb P^1_\ell$ induces a dominant, generically smooth $\ell$-morphism $V \rightarrow Y \times_\ell Y$ (the restriction of the composition $U \times_{\ell^\prime} U \rightarrow Y_{\ell^\prime} \times_{\ell^\prime} Y_{\ell^\prime} \rightarrow Y \times_\ell Y$). }
\epf

To apply \Cref{prop:W-R}, it helps to know that any two $k$-points are $R$-equivalent for retract rational $G$. This is not specific to groups thanks to the following result, which we could not locate in the literature in its natural generality, compare with \cite{CT19}*{théorème 3.7} and \cite{KS15}*{Corollary~8.6.3}. Our argument is inspired by the interpretation of $R$-equivalence in terms of weak approximation at two $k$-points of $\mathbb P^1_k$ (see \cite{Gil15a}*{section 8}), and we came to it by considering the homogeneous~space~case. %(the case we will use in arguing \Cref{prop:LkRt}).

\bthm \label{thm:RR}
For a smooth, \emph{retract rational} variety $X$ over a field $k$, all $x, x^\prime \in X(k)$ are directly $R$-equivalent, in fact, for dense opens $U \subset X$ and $V \subset \mathbb A^n_k$ and $k$-maps $i \colon U \rightarrow V$ and $r \colon V \rightarrow U$ with $r \circ i = \mathrm{id}_U$ there is a $\chi \in X(\mathscr O_{\mathbb A^1_k,\, 0,\, 1})$ that generically is in $U$ with $\chi|_0 = x$ and $\chi|_1 = x'$. Thus,~either 
\[
X(k) = \emptyset\text{ \up{with $k$ finite}} \quad \text{or} \quad X(k)/R = \{*\}.
\]
\ethm

\bpf
The final claim follows from the rest because if $k$ is infinite, then $V(k) \neq \emptyset$, and so also $X(k) \neq \emptyset$. Moreover, the density of $U$ and $V$ implies that $X$ is geometrically integral. For the rest, we equip $X(k(\!(t)\!))$ and $X(k(\!(t - 1)\!))$ with their valuation topologies (see \cite{GGMB14}*{beginning of Section 3.1}) and note that, by \cite{GGMB14}*{Proposition 3.3.4 (3)}, both
\be\label{eqn:O}
\{\widetilde x \in X(k[\![t]\!]) :  \widetilde x|_{t = 0} = x\} \subset X(k(\!(t)\!)) \ \ \ \text{and}\ \ \ \{\widetilde x' \in X(k[\![t - 1]\!]) :  \widetilde x'|_{t = 1} = x'\} \subset X(k(\!(t - 1)\!))
\ee
are nonempty opens. Consequently, it suffices to show that the pullback map
\be\label{eqn:D}
U(k(t)) \rightarrow X(k(\!(t)\!)) \times X(k(\!(t - 1)\!))
\ee
has dense image: any $\chi \in U(k(t))$ whose pullback is in the product of the opens \eqref{eqn:O} will extend to the sought $\chi \in X(\mathscr O_{\mathbb A^1_k,\, 0,\, 1})$ with $\chi|_0 = x$ and $\chi|_1 = x'$. On the other hand, $U(k(\!(t)\!)) \subset X(k(\!(t)\!))$ is dense, and likewise for $U(k(\!(t - 1)\!)) \subset X(k(\!(t -1)\!))$: to see this, by working locally on $X(k(\!(t)\!))$ and using the local structure of smooth varieties, we reduce to the case when $X(k(\!(t)\!))$ is still nonempty and there is an étale map $f \colon X \rightarrow \mathbb A^d_k$  and then note that, by \cite{GGMB14}*{Proposition 3.1.4} and the Hilbert Nullstellensatz, the image of any nonempty open of $X(k(\!(t)\!))$ is Zariski dense in $\mathbb A^d_{k(\!(t)\!)}$, whereas that of $(X \setminus U)(k(\!(t)\!))$ is not because it lies in the lower-dimensional closed $\overline{f(X\setminus U)} \subset \mathbb A^d_k$. 

For the remaining density of the image of \eqref{eqn:D}, we may therefore assume that $X = U$. Then, since maps of varieties induce continuous maps on $k(\!(t)\!)$-points (resp.,~$k(\!(t - 1)\!)$-points), by considering preimages under $r$, we may even assume that $X = V$ in \eqref{eqn:D}. A $k(t)$-point of $\mathbb A^n_k$ lies in a given nonempty open if and only if its induced $k(\!(t)\!)$-point does, so finally we need to show that $k(t) \subset k(\!(t)\!) \times k(\!(t - 1)\!)$ is dense for the product topology. This, however, follows from the ability to multiply by powers of $t(t - 1)$ and from the density of $k[t]$ in $k[\![t]\!] \times k[\![t - 1]\!]$. 
\epf

Conversely, vanishing of $K_1^G(-)$ detects retract rationality thanks to the following result of Gille.

\bprop[Gille] \label{prop:univK1}
A reductive group $G$ over a field $k$ satisfies $K_1^G(\ell) \cong 1$ for every field extension $\ell/k$ if and only if it is semisimple, simply connected, and in its canonical decomposition $G \cong \prod \mathrm{Res}_{k'/k}(G')$ as in \eqref{eqn:can-dec} each $G^\prime$ is isotropic and retract rational.
\eprop

%For the canonical decomposition $G \cong \prod \mathrm{Res}_{k'/k}(G')$ of a semisimple, simply connected $k$-group $G$, see \cite{SGA3IIInew}*{exposé XXIV, proposition 5.10 (i) (with section 5.3)}.

\bpf
The formation of $K_1^G(\ell)$ commutes with products of and restrictions of scalars of groups (see \S\ref{pp:K1Gk}), so the `if' is immediate from \Cref{prop:W-R} and \Cref{thm:RR}. For the `only if,' \cite{Gil09}*{proposition 3.1} ensures that $G$ is semisimple, simply connected, and each $\mathrm{Res}_{k^\prime/k}(G^\prime)$ is isotropic, so the $G^\prime$ are all isotropic. Thus, we may assume that $G \cong \mathrm{Res}_{k^\prime/k}(G')$ and, since $G^\prime$ is a direct factor of $G_{k^\prime}$, that $G \cong G^\prime$. The retract rationality of $G$ then follows from \cite{Gil09}*{théorème~5.9}.
\epf

For the triviality of $K_1^G(k)$, it is thus helpful to know which groups are retract rational, or even rational. The following nonexhaustive result gathers some such cases, compare with \cite{Gil09}*{théorème~6.1}.

\bthm\label{prop:retract}
Let $G$ be a semisimple, simply connected, absolutely indecomposable group over~a~field~$k$.
\benum
\m \label{m:R-a}
If $G$ is quasi-split \up{this is automatic when $k$ is finite or, more generally, is of dimension $\le 1$ in the sense that $\mathrm{Br}(k^\prime) \cong 0$ for all finite $k^\prime/k$}, then $G$ is rational over $k$.

\m
Suppose that $G$ is isotropic and $\mathrm{char}(k) \neq 2$. In the following cases, $G$ is rational over $k$\ucolon
\benumr
\m
\label{m:R-i} if $G$ has a multiple Dynkin edge \up{so is of type either $B_n$, $C_n$, $F_4$, or $G_2$}\uscolon or

\m
\label{m:R-ii} if $G$ is a trialitarian $D_4$\uscolon or 

\m
\label{m:R-iii} if $G$ is a nontrialitarian $D_n$ and its central simple algebra of \uS\uref{pp:conv} \ref{m:conv-ii} has index $\le 4$\uscolon or

\m
\label{m:R-iv} if $G$ is an inner type $E_6$\uscolon or

\m
\label{m:R-v} if $G$ is an outer type $E_6$ and has isotropic rank $\ge 2$ over $k$\uscolon or

\m
\label{m:R-vi} if $G$ is an $E_7$ and has isotropic rank $\ge 3$ over $k$\uscolon or

\m
\label{m:R-vii} if $G$ is an $E_8$ and has isotropic rank $\ge 3$ over $k$.
\eenum

\m\label{m:R-c} 
In the following cases, $G$ is retract rational over $k$\ucolon
\benumr \addtocounter{enumii}{7}
\m
\label{m:R-viii} if $G$ is an $A_n$ and its central simple algebra of \uS\uref{pp:conv} \ref{m:conv-i} has squarefree index \up{this is automatic if $n + 1$ is squarefree}\uscolon or

\m
\label{m:R-ix} if $G$ is isotropic and is as in either of \ref{m:R-i}–\ref{m:R-vii} over any field $k$\uscolon or

\m
\label{m:R-xi} if $G$ is an $E_n$ and has isotropic rank $\ge 2$ over $k$\uscolon or

\m
\label{m:R-x} if $G$ has Tits index ${}^2E_{6,1}^{29}$\uscolon or

\m
\label{m:R-xii} if $G$ has Tits index $E_{7, 1}^{78}$. %or $E_{8, 2}^{78}$.

\eenum
\eenum
In particular, in all of these cases \ref{m:R-a}–\ref{m:R-c} we have $G(\ell)/R \cong 1$ for every field extension $\ell/k$. %, and if $G$ is isotropic, then also $K_1^G(\ell) \cong 1$. 
\ethm

\bpf
The final assertion follows from the rest and from \Cref{thm:RR}.

In \ref{m:R-a}, $G$ is rational by  the big cell open immersion and by the inducedness of the maximal torus of a Borel (our $G$ is simply connected) \cite{SGA3IIInew}*{exposé XXIV, proposition 3.13; exposé~XXVI, théorème 4.3.2 (vi)}, compare with \cite{Gil09}*{bottom of p.~63}. The parenthetical follows from \cite{BS68}*{Section 8.6} (with \cite{Con14}*{Proposition 7.2.12}). 

In \ref{m:R-i}, $G$ is rational by \cite{CP98}*{Proposition 4, beginning of Section 5, Proposition 12, Theorem 2 (note the characteristic $\neq 2$ blanket assumption on the bottom of page 5)}, see also \cite{CM98}*{Corollary~2.3}. In \ref{m:R-ii}, $G$ is rational by \cite{CP98}*{Proposition 13}. In \ref{m:R-iii}, $G$ is rational by \cite{CP98}*{Corollaries~5 and 6, Propositions 10 (1) and 11} (for agreement with the explicit description \S\ref{pp:conv} \ref{m:conv-ii}, use \cite{KMRT98}*{Theorem 3.1 (with Theorem 1.1 and Corollary 2.8) and Theorem 4.2 (1)}). In \ref{m:R-iv}, $G$ is rational by \cite{CP98}*{Proposition 14}. In \ref{m:R-v}, $G$ is rational by \cite{CP98}*{Proposition 16} and \cite{Tit66}*{page~59}. In \ref{m:R-vi}, $G$ is rational by \cite{CP98}*{Proposition 17 (a)–(b)} and \cite{Tit66}*{page 59} (and \ref{m:R-a}). In \ref{m:R-vii}, $G$ is rational by \cite{CP98}*{Proposition 17 (c)} and \cite{Tit66}*{page 60} (and \ref{m:R-a}). 

The inner case of \ref{m:R-viii} was treated in \cite{Gil09}*{Proposition 2.4}. The outer case of \ref{m:R-viii} follows by first using \cite{CM98}*{Corollary 6.7} and \cite{KMRT98}*{Theorem 3.1 (2)} to reduce to when $G$ is isotropic and then combining \Cref{prop:univK1} with \Cref{prop:W-R} and \cite{CM98}*{Corollary 5.6} (we explain some of the identifications that are used in the latter in the proof of \Cref{thm:M} below).

For \ref{m:R-ix}–\ref{m:R-xii}, we will now use a standard specialization argument to reduce to characteristic $0$. By \cite{SP}*{Lemma \href{https://stacks.math.columbia.edu/tag/0328}{0328}}, there is a complete discrete valuation ring $\mathcal O$ with residue field $k$ and fraction field $K$ of characteristic $0$. By \cite{Mat89}*{Theorem 29.1}, for every field extension $\ell/k$, there is a complete discrete valuation ring extension $\mathcal O^\prime$ of $\mathcal O$ whose residue field is $\ell$ (and whose fraction field $L$ is an extension of $K$). By \cite{torsors-regular}*{Proposition 6.1.3 (b)} and, in the subcase \ref{m:R-iii} of \ref{m:R-ix}, instead by \cite{CF15}*{corollaire 2.7.0.32, propositions 4.4.0.34 and 4.5.2.7} and \cite{Hitchin-torsors}*{Proposition 2.1.4}, our $G$ lifts to a semisimple, simply connected $\mathcal O$-group scheme $\mathcal G$ that in the subcase \ref{m:R-iii} of \ref{m:R-ix} is of the form $\mathrm{Spin}( A, \sigma, f)$ for data $(A, \sigma, f)$ over $\mathcal O$ that lift the corresponding data for $G$ over $k$ from \S\ref{pp:conv}~\ref{m:conv-ii}.  By \cite{SGA3IIInew}*{exposé XXII, Proposition 2.8}, the generic fiber $\mathcal G_K$ geometrically has the same Dynkin diagram as $G$ and, by \cite{SGA3IIInew}*{exposé XXIV, sections 3.1–3.7}, the Galois actions on these Dynkin diagrams are compatible. By \cite{SGA3IIInew}*{exposé XXVI, section 7.15} and \cite{Guo24}*{Proposition 4.3}, the isotropic $K$-rank of $\mathcal G_K$ equals the isotropic $k$-rank of $G$: every parabolic $k$-subgroup of $G$ lifts to a parabolic $\mathcal O$-subgroup of $\mathcal G$ and, conversely, every parabolic $K$-subgroup of $\mathcal G_K$ extends to a parabolic $\mathcal O$-subgroup of $\mathcal G$. In particular, by considering the dimension of the anisotropic kernel and using \cite{Tit66}*{pages 58–60}, we see that in each case $\mathcal G_K$ is of the corresponding type \ref{m:R-ix}–\ref{m:R-xii}. Finally, by \cite{Gil09}*{corollaire 7.3}, for every field extension $\ell/k$, in the notation above,
\[
K_1^G(\ell) \cong K_1^{\mathcal G_K}(L).
\]
In conclusion, by \Cref{prop:univK1}, each \ref{m:R-ix}–\ref{m:R-xii} reduces to its case when $k$ is of characteristic $0$. The case \ref{m:R-ix} then follows from \ref{m:R-i}–\ref{m:R-vii}. The case \ref{m:R-x} follows from \cite{Gar07}*{Proposition 0.1}. 
 In the case \ref{m:R-xi}, by \cite{Tit66}*{page 59} and \ref{m:R-iv}--\ref{m:R-vii}, we are left with $G$ having one of
\[
\begin{array}{@{}c@{}}
\begin{tikzpicture}[
  x=0.8cm,
  y=0.8cm,
  line width=0.5pt,
  font=\Large,
  baseline={(0,0)}
]
  % Tits index.
  \node[anchor=east, inner sep=0pt] at (-1.05,0) {\scalebox{0.8}{$E_{7,2}^{31}$}};

  % Edges.
  \draw (0,0) -- (5,0);
  \draw (3,0) -- (3,1);

  % Vertices.
  \foreach \x in {0,...,5} {
    \fill (\x,0) circle[radius=2.1pt];
  }
  \fill (3,1) circle[radius=2.1pt];

  % Distinguished vertices: second from the left and rightmost.
  \foreach \x in {1,5} {
    \draw (\x,0) circle[radius=0.2cm];
  }
\end{tikzpicture}
\hspace{1.5cm}
\begin{tikzpicture}[x=0.8cm, y=0.8cm, line width=0.5pt, font=\Large, baseline={(0,0)}]
  \node[anchor=east, inner sep=0pt] at (-1.05,0) {\scalebox{0.8}{$E_{8,2}^{66}$}};

  \draw (0,0) -- (6,0);
  \draw (4,0) -- (4,1);

  \foreach \x/\i in {0/8,1/7,2/6,3/5,4/4,5/3,6/1} {
    \fill (\x,0) circle[radius=2.1pt];
    \node[below=6pt] at (\x,0) {\scalebox{0.8}{\ }};
  }

  \fill (4,1) circle[radius=2.1pt];
  \node[right=6pt] at (4,1) {\scalebox{0.8}{\ }};

  \draw (0,0) circle[radius=0.2cm];
  \draw (6,0) circle[radius=0.2cm];
\end{tikzpicture}
\\[4ex]
\hspace{-0.8cm}% shifts the diagram below 0.4cm to the left
\begin{tikzpicture}[x=0.8cm, y=0.8cm, line width=0.5pt, font=\Large, baseline={(0,0)}]
  % Tits index.
  \node[anchor=east, inner sep=0pt] at (-1.05,0) {\scalebox{0.8}{$E_{8,2}^{78}$}};

  % Edges.
  \draw (0,0) -- (6,0);
  \draw (4,0) -- (4,1);

  % Vertices.
  \foreach \x in {0,...,6} {
    \fill (\x,0) circle[radius=2.1pt];
  }
  \fill (4,1) circle[radius=2.1pt];

  % Distinguished vertices: the two leftmost (vertices 8 and 7);
  % the remaining vertices form the anisotropic kernel E_6.
  \foreach \x in {0,1} {
    \draw (\x,0) circle[radius=0.2cm];
  }
\end{tikzpicture}
\end{array}
\]
as its Tits index over $k$. For these, we start with $E_{7,2}^{31}$ and, by \Cref{prop:univK1}, need to argue that $K_1^G(\ell) \cong 1$ for every field extension $\ell/k$. Moreover, \ref{m:R-vi} lets us  assume that the Tits index of $G_\ell$ is still $E_{7, 2}^{31}$, in other words, we may assume that $\ell = k$. Then the vanishing of $K_1^G(k)$ was observed already in \cite{Tit78}*{section 1.1.7} and may be argued as follows. By the Prasad–Raghunathan reduction \cite{PR85}*{Theorem A and Section 1.6} (conjectured by Tits) about removing all circled nodes except for a single Galois orbit (see also \cite{Gil09}*{Section 3.3}), in the $E_{7, 2}^{31}$ case $K_1^G(k)$ is generated by the images of the $K_1^{G_{\Theta}}(k)$, where $G_\Theta$ is a semisimple, simply connected, absolutely indecomposable $k$-subgroup of isotropic rank $1$ and either of type $D_5$ or $D_6$. In the presentation \S\ref{pp:conv}~\ref{m:conv-ii} of $G_\Theta$ as $\mathrm{Spin}(A, \sigma, f)$, the central simple $k$-algebra $A$ has degree $10$ or $12$, so its index is $\le 4$. Thus, $K_1^{G_\Theta}(k) \cong 1$ by \ref{m:R-iii}.

In the case of $E_{8,2}^{66}$, the retract rationality follows from \cite{PTW12}*{Theorem 1.2}, \ref{m:R-vii}, and \Cref{prop:univK1}. Similarly, for $E_{8,2}^{78}$, the retract rationality follows from either \cite{ACP21}*{Theorem 8.1} or \cite{Tha22}*{Theorem 4.2}, \ref{m:R-vii}, and \Cref{prop:univK1}. Finally, in \ref{m:R-xii}, the retract rationality follows from either \cite{ACP21}*{Theorem 8.1} or \cite{Tha22}*{Theorem~5.1}, \ref{m:R-xi}, and \Cref{prop:univK1}.
\epf

\section{Lifting $R$-equivalence from a homogeneous space} \label{sec:R-homog}

A broadly useful technique for studying $R$-equivalence is controlling it in short exact sequences, see, for instance, \cite{BKG04}*{Appendix, Theorem 1}. However, interesting short exact sequences of reductive groups are scarce, so we now extend such control to quotients by possibly nonnormal algebraic subgroups. This is a bit similar to \cite{GS25}*{Lemma~2.6~(2)}, however, our assumptions seem sharper, for instance, in \ref{m:nR-iii} we require much less than the vanishing of  the entire $H^1_{\mathrm{\text{ét}}}(k(t), M)$. Imposing the latter vanishing would, critically, give a much less practical statement.

\bprop\label{prop:nonabelian-R}
Let $k$ be a field, let $M \subset L$ be an inclusion of $k$-group schemes locally of finite type, and consider the homogeneous space $X \colonequals L/M$. Suppose that
\benumr
\m
\label{m:nR-i} $M$ is smooth over $k$\uscolon

\m
\label{m:nR-ii} $M(k)/R = 1$\uscolon
\eenum
and that either 
\benumr \addtocounter{enumi}{2}
\m\label{m:nR-iii} $\mathrm{Ker}(H^1_{\mathrm{\text{\upshape{ét}}}}(k(t), M) \rightarrow H^1_{\mathrm{fppf}}(k(t), L)) = \{*\}$\uscolon or

\primeitem %[\upshape{(iii${}^\prime$)}]
\label{iiipr-0} 
$k$ is perfect and $\mathrm{Ker}(H^1_{\mathrm{\text{\upshape{ét}}}}(k'(\!(t)\!), M) \rightarrow H^1_{\mathrm{fppf}}(k'(\!(t)\!), L)) = \{*\}$ for each finite extension $k^\prime/k$\uscolon~or

\primeprimeitem %[\upshape{(iii${}^\prime$)}]
\label{iiipr} 
$M$ is \up{connected} reductive and $\mathrm{Ker}(H^1_{\mathrm{\text{\upshape{ét}}}}(K, M) \rightarrow H^1_{\mathrm{\text{\upshape{fppf}}}}(K, L)) = \{*\}$ for each extension $K/k$ such that there are a finite extension $k^\prime/k$ generated by a single element and isomorphisms $\iota \colon k^\prime(\!(t)\!) \isomto K$ \up{possibly not as $k$-algebras} and $\iota_*(M_{k^\prime(\!(t)\!)}) \simeq M_K$ \up{equivalently, instead of the triviality of the kernel, one may require the triviality of its $\iota$-preimage in $H^1(k'(\!(t)\!), M)$}.
\eenum
Then $L(k)/R \xrightarrow{\sim} X(k)/R$, more precisely, then two elements $l, l^\prime \in L(k)$ are \up{resp., directly} $R$-equivalent if and only if their images in $X(k)$ are \up{resp., directly} $R$-equivalent.
\eprop

%As we will show in the proof, \ref{iiipr}$\Rightarrow$\ref{iiipr-0}$\Rightarrow$\ref{m:nR-iii}. 
We will arrange \ref{iiipr} via the Rost invariant, using its insensitivity to field or group isomorphisms.

\bpf
We will argue that every $\chi \in X(\mathscr O_{\mathbb A^1_k,\, 0,\, 1})$ lifts to some $\lambda \in L(\mathscr O_{\mathbb A^1_k,\, 0,\, 1})$, even with the additional requirement that we may prescribe the $k$-points $\lambda|_0 = l$ and $\lambda|_1 =l'$ that lift $\chi|_0$ and $\chi|_1$. By considering those $\chi$ that are pullbacks of $k$-points of $X$, this will show that $L(k) \rightarrow X(k)$ is surjective, and then, since $R$-equivalence is functorial, we will also get the desired claim.

By the long exact nonabelian cohomology sequence, \ref{m:nR-iii} amounts to the map $L(k(t)) \rightarrow X(k(t))$ being surjective and, as we will argue at the end, either \ref{iiipr-0} or  \ref{iiipr} implies \ref{m:nR-iii}. The preimage of $\chi$ in $L$ is a torsor under $M$ over $\mathscr O_{\mathbb A^1_k,\, 0,\, 1}$, so \ref{m:nR-iii} ensures that this $M$-torsor is generically trivial. Thus, by \ref{m:nR-i} and the Grothendieck–Serre type result \cite{torsors-algebraic}*{Theorem 1.1.1} (already \cite{Gil97}*{corollaire~I.2.2~b)} suffices when $M$ is reductive), the torsor is trivial already over $\mathscr O_{\mathbb A^1_k,\, 0,\, 1}$. Concretely, this means that $\chi$ lifts to some $\nu \in L(\mathscr O_{\mathbb A^1_k,\, 0,\, 1})$, and this already gives the existence of $\lambda$ without the finer aspect about prescribing $\lambda|_0$ and $\lambda|_1$. To arrange the latter, we first note that the preimages of $\chi|_0$ and $\chi|_1$ in $L$ are $M$-torsors over $k$ trivialized by $l$ and $l^\prime$, respectively, so 
$$
l = \nu|_0m \quad \text{and}\quad l' = \nu|_1m' \quad \text{for some} \quad m,m' \in M(k).
$$
The assumption \ref{m:nR-ii} and the agreement with direct $R$-equivalence (see \S\ref{pp:R-equivalence}) give a $\mu \in M(\mathscr O_{\mathbb A^1_k,\, 0,\, 1})$ with $\mu|_0 = m$ and $\mu|_1 = m'$. We then set $\lambda \colonequals \nu \mu$.

To show that either \ref{iiipr-0} or \ref{iiipr} implies \ref{m:nR-iii}, we use a trick of Gabber from the proof of \cite{Gil97}*{lemme~II.2.7}. For this, an $M$-torsor $E$ over $k(t)$ is trivial if and only if it trivializes over the fraction field $K_v$ of the completion $\mathcal O_v$ of $\mathbb P^1_k$ at every closed point $v$: by spreading out and patching, $E$  then extends to an $M$-torsor over $\mathbb P^1_k$ that is trivial at $\infty$, and hence, by the classification of $M$-torsors over $\mathbb P^1_k$ from \cite{Gil02}*{théorème 3.8 (a)} (or \cite{torsors-algebraic}*{Theorem~5.2.4}), it is trivial over $k(t)$. By the Cohen structure theorem \cite{Mat89}*{Theorem 29.7 and its proof, as well as Theorem 29.6}, letting $k_v$ be the residue field of $\mathbb P^1_k$ at $v$, we have an isomorphism $\iota \colon k_v[\![t]\!] \isomto \mathcal O_v$ compatible with the reduction modulo the maximal ideal, and even as $k$-algebras if $k_v/k$ is separable (but possibly not if it is inseparable). Thus, if $k$ is perfect, then $k_v/k$ is separable, $K_v \cong k_v(\!(t)\!)$ as $k$-algebras, and~\ref{iiipr-0}$\Rightarrow$\ref{m:nR-iii}. 

When $k$ is not necessarily perfect but $M$ is reductive, at least $\iota$ is compatible with the reduction modulo the maximal ideal, so the reductive $\cO_v$-group schemes $M_{\mathcal O_v}$ and $\iota_*(M_{k_v[\![t]\!]})$ are residually isomorphic. Thus, by \cite{torsors-regular}*{Proposition 6.1.3 (a)}, there is an $\cO_v$-group isomorphism $\iota_*(M_{k_v[\![t]\!]}) \simeq M_{\mathcal O_v}$. Since $k_v$ is generated by a single element over $k$, by taking $K = K_v$ in \ref{iiipr} we get that \ref{iiipr}$\Rightarrow$\ref{m:nR-iii}.
\epf

%\brem
%In practice it is easier to check the condition \ref{iiipr} rather than \ref{m:nR-iii} (which it implies). This is because a complete discrete valuation makes $M$-torsors more accessible over $k^\prime(\!(t)\!)$ than over $k(t)$, see \Cref{thm:Rost-A2} below and its proof for a concrete instance of this.
%\erem

\beg
Although this is not its main case of interest, \Cref{prop:nonabelian-R} applies when $M$ is a quasi-trivial torus. In particular, it recovers the possibility of passing to the quotient by the maximal split central torus when studying $R$-equivalence (such a passage is sometimes useful for dealing with Levi subgroups), compare with \cite{GS25}*{Lemma 2.8, Corollary 3.10}.
\eeg

%%%%%%%%%%%%%%%%%%%%%%%%%%%%%%%%%%%%%%%%%%%
 
\section{The kernel of the Rost invariant for simply connected type $A_n$} \label{sec:Rost}

In \Cref{prop:nonabelian-R}, the most delicate condition to arrange is the vanishing of the pointed set in either \ref{m:nR-iii} or \ref{iiipr}. We will ensure it by using the Rost invariant and applying \Cref{thm:Rost-A1,thm:Rost-A2}. 

\bpp[The Rost invariant] \label{pp:Rost}
For a semisimple, simply connected, absolutely indecomposable group $G$ over a field $k$, its Rost invariant is a natural transformation of functors
\be \label{eqn:Rost}
r_{G, -} \colon H^1_{\mathrm{\text{ét}}}(-, G) \rightarrow H^3_{\mathrm{\text{ét}}}(-, (\mathbb Q/\mathbb Z)(2))
\ee
on the category of field extensions $k^\prime/k$ that is normalized, meaning $r_{G,\, k}(*) = 0$, and that generates the finite cyclic group  of all such normalized natural transformations, see \cite{GMS03}*{page 129, Theorem~9.11 and before Theorem 9.10 (as well as footnote a on p. 152)} (and \S\ref{pp:conv}). \emph{Loc.~cit.}~shows that the formation of the Rost invariant commutes with restriction to an overfield, in fact, it constructs an explicit generator of this finite cyclic group for which such base change is clear. This finer construction is delicate and comparing it with alternative (no more direct) definitions of the Rost invariant is serious work, see \cite{GM17}*{Sections 1 and 6C} and the references there. We will only be interested in the vanishing (or not) of some values of the Rost invariant, so the cruder definition with ambiguity up to an invertible multiple will suffice. A useful criterion for recognizing when a normalized natural transformation $r_{G, -}$ as in \eqref{eqn:Rost} is a Rost invariant of $G$ is given by \cite{GMS03}*{page 109, Theorem 3.3 (with page 154, Proposition A.9) and page 129, Theorem~9.10}: it suffices that for every $G$-torsor $E$ over a field extension $k^\prime/k$, we have
\[
\langle r_{G, k^\prime}([E])\rangle \cong \mathrm{Ker}\left( H^3_{\mathrm{\text{ét}}}(k^\prime, (\mathbb Q/\mathbb Z)(2)) \rightarrow  H^3_{\mathrm{\text{ét}}}(k^\prime(E), (\mathbb Q/\mathbb Z)(2)) \right).
\]
For example, by \cite{Gil00}*{théorème 1}, the Rost invariant used in \emph{op.~cit.}~satisfies this.

For a $k$-homomorphism $f\colon G \rightarrow G^\prime$ of semisimple, simply connected, absolutely indecomposable $k$-groups, its \emph{Rost multiplier} is the unique $n_f \in \mathbb Z_{\ge 0}$ such that a generator of the rank $1$, free $\mathbb Z$-module of Weyl-group invariant integral quadratic forms on the cocharacter lattice of $G^{\prime}_{\overline{k}}$ restricts to $n_f$ times an analogous generator for $G$, see \cite{GMS03}*{page 122, before Proposition 7.9}. By \cite{GMS03}*{page~129, Theorem 9.11}, the restriction of a Rost invariant of $G'$ is $n_f$ times a Rost invariant of $G$.
\epp

\beg \label{eg:Rost}
Let $G$ be of inner type $A_{n - 1}$, so $G \simeq \mathrm{SL}_1(A)$ for a central simple algebra $A$ of degree $n$ over $k$ (see \S\ref{pp:conv}~\ref{m:conv-i}). Then $H^1_{\mathrm{\text{ét}}}(k, G) \cong k^\times/\mathrm{Nrd}(A^\times)$ and, by \cite{GQM11}*{théorème 1.1} (and the end of \S\ref{pp:conv}), the Rost invariant of a $G$-torsor corresponding to the class of an $x \in k^\times$ is the cup product 
\[
[x] \cup [A] \in H^3_{\mathrm{\text{ét}}}(k, (\mathbb Q/\mathbb Z)(2))
\]
(for a good choice in \eqref{eqn:Rost}), where $[x] \in k^\times/k^{\times n} \cong H^1_{\mathrm{\text{fppf}}}(k, \mu_n)$ and $[A] \in \mathrm{Br}(k) \cong H^2_{\mathrm{\text{ét}}}(k, (\mathbb Q/\mathbb Z)(1))$.
\eeg

\bpp[The Rost invariant and the center] \label{pp:R-center}
For a semisimple, simply connected, absolutely indecomposable group $G$ over a field $k$, Garibaldi–Merkurjev \cite{GM17}*{Theorem 1.2} gave a formula for the Rost invariant on the image of $H^1_{\mathrm{fppf}}(k, Z_G) \rightarrow H^1_{\mathrm{\text{ét}}}(k, G)$  (see also the earlier \cite{MPT02}). They let
\benumr
\m
\label{RCi} $t_G \in H^2_{\mathrm{fppf}}(k, Z_G)$ be the Tits class of $G$, namely, the negative of the boundary image of the unique class in $H^1_{\mathrm{\text{ét}}}(k, G^{\mathrm{ad}})$ for which the corresponding inner form of $G$ is quasi-split, see, for instance, \cite{KMRT98}*{Proposition 31.5 and before Proposition 31.7};

\m
\label{RCii} $Z_G^\circ$ be the finite $k$-group of multiplicative type that is “dual” to $Z_G$ in the sense that their character groups are dual: $X^*(Z_G^\circ) \cong \mathrm{Hom}(X^*(Z_G), \mathbb Q/\mathbb Z)$;

\m
\label{RCiii} $\widehat{\rho}\colon Z_G \rightarrow Z_G^\circ$ be the $k$-homomorphism whose dual is induced by the unique map from the coweight to the weight lattice of a maximal torus of $G$ that sends each short coroot to the corresponding root, see \cite{GM17}*{Proposition 5.2, Section 6A}; by \emph{loc.~cit.},~$\widehat{\rho}$ is an isomorphism when the Dynkin diagram of $G$ is simply laced (the only case we will use).
\eenum
By \cite{GM17}*{Theorem 1.2}, the Rost invariant of $G$ composes with $H^1_{\mathrm{fppf}}(k, Z_G) \rightarrow H^1_{\mathrm{\text{ét}}}(k, G)$ to give
\be\label{eqn:R-center}
H^1_{\mathrm{fppf}}(k, Z_G) \xrightarrow{- \cup H^2_{\mathrm{fppf}}(k, \widehat \rho)(-t_G)} H^3_{\mathrm{\text{ét}}}(k, (\mathbb Q/\mathbb Z)(2))
\ee
(for a good choice in \eqref{eqn:Rost}), where the cup product map is as in \cite{GM17}*{(2.9)}. For a $k$-subgroup $i\colon C \subset Z(G)$, by \cite{GM17}*{Lemma~2.16},  the Rost invariant of $G$ on $\mathrm{Im}(H^1_{\mathrm{fppf}}(k, C) \rightarrow H^1_{\mathrm{\text{ét}}}(k, G))$ is given by the same formula but with $H^2_{\mathrm{fppf}}(k, \widehat \rho)(-t_G)$ replaced by its $H^2_{\mathrm{fppf}}(k, i^\circ)$-image in $H^2_{\mathrm{fppf}}(k, C^\circ)$.

\epp

To make our results sharp in positive characteristic, we use the following cohomological dimension.

\bpp[The separable dimension] \label{pp:sep-dim}
For a field $k$ and a prime $p$, the \emph{separable $p$-dimension} $\mathrm{dim}^{\mathrm{sep}}_p(k)$~is 
\benumr
\m \label{m:dimsep-i} the $p$-cohomological dimension $\mathrm{cd}_p(k)$ of $\mathrm{Gal}(k^s/k)$ in the usual sense of Serre if $p \neq \mathrm{char}(k)$;

\m
\label{m:dimsep-ii} $\mathrm{inf}\{d \ge 0\, \vert\, H^1_{\mathrm{\text{ét}}}(k', \underline{\Omega}^d_{k',\, \mathrm{log}}) \cong 0\text{ for all finite separable extensions }k^\prime/k\}$ if $p = \mathrm{char}(k)$;
\eenum
it was introduced by Gille in \cite{Gil00}*{Section 1.1 (c)}. %(see also \eqref{eqn:Cartier} and the end of \S\ref{pp:conv} for the agreement between the \ref{m:dimsep-ii} above and \emph{loc.~cit.}). 
It refines the \emph{$p$-dimension} $\mathrm{dim}_p(k)$ defined as
\benumr \addtocounter{enumi}{1}
\primeitem \label{m:dim-i}
the $p$-cohomological dimension $\mathrm{cd}_p(k)$ of $\mathrm{Gal}(k^s/k)$ in the usual sense of Serre if $p \neq \mathrm{char}(k)$;

\addtocounter{enumi}{1}
\primeitem \label{m:dim-ii}
$\mathrm{inf}\{d \ge 0\, \vert\, [k : k^p] \le p^d \text{ and } H^1_{\mathrm{\text{ét}}}(k^\prime, \underline{\Omega}^d_{k^\prime,\, \mathrm{log}}) \cong 0 \text{ for all finite extensions }k^\prime/k\}$ if $p = \mathrm{char}(k)$;
\eenum
it was introduced by Kato in \cite{Kat82}*{p. 220, Definition} (by \cite{SP}*{\href{https://stacks.math.columbia.edu/tag/07P2}{Lemma 07P2}}, if $[k : k^p] = p^d$, then even $\underline{\Omega}^{d + 1}_{k^\prime,\, \mathrm{log}} \cong 0$ for all finite extensions $k^\prime/k$). The \emph{separable dimension} and the \emph{dimension} of $k$ are
\[
\mathrm{dim}^{\mathrm{sep}}(k) \colonequals \mathrm{sup}_p(\mathrm{dim}^{\mathrm{sep}}_p(k)) \qxq{and} \mathrm{dim}(k) \colonequals \mathrm{sup}_p(\mathrm{dim}_p(k)),
\]
respectively. To recast without a variable $k^\prime$, by \Cref{thm:dimp} (and \S\ref{pp:conv}), $\mathrm{dim}_p(k)$ (resp., $\mathrm{dim}_p^{\mathrm{sep}}(k)$) is the infimum of all $d \ge 0$ such that $H^i_{\mathrm{\text{ét}}}(k, \mathscr F(j)) \cong 0$ for all $p^\infty$-torsion abelian sheaves $\mathscr F$ on $k_{\mathrm{\text{ét}}}$ and all $i, j \ge0$ with $i > d$ (resp., and all $i > j \ge 0$ with $i > d$); by a limit argument, we can even restrict to  $\mathscr F$ that are commutative, finite, étale $k$-groups of $p$-power order. In addition, by \Cref{thm:dimp}, all the infima above may be replaced by minima if the corresponding sets are~nonempty.

We have $\mathrm{dim}^{\mathrm{sep}}_p(k) \le \mathrm{dim}_p(k)$ and, if $k$ is perfect, then $\mathrm{dim}^{\mathrm{sep}}_p(k) = \mathrm{dim}_p(k) = \mathrm{cd}_p(k)$: in the key case $p = \mathrm{char}(k)$, these dimensions are $1$ granted that some finite separable $k'/k$ has a $\mathbb Z/p\mathbb Z$-extension, else they are all $0$, see \cite{Gil00}*{Section 1.1 (b) and (c)} and \cite{Ser02}*{Chapter~I, Section~4.1, Proposition~21${}^\prime$}. 

The condition $\mathrm{dim}_p^{\mathrm{sep}}(k) \le 2$ is significant for us because it forces the vanishing of the $p$-primary part of the Rost invariant (even of its target). In addition, by a result of Suslin and Gille \cite{Gil00}*{théorème~7}, the bound $\mathrm{dim}_p^{\mathrm{sep}}(k) \le 2$ is equivalent to the surjectivity of the reduced norm map $\mathrm{Nrd} \colon D^\times \rightarrow k^{\prime \times}$ for  every central division algebra $D$ of $p$-power degree over every finite separable extension $k^\prime/k$ (or even over every algebraic extension $k'/k$ thanks to the following remark). 
\epp

\brem \label{rem:dimp-finite-extn}
By \Cref{thm:dimp}~\ref{m:dimp-a}, if $\mathrm{dim}^{\mathrm{sep}}_p(k) \le d$, then $\mathrm{dim}^{\mathrm{sep}}_p(k') \le d$ for every algebraic $k^\prime/k$. 
\erem

\begin{caution} \label{caution} 
In contrast, the behavior of $\mathrm{dim}_p^{\mathrm{sep}}(-)$ with respect to nonalgebraic field extensions $k^\prime/k$ is less intuitive. For example, suppose that $\mathrm{dim}_p^{\mathrm{sep}}(k) = 0$, and even that $k$ is separably closed but imperfect of characteristic $p$ with $p$-rank $\ge d$ and $p$-independent $x_1, \dotsc, x_d \in k$. 
\benuma
\m \label{m:caution-1}
By \cite{SP}*{Lemma \href{https://stacks.math.columbia.edu/tag/07P2}{07P2}} and  Izhboldin’s residue calculation \cite{Tot22}*{Theorem 4.3}, in terms of the explicit description \eqref{eqn:Cartier}, the class of $\frac 1t \frac{\mathrm{d}x_1}{x_1} \wedge \dotsc \wedge \frac{\mathrm{d}x_d}{x_d}$ in $H^{d + 1}_{\mathrm{\text{ét}}}(k(t), (\mathbb Z/p\mathbb Z)(d))$ is nonzero. In particular, $\mathrm{dim}^{\mathrm{sep}}_p(k(t)) \ge d + 1$, contrary to the transition result \cite{Ser02}*{Section~II.4.2, Proposition 11} for the usual Galois-theoretic cohomological dimension. 

\m \label{m:caution-2}
Similarly, for every complete, discretely valued field $K$ of characteristic $0$ with residue field $k$ and every algebraic extension $K^\prime/K$ with residue field $k$, by combining \cite{GS17}*{Proposition~7.1.4, Lemma 9.5.1} and \cite{Voe11}*{Theorem~6.16}, we have $\mathrm{cd}_p(K') \ge d$. In particular, this shows that various reductions to characteristic $0$ that are based on transfer principles developed in \cite{ILA25} tend to be toothless when one works with the “optimal” hypotheses that involve the separable $p$-dimension $\mathrm{dim}_p^{\mathrm{sep}}(k)$ instead of the $p$-dimension $\mathrm{dim}_p(k)$.
\eenum
\end{caution}

The triviality of the kernel of the Rost invariant under various assumptions is a vast subject that we do not attempt to survey here. For our purposes, we need such triviality over $k(\!(t)\!)$ for a group of type $A_n$ defined over a field $k$ with $\mathrm{dim}^{\mathrm{sep}}_2(k) \le 2$ (even $n = 5$ would suffice, though we see no point in being this specific). This does not seem to be in the literature, especially, in view of \Cref{caution}, although \cite{Pre13}*{Theorem 1.7} and \cite{PS24}*{Theorem 9.4} come close. The inner type case, which we review in \Cref{thm:Rost-A1} for convenience, is an immediate consequence of the Merkurjev–Suslin theorem and of its extension to every characteristic due to Gille. The outer type case, which we develop in \Cref{thm:Rost-A2}, rests on \cite{BFT07} and lightly draws on some ideas from \cite{PS24}. Statements similar to \Cref{thm:Rost-A1,thm:Rost-A2} may hold for wider classes of groups or fields but we do not try to explore this here for the sake of our focus on \Cref{thm:MA}.

\bthm[Merkurjev–Suslin, Gille] \label{thm:Rost-A1}
For a central division algebra $D$ of squarefree index over a field $k$ and an $r \ge 1$, an $\mathrm{SL}_r(D)$-torsor over $k$ is trivial if and only if  its Rost invariant vanishes.
\ethm

\bpf
If the index of $D$ is invertible (resp., vanishes) in $k$, \Cref{eg:Rost} and \cite{Sus84}*{Theorem~24.4, Remark 24.5 (and the paragraph after Lemma 10.8)} (resp., and \cite{Gil00}*{théorème~6~(b)}) suffice.
\epf

\brem
The squarefreeness of the index ensures that the period and the index of $D$ are equal over $k$ (and also over every extension). Without this, even with the cohomological dimension $\le 3$ assumption on $k$, one has the Merkurjev counterexample \cite{CTPS12}*{Remark 5.1} and its variants.
\erem

\bthm \label{thm:Rost-A2}
Let $M$ be a semisimple, simply connected, absolutely indecomposable group of type $A_n$ over a field $k$\ucolon either $M \simeq \mathrm{SL}_r(D)$ or $M \simeq \mathrm{SU}_r(D, h)$ with $D$, $k^\prime$, $\sigma$, $r$, $h$ as in \uS\uref{pp:conv}~\uref{m:conv-i}. If
\benumr
\m\label{m:RA-i} $\mathrm{dim}_2^{\mathrm{sep}}(k) \le 2$\uscolon

\m
\label{m:RA-iii} the odd part of the index of $D$ is squarefree \up{for instance, $n + 1$ is squarefree}\uscolon and

\m \label{m:RA-ii} 
the $2$-adic valuations of the period and the index of $D$ agree \up{for instance, $n + 1$ is squarefree}\uscolon
\eenum
then an $M$-torsor over $k(\!(t)\!)$ is trivial if and only if its Rost invariant vanishes.
\ethm

%Both \ref{m:RA-ii} and \ref{m:RA-iii} are automatic if $n + 1$ is squarefree, see \S\ref{pp:conv}.

\bpf
The index of $D$ divides $n + 1$ (see \uS\uref{pp:conv}), which gives the parentheticals. By \cite{Bla11}*{Theorem~0.2} (whose restriction $\mathrm{char}(k) \neq 2$ is bypassed by lifting to characteristic $0$, as in the proof of \Cref{prop:retract} and as explained in \cite{GP26}*{Section 6.1 (1)}), an $M$-torsor with a zero cycle of degree~$1$ is trivial.

In the case when $M\simeq \mathrm{SL}_r(D)$, by the above, \Cref{rem:dimp-finite-extn} and passage to suitable finite extensions splitting the primary factors of $D$ supplied by \cite{Dra83}*{Section~9, Corollary 11} let us assume that $D$ is of prime power index. \Cref{thm:Rost-A1} and \ref{m:RA-iii} then let us assume that the index of $D$ is a power of~$2$. Thus, by \ref{m:RA-i}, the reduced norms of $D^\times$ exhaust $k^\times$ (see \S\ref{pp:sep-dim}). Consequently, by the smoothness of $\mathrm{SL}_1(D)$ and Hensel lifting, reduced norms of $D[\![t]\!]^\times$ exhaust $k[\![t]\!]^\times$. It follows that every class in 
\[
H^1_{\mathrm{\text{ét}}}(k(\!(t)\!), \mathrm{SL}_r(D)) \cong k(\!(t)\!)^\times/\mathrm{Nrd}_{D(\!(t)\!)/k(\!(t)\!)}(D(\!(t)\!)^\times)
\]
is that of $t^a$ for some $a \in \mathbb Z$ . By \Cref{eg:Rost}, its Rost invariant is $[t^a]  \cup [D] \in H^3_{\mathrm{\text{ét}}}(k(\!(t)\!), (\mathbb Q/\mathbb Z)(2))$. This cup product dies over $k(\!(t)\!)^{\mathrm{nr}}$ because $D_{k^s}$ splits and, by \cite{GM17}*{Example 3.3, Lemma~3.4}, its residue in $H^2_{\mathrm{\text{ét}}}(k, (\mathbb Q/\mathbb Z)(1)) \cong \mathrm{Br}(k)$ is $a[D]$. By \ref{m:RA-ii}, this residue vanishes precisely when $\mathrm{deg}_k(D) \mid a$, which is thus also necessary for the Rost invariant to vanish. Since $t^{\mathrm{deg}_k(D)}$ is a reduced norm of $t \in D(\!(t)\!)^\times$, we get that an $\mathrm{SL}_r(D)$-torsor over $k(\!(t)\!)$ with a trivial Rost invariant is trivial.

In the remaining case when $M \simeq \mathrm{SU}_r(D, h)$ with $D$, $k^\prime$, $\sigma$, $r$, $h$ as in \S\ref{pp:conv}~\ref{m:conv-i}, we have $M_{k^\prime} \simeq \mathrm{SL}_r(D)$. Thus, \Cref{rem:dimp-finite-extn} and the $\mathrm{SL}_r(D)$ case show that every $M$-torsor over $k(\!(t)\!)$ with a trivial Rost invariant trivializes over the quadratic extension $k'(\!(t)\!)/k(\!(t)\!)$. By the opening paragraph, \Cref{rem:dimp-finite-extn} again, and a restriction-corestriction argument to retain \ref{m:RA-ii}, we may therefore replace $k$ by any finite separable extension of odd degree. By \cite{BFT07}*{Lemma 3.3 (1)}, some such extension kills the odd primary part of $D$, so we are left with the case when the index of $D$ is a power of $2$.

We let $\mathrm U_r(D, h)$ be the $k$-group of those $D$-linear automorphisms of $D^r$ that preserve $h$, so
\be\label{ses}
1 \rightarrow \mathrm{SU}_r(D, h) \rightarrow \mathrm{U}_r(D, h) \xrightarrow{\mathrm{Nrd}_{/k^\prime}} \mathrm{Res}^1_{k^\prime/k}(\mathbb G_m) \rightarrow 1
\ee
is an exact sequence, where $\mathrm{Res}^1_{k^\prime/k}(\mathbb G_m) \colonequals \mathrm{Ker}\left(\mathrm{Norm}_{k'/k}\colon \mathrm{Res}_{k^\prime/k}(\mathbb G_m) \rightarrow \mathbb G_m\right)$, see \cite{KMRT98}*{page~402 (with Proposition 4.1)}. By \cite{Knu91}*{Chapter I, Proposition 6.2.4 1)}, our $h$ is an orthogonal sum of $D$-lines. By \ref{m:RA-i} and \cite{BFT07}*{Corollary 3.12} applied to one such line, the map $\mathrm{Nrd}_{/k^\prime}$ is surjective on $k$-points. Its smoothness and Hensel lifting then imply its surjectivity on $k[\![t]\!]$-points, and hence also on $k(\!(t)\!)$-points because the norm $1$ condition implies that $k(\!(t)\!)$-points of $\mathrm{Res}^1_{k^\prime/k}(\mathbb G_m)$ are its $k[\![t]\!]$-points. By the nonabelian cohomology sequence of \eqref{ses}, we are left with showing that an $\mathrm{SU}_r(D, h)$-torsor $E$ over $k(\!(t)\!)$ with a trivial Rost invariant induces a trivial $\mathrm{U}_r(D, h)$-torsor. 

Isomorphism classes of $\mathrm U_r(D, h)$-torsors amount to $D$-isometry classes of nondegenerate $\sigma$-Hermitian forms on $D^r$, and similarly over $k(\!(t)\!)$, see \cite{KMRT98}*{after (29.17)}. Thus, since $D_{k'(\!(t)\!)} \cong D(\!(t)\!)$ is a central division algebra over $k'(\!(t)\!)$, our task is to show that the nondegenerate $\sigma$-Hermitian form
\[
h' \colon D(\!(t)\!)^r \times D(\!(t)\!)^r \rightarrow D(\!(t)\!) 
\]
associated to $E$ is $D(\!(t)\!)$-isometric to $h_{D(\!(t)\!)}$.  For this, similarly to \cite{PS24}*{proofs of Propositions 8.3 and 8.4}, we use \cite{Knu91}*{Chapter I, Proposition 6.2.4 1)} again to diagonalize $h^\prime$, then rescale the diagonal basis by powers of $t$ to make its diagonal self-pairing values have $t$-valuation either $0$ or $1$, and finally combine the smoothness of the isometry scheme (with $r = 1$) with Hensel lifting to get
\be \label{eqn:h}
h^\prime \simeq f_{D(\!(t)\!)} \perp tg_{D(\!(t)\!)} \qxq{with} f \colon D^{a} \times D^{a} \rightarrow D \qxq{and} g \colon D^{b} \times D^{b} \rightarrow D,
\ee
where $f, g$ are nondegenerate $\sigma$-Hermitian forms that are orthogonal direct sums of $D$-lines. 

We recall from \cite{BFT07}*{Section 3.3} that the discriminant $\mathrm{disc}(u)$ of a $\sigma$-Hermitian form $u$ on some $D^d$ is well defined in $k^\times /\mathrm{Norm}_{k^\prime/k}(k^{\prime\times})$ as $(-1)^{m(m - 1)/2}\mathrm{Nrd}_{\mathrm{Mat}_{d \times d}(D)/k^\prime}(M(u))$ where $m \colonequals d \cdot \mathrm{deg}_{k^\prime}(D)$ and $M(u)$ is the pairing matrix of $u$ in some basis of $D^d$. Since $h^\prime$ comes from an $\mathrm{SU}_r(D, h)$-torsor and $\mathrm{SU}_r(D, h) = \mathrm{U}_r(D, h)^{\mathrm{Nrd} = 1}$, we have $\mathrm{disc}(h^\prime) = \mathrm{disc}(h)$. Thus, $\mathrm{disc}(f_{D(\!(t)\!)} \perp tg_{D(\!(t)\!)}) = \mathrm{disc}(h_{D(\!(t)\!)})$, that is, $b \cdot \mathrm{deg}_{k^\prime}(D)$ is even and $\mathrm{disc}(f) \mathrm{disc}(g) = \mathrm{disc}(f \perp g) = \mathrm{disc}(h)$ in $k^\times /\mathrm{Norm}_{k^\prime/k}(k^{\prime\times})$. By \ref{m:RA-i} and \cite{BFT07}*{Corollary 3.9}, it follows that $f \perp g$ is $D$-isometric to $h$. To conclude, it suffices to show that $g$ is hyperbolic, that is, is a sum of hyperbolic $D$-planes: by scaling one of the two basis vectors of each relevant hyperbolic plane over $D(\!(t)\!)$ by $t$, this will show that $tg_{D(\!(t)\!)}$ is $D(\!(t)\!)$-isometric to $g_{D(\!(t)\!)}$, so $h^\prime$ will be $D(\!(t)\!)$-isometric to $h_{D(\!(t)\!)}$, and so $E$ will be trivial. Thus, we may restrict to~$b > 0$. 

To get our hands on $g$, the idea is similar to the case of $\mathrm{SL}_r(D)$: we wish to extract concrete information about $g$ by computing the residue at $t$ of the (vanishing) Rost invariant. One may approach this by adapting \cite{PS24}*{Lemma 4.6} or by reducing to a situation where Gille’s general \cite{Gil00}*{théorème 4} applies. For the sake of brevity, we devise a more indirect but shorter approach.

The $D$-isometry $f \perp g \simeq h$  above gives an inclusion $\iota \colon \mathrm{SU}_b(D, g) \hookrightarrow \mathrm{SU}_r(D, h)$, which over $k^\prime$ is an inclusion $\mathrm{SL}_b(D) \hookrightarrow \mathrm{SL}_r(D)$ of a factor of the derived group of a block Levi. Since $b \cdot \mathrm{deg}_{k^\prime}(D)$ is even, we have a central $\mu_2 \subset \mathrm{SU}_{b}(D, g)$. The class of $t$ in $H^1_{\mathrm{\text{fppf}}}(k(\!(t)\!), \mu_2) \cong k(\!(t)\!)^\times/k(\!(t)\!)^{\times 2}$ gives an $\mathrm{SU}_{b}(D, g)$-torsor $E_t$ over $k(\!(t)\!)$ whose associated $\sigma$-Hermitian form is $tg_{D(\!(t)\!)}$. Consequently, $\iota_*(E_t)$ is an $\mathrm{SU}_r(D, h)$-torsor over $k(\!(t)\!)$ whose associated $\sigma$-Hermitian form is $f_{D(\!(t)\!)} \perp tg_{D(\!(t)\!)}$, which, by \eqref{eqn:h}, is isometric to $h^\prime$. On the other hand, twisting \eqref{ses} by the $\mathrm{SU}_r(D, h)$-torsor $E$ gives
\[
1 \rightarrow \mathrm{SU}_r(D(\!(t)\!), h') \rightarrow \mathrm{U}_r(D(\!(t)\!), h') \xrightarrow{\mathrm{Nrd}_{/k^\prime(\!(t)\!)}} \left(\mathrm{Res}^1_{k^\prime/k}(\mathbb G_m)\right)_{k(\!(t)\!)} \rightarrow 1.
\]
By the same Hensel lifting argument as for \eqref{ses}, $D(\!(t)\!)$-isometries of a single $D$-line summand of either $f$ or $g$ show that the reduced norm map $\mathrm{Nrd}_{/k^\prime(\!(t)\!)}$ is surjective on $k(\!(t)\!)$-points. Thus, up to isomorphism, $E$ is the unique $\mathrm{SU}_r(D, h)$-torsor over $k(\!(t)\!)$ whose associated $\sigma$-Hermitian form is $h^\prime$, so $E \simeq \iota_*(E_t)$. By \cite{GMS03}*{pages 122–123, Proposition 7.9 (4), Example 7.10}, the Rost multiplier of $\iota$ is $1$, so the Rost invariant of $E_t$ vanishes. In effect, by replacing $\mathrm{SU}_r(D, h)$ by $\mathrm{SU}_b(D, g)$ and $E$ by $E_t$, we reduce to when $f = 0$ and $E$ comes from the $\mu_2$-torsor over $k(\!(t)\!)$ associated to $t$.

To now compute the Rost invariant of $E$, we let $i \colon \mu_2 \hookrightarrow Z_M$ be the inclusion into the center. The character group of $Z_M$ is cyclic over $\overline{k}$, so, since the Dynkin diagram of $M$ is simply laced, by \S\ref{pp:R-center}~\ref{RCiii}, the composition $Z_M \xrightarrow{\widehat \rho} Z_M^\circ \xrightarrow{i^\circ} \mu_2^{\circ} \cong \mu_2$ is the unique character of $Z_M$ of order $2$. Thus, by the formula for the Rost invariant reviewed in \S\ref{pp:R-center}, \cite{KMRT98}*{Proposition 31.7 (with Lemma 27.6) and the Type $A_n$ case of Section 27.B}, and the residue formula \cite{GM17}*{Example 3.3, Lemma 3.4} that we have used in the inner type case, the discriminant algebra of $M$ splits. By restricting it to $\mathrm{Br}(k^\prime)$ and using \cite{KMRT98}*{Definition 10.4 (a) and Proposition 10.30}, we find that $(r \cdot \mathrm{deg}_{k^\prime}(D))/2$ kills the class of $D$ in $\mathrm{Br}(k^\prime)$. However, by \ref{m:RA-ii}, the class of $D$ in $\mathrm{Br}(k^\prime)$ has order $\mathrm{deg}_{k^\prime}(D)$, so $r$ must be even (as is necessary for the sought hyperbolicity). Thus, we have verified the conditions of \cite{BFT07}*{Theorem 3.2}, namely, \ref{m:RA-i}, the splitting of the discriminant algebra, and the evenness of $r$, so \emph{loc.~cit.}~(with \cite{BFT07}*{Section 1.4}) shows that $h$, and hence also $g$, is hyperbolic, as desired.
\epf

\brem
In \Cref{thm:Rost-A2}, the vanishing of the Rost invariant is critical: we have only assumed that $\mathrm{dim}^{\mathrm{sep}}_2(k) \le 2$, so nontrivial $M$-torsors (with a nontrivial Rost invariant) may exist even over $k$.
\erem

\section{An $E_{7, 1}^{66}$ determines a $k$-form of a $32$-dimensional half-spin representation} \label{sec:E7166}

We start building the homogeneous space to which we will apply \Cref{prop:nonabelian-R} to settle the key $E_{7,1}^{66}$ case of \Cref{thm:MA}. For this, $E_{7,1}^{66}$ will only matter through representation-theoretic data that it determines for its Levi in \Cref{prop:E7166-Levi}. To identify these data we will use highest weight theory.

\bpp[Irreducible representations via highest weights] \label{pp:high-wt}
Let $G$ be a split reductive group over a field $k$, let $T \subset B$ be a split maximal $k$-torus and a $k$-Borel of $G$, and let $X^*(T)^+$ be the resulting positive cone of the $T$-weights $\lambda$, with $\langle \lambda, \alpha^\vee\rangle \ge 0$ for each positive simple root $\alpha$. By \cite{Jan03}*{Part II, Propositions 2.4 and 2.6}, having fixed these data, the irreducible (and thus nonzero) $k$-representations $V$ of $G$ are classified by their highest $T$-weight: we have a bijection
\be\label{eqn:w}
\{\text{irreducible $k$-representations }V\}/\sim\, \xleftrightarrow{V \, \mapsto \text{ its highest weight } \lambda} X^*(T)^+,
\ee
where the \emph{highest weight} of $V$ is its $T$-weight $\lambda$ such that $\lambda - \lambda^\prime$ is a nonempty sum of positive simple roots for any $T$-weight $\lambda^\prime \neq \lambda$ in $V$. By \emph{loc.~cit.},~the highest weight $\lambda$ exists and its weight space $V_\lambda$ is a line equal to $V^{\mathscr R_{\mathrm u}(B)}$. Thus, the weight space $V_{w\lambda}$ is a line for every $w$ in the Weyl~group~$W$~of~$G$. 

By \cite{Jan03}*{Part II, Corollary 2.9}, the irreducible $k$-representations of $G$ stay irreducible over every field extension, so the bijection \eqref{eqn:w} commutes with base change. Thus, we have the Schur lemma $\mathrm{End}_G(V) \cong k$, and the center $Z_G$ acts on $V$ via $\lambda|_{Z_G}$. In contrast, irreducible representations in characteristic $0$ may become reducible when reduced modulo $p$, as seen already for the symmetric square of the standard representation of $\mathrm{SL}_2$. As we now explain, this does not happen for those irreducible representations $V$ whose highest weight $\lambda$ is \emph{minuscule} in the sense that $\langle \lambda, \alpha^\vee\rangle \in \{0, \pm 1\}$ for every root $\alpha$.\footnote{\label{foot:minuscule}Conversely, as noticed by Gross, except for essentially a single case consisting of the adjoint representation of $E_8$, the irreducibility of all the reductions modulo $p$ implies that the highest weight is minuscule, see \cite{GGN17}*{Theorem~1.1}.} This condition ensures the translates of $V_\lambda$ by the $\mathrm{SL}_2$ that corresponds to the coroot $\alpha^\vee$ generate either the trivial or the standard $2$-dimensional representation of $\mathrm{SL}_2$. These $\mathrm{SL}_2$’s, and so also $G$ that they generate with $T$, then preserve the subspace $\sum_{w \in W} V_{w\lambda}$. Thus, $V$ is a direct sum of lines of the form $V_{w\lambda}$, whose formation, of course, commutes with reduction modulo $p$.

When $G$ is semisimple, simply connected, the positive simple coroots $\alpha^\vee$ form a basis for the cocharacter lattice $X_*(T)$. The elements of the corresponding dual basis for $X^*(T)$ are the \emph{fundamental weights} $\omega$ of $G$, they correspond to the nodes of the Dynkin diagram of $G$. The corresponding irreducible $k$-representations are the fundamental representations of $G$ and, since $X^*(T)^+ \cong \bigoplus_\omega \mathbb Z_{\ge 0} \omega$, every irreducible $k$-representation $V$ of $G$ is a subquotient of a tensor product of the fundamental ones. The classification of those fundamental weights that are minuscule is recorded in \cite{GGN17}*{Table 1}. 
\epp

\beg\label{eg:half-spin}
In \S\ref{pp:high-wt}, suppose that $G$ is a semisimple, simply connected, absolutely indecomposable group of type $D_n$. By \cite{BouLGLA}*{chapitre VI, planche IV, (IV)}, in the highest root expression via simple roots, the two fork node simple roots contribute with coefficient $1$. Thus, since, by \cite{BouLGLA}*{chapitre~VI, section 1, numéro 1, lemme 2}, the map $\alpha \mapsto \alpha^\vee$ is additive for the roots of any simply-laced root system, such as $D_n$, \cite{BouLGLA}*{chapitre VI, section 1, numéro 8, proposition 25} implies that the fork node fundamental weights are minuscule. By \cite{FH91}*{Proposition 20.15} (with \cite{Mil17}*{Theorem 22.53}), \cite{BouLGLA}*{chapitre VI, planche IV, (VI)}, and \S\ref{pp:high-wt}, the corresponding fundamental representations of $G$ are the two half-spin representations (constructed over $\mathbb Z$ in \cite{Con14}*{Remarks C.2.9 and C.4.11}), their dimensions are $2^{n - 1}$. By \cite{BouLGLA}*{chapitre VI, planche IV, (XI)} and the description of minuscule representations reviewed in \S\ref{pp:high-wt}, these half-spin representations are self-dual when $n$ is even.
\eeg

\bpp[Tits algebras] \label{pp:TA} 
For a possibly nonsplit semisimple group $G$ over a field $k$ and a maximal $k$-torus $T \subset G$, let $\lambda$ be a dominant weight of $T_{k^s}$ for some Borel of $G_{k^s}$ containing $T_{k^s}$. By the uniqueness of the highest weight and a result of Tits \cite{Tit71}*{théorème 3.3}, $\lambda$ is fixed by the $*$-action of $\mathrm{Gal}(k^s/k)$ on $X^*(T_{k^s})^+$ (described in \cite{KMRT98}*{beginning of Section 27.A}) if and only if the $k^s$-representation $V_\lambda$ of $G_{k^s}$ of highest weight $\lambda$ descends to a representation of $G$ on a finite dimensional vector space over a central division algebra $D_\lambda$ over $k$. In this case, by the cited result of Tits, this Tits algebra $D_\lambda$ of a $*$-invariant $\lambda$ is uniquely determined, and it is split precisely when $V_\lambda$ descends to a $k$-representation of $G$, which happens when $\lambda$ lies in the root lattice or when $G$ is quasi-split. The uniqueness forces the formation of $D_\lambda$ to commute with the restriction to the simply connected cover of $G$ (with $T$ replaced by its preimage there). By \cite{Tit71}*{corollaire 3.5}, the map $\lambda \mapsto [D_\lambda] \in \mathrm{Br}(k)$ is additive, in fact, \cite{Tit71}*{section 4} (alternatively, \cite{KMRT98}*{Proposition 31.7}) ensures that 
\be \label{eqn:TF}
[D_\lambda] = H^2(k, \lambda|_{Z_G})(t_G) \quad \text{with the Tits class} \quad t_G \in H^2(k, Z_G)
\ee
defined analogously to \S\ref{pp:R-center}~\ref{RCi}. For the derived group $H$ of a $k$-Levi $T \subset L$ of a $k$-parabolic $P \subset G$ and any $*$-invariant $\lambda$, as we now argue, \eqref{eqn:TF} implies that, similarly to \cite{PS11}*{Theorem 2 (1)},
\be \label{eqn:DL}
D_\lambda \simeq D_{\lambda|_{(T \cap H)_{k^s}}}.
\ee
Indeed, by \cite{SGA3IIInew}*{exposé XXIV, 3.4 (iv); exposé XXVI, lemme 3.8}, the unique quasi-split inner twist of $G$ has a parabolic $k$-subgroup of the same geometric type as $P$. Thus, by \cite{Ser02}*{Chapter I, Proposition~37}, this inner twist is realized by a cocycle valued in $P/Z_G$, so, by \cite{SGA3IIInew}*{exposé XXVI, corollaire 2.3}, even by a cocycle valued in $L/Z_G$. Twisting by it also makes $H$ quasi-split, so we get that the Tits classes $t_{H}$ and $t_{G}$ have the same image in $H^2(k, Z_L)$. The sought $D_\lambda \simeq D_{\lambda|_{(T \cap H)_{k^s}}}$ now follows by applying $H^2(k, \lambda|_{Z_L})$ to this common~image. 
\epp

The following well-known lemma is the final preparation for extracting the promised data from $E_{7,1}^{66}$.

\blem\label{lem:Levi-sc}
For a semisimple, simply connected group $G$ over a field $k$, the derived subgroup of a Levi subgroup of a parabolic subgroup of $G$ is again simply connected.
\elem

\bpf
This is well known, see, for instance, \cite{MT11}*{Proposition 12.14} or \cite{Con14}*{Remark 5.4.7}.
\epf

\bprop\label{prop:E7166-Levi} Let $G$ be a semisimple, simply connected, absolutely indecomposable group of Tits index $E_{7,\, 1}^{66}$ over a field $k$, let $\iota\colon \mathbb G_m \simeq S \hookrightarrow G$ be a maximal split $k$-torus, and let $L \colonequals Z_G(S)$ be the corresponding Levi subgroup of a minimal parabolic $k$-subgroup $P \subset G$.
\benum
\m \label{m:E7-a}
The derived group $H$ of $L$ is semisimple, simply connected, absolutely indecomposable~of~type~$D_6$.

\m \label{m:E7-b}
One of the half-spin representations of $H_{\overline k}\simeq\mathrm{Spin}_{12}$ descends to $k$, and $L \cong (H \times S)/\mu_2$ where $\mu_2 \cong H \cap S$ is the kernel of the other half-spin representation of $H_{\overline k}$. 

\eenum

\eprop

\bpf
We fix a maximal $k$-torus $T \subset G$ containing our $S$ and a Borel $k^s$-subgroup $B \subset P_{k^s} \subset G_{k^s}$ containing $T_{k^s}$, and, in Bourbaki numbering, consider the resulting Tits index
\be\ba\label{eqn:E7}
%\begin{comment}
\begin{tikzpicture}[
  x=0.8cm,
  y=0.8cm,
  line width=0.5pt,
  font=\Large
]
  \node[anchor=east, inner sep=0pt] at (-1.05,0) {\scalebox{0.8}{$E_{7,1}^{66}$}};

  \draw (0,0) -- (5,0);
  \draw (3,0) -- (3,1);

  \foreach \x/\i in {0/7,1/6,2/5,3/4,4/3,5/1} {
    \fill (\x,0) circle[radius=2.1pt];
    \node[below=6pt] at (\x,0) {\scalebox{0.8}{$\alpha_{\i}$}};
  }

  \fill (3,1) circle[radius=2.1pt];
  \node[right=6pt] at (3,1) {\scalebox{0.8}{$\alpha_2$}};

  \draw (5,0) circle[radius=0.2cm];
\end{tikzpicture}
%\end{comment}
\ea \ee
of $G$,  see \cite{Tit66}*{page 59} or \cite{PS11}*{page 215}. Thus,  $\alpha_1$ is the unique simple root $\alpha_i$ that does not vanish on $S$, the $\alpha_2, \dotsc, \alpha_7$ span a subsystem of type $D_6$, and \Cref{lem:Levi-sc} already implies \ref{m:E7-a}.

In terms of the corresponding fundamental weights $\omega_i$, the definition of the Dynkin diagram \cite{BouLGLA}*{chapitre VI, section 4, numéro 2} gives $\alpha_1 = 2\omega_1 - \omega_3$ (see also \S\ref{pp:high-wt}). Moreover, the Dynkin diagram of $E_7$ has no nontrivial automorphisms, so all the $\alpha_i$ and $\omega_i$ are $*$-invariant. Thus, by \S\ref{pp:TA}, specifically, since $D_{\alpha_1}$ splits, $D_{\alpha_1} \simeq D_{\alpha_1|_{(T \cap H)_{k^s}}}$,  and $\omega_1|_{(T \cap H)_{k^s}}$ vanishes, \eqref{eqn:TF} shows that the Tits algebra of $H$ of the fundamental weight of the $\alpha_3$ fork node of $D_6$ in \eqref{eqn:E7} splits. \Cref{eg:half-spin} and \S\ref{pp:TA} then ensure that the corresponding half-spin representation of $H_{\overline k}$ descends to $k$.  

By considering the $\alpha_i$, we see that $S$ is the maximal central torus of $L$ (this is also a general fact because $G$ is of inner type, see \cite{CP98}*{Section 2}). Thus, it remains to argue that the central $H \cap S \subset H$ is the $\mu_2$ that is the kernel of the other half-spin representation of $H_{\overline k}$. 

Since the inclusion $\iota \colon \mathbb G_m \simeq S \hookrightarrow G$ is a primitive generator of the annihilator of the $\alpha_2, \dotsc, \alpha_7$ in $X_*(T_{k^s})$, the definition and the shape of the extended Dynkin diagram of $E_7$ given in \cite{BouLGLA}*{chapitre~VI, section 4, numéro~3; chapitre VI, planche VI, (IV)} show that the highest root is $\widetilde \alpha = \omega_1$, that $\iota = \pm\widetilde \alpha^\vee$, and that $\langle \alpha_1, \iota\rangle = \pm 1$ (see \cite{BouLGLA}*{chapitre VI, section 1, no. 8, proposition 25} for the properties of the highest root). It follows that $\langle \omega_1, \iota\rangle = \pm 2$, so, since $\omega_1$ factors through $L/H$ and induces a primitive character of this torus, $H \cap S \cong \mu_2$ and it remains to identify this central $\mu_2$.

By \cite{BouLGLA}*{chapitre~VI, planche IV, (VIII)}, $Z_{H_{\overline k}} \simeq \mu_2 \times \mu_2$ with a character basis given by the restrictions of the two fork node fundamental weights. Thus, we only need to argue that 
\[
\langle \omega_3, \iota\rangle \equiv  1 \bmod 2 \quad \text{and} \quad \langle \omega_2, \iota\rangle \equiv  0 \bmod 2.
\]

The relation $\alpha_1 = 2\omega_1 - \omega_3$ now shows that $\langle \omega_3, \iota\rangle \equiv  1 \bmod 2$. For the remaining $\langle \omega_2, \iota\rangle$, we will also use the following equalities supplied by the Dynkin diagram \eqref{eqn:E7}:
\[
 \alpha_2 = 2\omega_2 - \omega_4 \quad \text{and} \quad \alpha_3 = 2\omega_3 - \omega_1 - \omega_4.
 \]
 As far as pairing with $\iota$ goes,  we may ignore the $\alpha_{\ge 2}$, so these equalities give
 \[
 \langle 2\omega_2, \iota\rangle = \langle 2\omega_3 - \omega_1, \iota\rangle = \langle 3\omega_1 - 2\alpha_1, \iota\rangle = \pm\langle 3\widetilde \alpha, \widetilde \alpha^\vee\rangle \mp 2 = \pm 4.
 \]
 By dividing both sides by $2$, we obtain the desired $\langle \omega_2, \iota\rangle \equiv 0 \bmod 2$.
\epf

\brems
\remi
We have argued  with Tits algebras, however, one can also build the half-spin $k$-representation of $H$ concretely as the $S$-weight $1$ piece of the Lie algebra $\mathfrak g$ of $G$. By analyzing the $S$-weight grading of $\mathfrak g$ and using the $G$-equivariant Lie bracket $[\, , ] \colon \mathfrak g \times_k \mathfrak g \rightarrow \mathfrak g$ one can then also construct its various invariant functions, see \cite{Hel12}*{especially, top of page 118 and~Proposition~38}.

\remi
Writing $H \simeq \mathrm{Spin}(A, \sigma, f)$ as in \S\ref{pp:conv}~\ref{m:conv-ii}, one may also analyze the central simple algebra $A$ in detail, see \cite{GQM09}*{Theorem 3.1 and Remark 3.3}.
\erems

The Tits algebra technique used in the proof of \Cref{prop:E7166-Levi} combines with various results in the literature to give the following nonexistence result for various exceptional Tits indices. For this, we focus on the needs of \Cref{thm:MA} and do not aim to be exhaustive: some other exceptional Tits indices also cannot occur when $\mathrm{dim}_2^{\mathrm{sep}}(k) \le 2$ (or when $\mathrm{dim}_3^{\mathrm{sep}}(k) \le 2$), but we will not need this.

\bprop
\label{prop:no} For a field $k$ with $\mathrm{dim}_2^{\mathrm{sep}}(k) \le 2$, there exists no semisimple, absolutely indecomposable $k$-group $G$ of Tits index either ${}^2E_{6,1}^{35}$, $E_{7, 1}^{48}$, $E_{8, 1}^{91}$, or $E_{8, 1}^{133}$.
\eprop
 
 \bpf
 We lose no generality by assuming that $G$ is simply connected. When the Tits index of $G$ is
 \be\label{eqn:E6135}\ba 
 \begin{tikzpicture}[
  x=0.8cm,
  y=0.8cm,
  line width=0.5pt,
  font=\Large
]
  \node[anchor=east, inner sep=0pt] at (-1.05,0) {\scalebox{0.8}{${}^2E_{6,1}^{35}$}};

  \draw (0,0) -- (4,0);
  \draw (2,0) -- (2,1);

  \foreach \x/\i in {0/6,1/5,2/4,3/3,4/1} {
    \fill (\x,0) circle[radius=2.1pt];
    \node[below=6pt] at (\x,0) {\scalebox{0.8}{$\ $}};
  }

  \fill (2,1) circle[radius=2.1pt];
  \draw (2,1) circle[radius=0.2cm];
  \node[right=6pt] at (2,1) {\scalebox{0.8}{$\ $}};
\end{tikzpicture}
\ea
\ee
we use the argument in the second paragraph of the proof of \Cref{prop:E7166-Levi} based on \eqref{eqn:DL} to get that the Tits algebra of the fundamental weight of the middle node in the $A_5$ Dynkin diagram of the anisotropic kernel splits. This anisotropic kernel is of outer type $A_5$, so, by \S\ref{pp:conv}~\ref{m:conv-i}, it is an $\mathrm{SU}(A, \sigma)$ for a separable quadratic extension $k^\prime/k$, a central simple algebra $A$ over $k^\prime$ of degree~$6$, and a $k^\prime/k$-unitary involution $\sigma$ on $A$. By \cite{KMRT98}*{Type $A_n$ case of Section 27.B}, the Tits algebra above is the discriminant algebra of $(A, \sigma)$ and, by \cite{KMRT98}*{Definition 10.4 (a) and Proposition~10.30}, its restriction to $k^\prime$ has Brauer class $3[A]$. In particular, this discriminant algebra splits and the index of $A$ divides $3$. Thus, since $\mathrm{dim}_2^{\mathrm{sep}}(k) \le 2$, \cite{BFT07}*{Theorem 3.2} (which we have used at the end of the proof of \Cref{thm:Rost-A2}) applies and contradicts the anisotropicity of $\mathrm{SU}(A, \sigma)$. 
 
In the case when the Tits index of $G$ is
\be\label{eqn:E7148}\ba
\begin{tikzpicture}[x=0.8cm, y=0.8cm, line width=0.5pt, font=\Large]
  \node[anchor=east, inner sep=0pt] at (-1.05,0) {\scalebox{0.8}{$E_{7,1}^{48}$}};

  \draw (0,0) -- (5,0);
  \draw (3,0) -- (3,1);

  \foreach \x/\i in {0/7,1/6,2/5,3/4,4/3,5/1} {
    \fill (\x,0) circle[radius=2.1pt];
    \node[below=6pt] at (\x,0) {\scalebox{0.8}{$\ $}};
  }

  \fill (3,1) circle[radius=2.1pt];
  \node[right=6pt] at (3,1) {\scalebox{0.8}{$\ $}};

  \draw (1,0) circle[radius=0.2cm];
\end{tikzpicture}
\ea\ee
the anisotropic kernel is of inner type $A_1 \times D_5$ and, by its anisotropy, the $A_1$ factor is $\mathrm{SL}_1(D)$ for a division algebra $D$ of index $2$ (which is necessarily a quaternion algebra, see \cite{EKM08}*{Proposition~98.9}). By \Cref{eg:half-spin} and the argument in the second paragraph of the proof of \Cref{prop:E7166-Levi} based on \eqref{eqn:DL}, $D$ is also the Tits algebra of the fundamental weight of one of the half-spin fork nodes of the $D_5$ factor. Thus, by \S\ref{pp:TA} and \cite{BouLGLA}*{chapitre VI, planche IV, (VIII)}, the Tits algebra of the fundamental weight of the rightmost node of the Dynkin diagram of the $D_5$ factor in \eqref{eqn:E7148} splits. By \cite{KMRT98}*{Type~$D_n$ case of Section 27.B} and \S\ref{pp:conv}~\ref{m:conv-ii}, this means that our $D_5$ factor of inner type is of the form $\mathrm{Spin}(q)$ for a nondegenerate quadratic form $q$ on $k^{10}$ of trivial discriminant. By \emph{loc.~cit.},~\cite{KMRT98}*{Proposition~8.8 and before Theorem 25.12}, and \cite{EKM08}*{Remark 13.9}, the Brauer class of the Clifford algebra of $q$ is represented by $D$. In particular, letting $\mathbb H$ be the hyperbolic plane and $\mathrm{Nrd}_D$ be the reduced norm form of $D$ over $k$, by \cite{EKM08}*{Example 9.6, Corollary 13.3, Lemma 13.4, Example 13.5, Example~14.1, and Lemma 14.2}, the nondegenerate quadratic forms
\[
q \quad \text{and} \quad \mathrm{Nrd}_D \perp\mathbb H \perp \mathbb H \perp \mathbb H
\]
have the same dimension, discriminant, and Clifford invariant. Thus, since $\mathrm{dim}_2^{\mathrm{sep}}(k) \le 2$, \cite{BFT07}*{Theorem 2.1 and Remark 2.2} gives $q \simeq \mathrm{Nrd}_D \perp\mathbb H \perp \mathbb H \perp \mathbb H$, which contradicts the anisotropicity.

In the case when the Tits index of $G$ is
\be\label{eqn:E8191}\ba
\begin{tikzpicture}[x=0.8cm, y=0.8cm, line width=0.5pt, font=\Large]
  \node[anchor=east, inner sep=0pt] at (-1.05,0) {\scalebox{0.8}{$E_{8,1}^{91}$}};

  \draw (0,0) -- (6,0);
  \draw (4,0) -- (4,1);

  \foreach \x/\i in {0/8,1/7,2/6,3/5,4/4,5/3,6/1} {
    \fill (\x,0) circle[radius=2.1pt];
    \node[below=6pt] at (\x,0) {\scalebox{0.8}{$\ $}};
  }

  \fill (4,1) circle[radius=2.1pt];
  \node[right=6pt] at (4,1) {\scalebox{0.8}{$\ $}};

  \draw (6,0) circle[radius=0.2cm];
\end{tikzpicture}
\ea\ee
the anisotropic kernel is of type $D_7$, so \cite{GG09}*{Lemma 9.1} forces the order of the Rost invariant of the torsor of the split form of $E_8$ that gives $G$ to be $2$. This is impossible because the condition $\mathrm{dim}_2^{\mathrm{sep}}(k) \le 2$ kills the target of the $2$-primary part of the Rost invariant, see \S\ref{pp:Rost} and \S\ref{pp:sep-dim}.

Since $E_8$ has a trivial center (see \cite{BouLGLA}*{chapitre VI, planche VII, (VIII)}), in the case of 
\be\label{eqn:E81133}\ba
\begin{tikzpicture}[x=0.8cm, y=0.8cm, line width=0.5pt, font=\Large]
  \node[anchor=east, inner sep=0pt] at (-1.05,0) {\scalebox{0.8}{$E_{8,1}^{133}$}};

  \draw (0,0) -- (6,0);
  \draw (4,0) -- (4,1);

  \foreach \x/\i in {0/8,1/7,2/6,3/5,4/4,5/3,6/1} {
    \fill (\x,0) circle[radius=2.1pt];
    \node[below=6pt] at (\x,0) {\scalebox{0.8}{$\ $}};
  }

  \fill (4,1) circle[radius=2.1pt];
  \node[right=6pt] at (4,1) {\scalebox{0.8}{$\ $}};

  \draw (0,0) circle[radius=0.2cm];
\end{tikzpicture}
\ea\ee
\Cref{lem:Levi-sc}, the same cocycle argument as at the end of \S\ref{pp:TA}, and the vanishing of $H^1(k, \mathbb G_m)$ show that the anisotropic kernel of type $E_7$ is purely inner. Thus, \cite{GPS16}*{Theorem 10.10} (with \cite{GMS03}*{(10.16) on page 135}) shows that the vanishing of the $2$-primary part of the Rost invariant, which is again imposed by $\mathrm{dim}_2^{\mathrm{sep}}(k) \le 2$, contradicts the anisotropicity.
\epf

\brem\label{rem:no} 
In contrast, $E_{7, 1}^{66}$ cannot be excluded by cohomological dimension $\le 2$. To see this, first of all, \cite{BH07}*{Theorem 1.3} (which is based on earlier constructions of Merkurjev) supplies a field $k$ of characteristic $0$ with $\mathrm{cd}_2(k) = 2$ and $\mathrm{cd}_p(k) = 0$ for every odd prime $p$, and a central division algebra $D$ over $k$ of degree $4$ that is a tensor product of two quaternion algebras over $k$. In fact, their $k$ is such that \cite{QMT20}*{Corollary 2.1 (especially, its final aspect) and before Section 1 on page 171} applies to $D$, so \emph{loc.~cit.}~(with \cite{KMRT98}*{Type $D_n$ case of Section 27.B}) produces a semisimple, simply connected $k$-group $H \colonequals \mathrm{Spin}(\mathrm{Mat}_{3\times 3}(D), \sigma)$ of inner type $D_6$ (see \S\ref{pp:conv}~\ref{m:conv-ii}) each of whose fork node fundamental weights has a Tits algebra that is either isomorphic to $D$ or splits. By \S\ref{pp:TA} and \cite{KMRT98}*{Type $D_n$ case of Section 27.B}, the sum of the Brauer classes of these two Tits algebras is the class of $D$, so precisely one of them splits. At this point  \cite{PS11}*{Theorem~2~(2)} applied to $H$ produces a semisimple, simply connected $k$-group $G$ of type $E_7$ equipped with a parabolic $k$-subgroup for which the derived subgroup of its Levi is $H$. Since the index of $D$ is $4$, \cite{DCG12}*{Section~4.5.1, especially, Table 8} implies that the Tits index of $G$ over $k$ must be $E_{7, 1}^{66}$. 
\erem

\section{Unique decompositions into sums of pure spinors give the rationality of $L/M$} \label{sec:half-spin}
 
 We analyze the geometry of $32$-dimensional half-spin representations over an arbitrary field $k$, the one from \Cref{prop:E7166-Levi}~\ref{m:E7-b} being the key example. Starting from any $k$-form $V$ of a $32$-dimensional half-spin representation, in \Cref{cor:LM-ratl} we build a certain homogeneous space, show that it is a rational variety, and then in \Cref{cor:GRtrivial} use it to obtain the $E_{7,1}^{66}$ case of \Cref{thm:MA}. The geometric engine behind the rationality is \Cref{thm:dt}, which gives unique decompositions of generic spinors into sums of two opposite pure spinors. These decompositions are special to dimension $32$, although with some restrictions on the characteristic and in slightly broader settings similar phenomena have been observed by Helenius in \cite{Hel12}*{Lemma 24} and by Slupinski–Stanton in~\cite{SS15}.
 
 \bpp[A $k$-descent of a $32$-dimensional half-spin representation] \label{pp:setup}
 We fix
 \bitem
 \m
 a field $k$;
 
 \m
 a semisimple, simply connected, absolutely indecomposable $k$-group $H$ of type $D_6$;
 
 \m
 a $k$-representation $V$ of $H$ such that $V \otimes_k \overline k$ is a half-spin representation of $H_{\overline k} \simeq \mathrm{Spin}_{12}$.
 \eitem
 By \Cref{eg:half-spin}, our $V$ has $k$-dimension $32$ but we cannot specify “which” of the two half-spin representations of $H_{\overline k}$ it descends because the nontrivial outer automorphism of $H_{\overline k}$ interchanges them. The “other” half-spin of $H_{\overline k}$ need not descend to $k$; relatedly, $V$ may not have a quadratic space Clifford model over $k$. Nevertheless, $V$ becomes isomorphic to a half-spin representation defined over $\mathbb Z$ over any splitting field of $H$ (see \S\ref{pp:high-wt}, \Cref{eg:half-spin}, and use Schur’s lemma).
 
 The existence of $V$ constrains $H$, for instance, by \S\ref{pp:high-wt}, the $*$-action must fix the fork nodes of the Dynkin diagram, so $H$ is of inner type. By \S\ref{pp:high-wt}, \Cref{eg:half-spin}, and \cite{BouLGLA}*{chapitre VI, planche~IV, (VIII)}, the kernel $\mu_2$ of the other half-spin representation of $H_{\overline{k}}$ descends to $k$ (even if the representation itself need not) and acts on $V$ via its unique nontrivial character. Thus, letting $\mathbb G_m$ act by scaling, the $H$-action on $V$ extends to the action of 
 \be\label{eqn:L}
 L \colonequals (H \times_k \mathbb G_m)/\mu_2,
 \ee
 where $\mu_2$ is embedded diagonally using the kernel of the other half-spin representation of $H_{\overline k}$. 
 \epp

\bpp[The variety of pure spinors] \label{pp:pure}
In the setting of \S\ref{pp:setup}, the $L$-action on $V$ is linear, so it induces an $L$-action on $\mathbb P(V)$. As we now argue (inspired by \cite{BGL14}*{Proposition 5.1}), there is a unique closed $H$-orbit $S \subset \mathbb P(V)$, and its formation commutes with base change. For this, by Galois descent, we may assume that $H$ is split, fix a split $k$-torus $T$ in a $k$-Borel $B \subset H$, and then argue as follows.
\benum
\m\label{m:p-a} 
$S$ exists: by \S\ref{pp:high-wt}, $V$ has a unique stable $B$-line---its highest $T$-weight line $V_\lambda$---and, by the minusculeness of half-spins (see \Cref{eg:half-spin}), the stabilizer of $V_\lambda$ is the parabolic $B \subset P \subset H$ corresponding to the complement of the half-spin node of $V$ in the Dynkin diagram of $H$;

\m\label{m:p-b} 
$S$ is unique: by the Borel fixed point theorem \cite{Bor91}*{Theorem 10.4}, every $B$-stable closed subvariety of $\mathbb P(V)$ has a $\overline{k}$-point fixed by $B$, which must then be the class of the unique stable $B$-line, whose $H$-orbit is the one from \ref{m:p-a}.
\eenum
By \ref{m:p-a}, in fact, $S \subset \mathbb P(V)$ is the smooth $k$-variety parametrizing the parabolic subgroups of $H$ of the same type as $P$ and, by \cite{BouLGLA}*{chapitre VI, planche I, (I), planche IV, (I)}, its dimension is $15$.

The \emph{variety of pure spinors} $\mathscr S \subset V \setminus \{0\}$ is the preimage of $S \subset \mathbb P(V)$, it is $L$-stable of dimension $16$. We say that a $v \in V \setminus \{0\}$ is a \emph{pure spinor} if it lies in $\mathscr S$. This agrees with the classical definition of Chevalley \cite{Che97}*{Chapter III, Section 3.1, page 136 and what follows}: in a Clifford model for the base change $V_{\overline{k}}$, the pure spinors are precisely the nonzero vectors $v$ whose annihilator subspaces in the quadratic space have the largest possible dimension (this dimension condition amounts to an upper bound on the rank of a linear map that depends on $v$, so it defines a closed subvariety in $\mathbb P(V)$, which is a single $H$-orbit by \emph{loc.~cit.},~and so is the unique closed $H$-orbit $S$ from \ref{m:p-a}). 
The $L$-action on $\mathscr S$ is transitive and the $H$-stabilizer of any pure spinor line in $V$ is a parabolic of the same type as $P$ in \ref{m:p-a}. By \Cref{lem:Levi-sc}, the derived subgroup of a Levi of any such stabilizer is semisimple, simply connected, of type $A_5$. Thus, if $H$ is anisotropic, then $\mathscr S(k) = \emptyset$, although \Cref{thm:dt} will show that $\mathscr S$ has many points over separable, at most quadratic extensions of $k$.

Two pure spinors are \emph{opposite} if their line stabilizers in $H$ are opposite as parabolic subgroups in the sense of \cite{SGA3IIInew}*{exposé XXVI, définition 4.3.3}. Pairs of such are parametrized by a subfunctor
\[
\mathscr O \subset \mathscr S \times_k \mathscr S
\]
which, by its definition, is the preimage in $\mathscr S \times_k \mathscr S$ of its image in $S \times_k S$. By \cite{SGA3IIInew}*{exposé~XXVI, corollaire 4.3.4 (ii)}, $\mathscr O$ is an $L$-stable dense open of $\mathscr S \times_k \mathscr S$ (of dimension $32$). 
\epp

By the following lemma, once a rational pure spinor exists, the varieties $S$, $\mathscr S$, and $\mathscr O$  become~rational.

\blem\label{lem:S-ratl}
If $S(k) \neq \emptyset$, equivalently, if there is an $s \in \mathscr S(k)$, then $S$, $\mathscr S$, and $\mathscr O$ are rational over $k$\ucolon then for any $k$-parabolic $P^- \subset H$ opposite to the $k$-parabolic stabilizing the line $ks \subset V$, the action of $\mathscr R_{\mathrm u}(P^-)$ on the class of $s$ \up{resp., of $\mathbb G_{m,\, k} \times \mathscr R_{\mathrm u}(P^-)$ on $s$} gives an open immersion
\[
\mathscr R_{\mathrm{u}}(P^-) \subset S \quad \text{{\upshape(}resp.,} \quad\mathbb G_m \times_k \mathscr R_{\mathrm{u}}(P^-) \subset \mathscr S).
\]
\elem 

\bpf
The map $\mathscr S \rightarrow S$ is a $\mathbb G_m$-torsor, so $S(k) \neq \emptyset$ if and only if $\mathscr S(k) \neq \emptyset$. Moreover, by \S\ref{pp:pure} and \cite{SGA3IIInew}*{exposé XXVI, corollaires~4.3.5~(i),~2.5}, a $P^-$ exists and $\mathscr R_{\mathrm u}(P^-) \simeq \mathbb A^n_k$. It remains to note that the open immersions follow from the description of $S$ in \S\ref{pp:pure}, from \cite{SGA3IIInew}*{exposé XXVI, théorème~4.3.2~(vi${}^\prime$), définition 4.3.3}, and from $\mathscr S \rightarrow S$ being a $\mathbb G_m$-torsor.
\epf

The following result, which is critical for our purposes, generalizes \cite{Hit03}*{Proposition 2} to arbitrary characteristic and simultaneously gives a synthetic, computation-free reproof of \emph{loc.~cit.}~itself.

\bthm \label{thm:dt} For $k$, $H$, and $V$ as in \uS\uref{pp:setup}, there is a unique open $L$-orbit $\mathcal O \subset V$, the $L$-morphism $\mathscr O \rightarrow V$ given by summing two pure spinors is étale and induces a $\mathbb Z/2\mathbb Z$-cover $\mathscr O \rightarrow \mathcal O$ \up{whose nontrivial deck transformation interchanges the summands}, and the $L$-action on $\mathscr O$ is transitive. Moreover, for any expression of a  $v \in \mathcal O(\overline{k})$ as a sum of two pure spinors, the latter are opposite. 
\ethm

\bpf
The sum map $\mathscr S \times_k \mathscr S \rightarrow V$ is $L$-equivariant and its source and target are geometrically irreducible of dimension $32$ (see \S\ref{pp:setup}--\ref{pp:pure}). Thus, the claim about $v$ follows from the rest: the $L$-stable boundary $(\mathscr S \times_k \mathscr S) \setminus \mathscr O$ has dimension $\le 31$, so its image cannot meet an open $L$-orbit $\mathcal O \subset V$.

For the rest, it is enough to show that $L$ acts transitively on $\mathscr O$ and that $\mathscr O \rightarrow V$ is étale, of fibral degrees $\le 2$: its image will then be the (unique) open $L$-orbit $\mathcal O \subset V$ and, since interchange of opposite pure spinors has no fixed points (a proper parabolic of $H$ is never opposite to itself), the induced map $\mathscr O \rightarrow \mathcal O$ will be a $\mathbb Z/2\mathbb Z$-cover. The transitivity, étaleness, and the fibral degree claims may be checked over $\overline k$, so for the rest of the argument we assume that $k$ is algebraically closed.

For the étaleness, let $(s, s^\prime) \in \mathscr O(k)$ be a pair of opposite pure spinors and let $P_s, P_{s'} \subset H$ be the opposite parabolics over $k$ that are the $H$-stabilizers of the lines $ks, ks^\prime \subset V$ (see \S\ref{pp:pure}). The intersection $P_s \cap P_{s^\prime}$ is the common Levi of $P_s$ and $P_{s^\prime}$, its derived group $M$ is semisimple, simply connected of type $A_5$ (see \S\ref{pp:pure}). Since $M$ has no nontrivial characters, it fixes both $s$ and $s^\prime$, and not only their classes in $S$. By \Cref{lem:S-ratl}, the $(\mathbb G_m \times H)$-action on $s$ and $s^\prime$ induces the maps
\be\label{eqn:sS-ratl}
\mathbb G_m \times_k \mathscr R_{\mathrm{u}}(P_{s^\prime}) \subset \mathscr S \quad \text{and} \quad \mathbb G_m \times_k\mathscr R_{\mathrm{u}}(P_{s}) \subset \mathscr S,
\ee
respectively,  that are $M$-equivariant open immersions. Therefore, $M$-equivariantly, we get
\be\label{eqn:Lie-Lie}
\mathrm{Tan}_s(\mathscr S) \cong ks \oplus \mathrm{Lie}(\mathscr R_{\mathrm{u}}(P_{s^\prime}))\quad \text{and} \quad \mathrm{Tan}_{s^\prime}(\mathscr S) \cong ks^\prime \oplus \mathrm{Lie}(\mathscr R_{\mathrm{u}}(P_{s})),
\ee
so that also
\be\label{eqn:Lievv}
\mathrm{Tan}_{(s, s^\prime)}(\mathscr O) \cong ks \oplus \mathrm{Lie}(\mathscr R_{\mathrm{u}}(P_{s^\prime}))  \oplus \mathrm{Lie}(\mathscr R_{\mathrm{u}}(P_{s})) \oplus ks^\prime.
\ee
For the étaleness of $\mathscr O \rightarrow V$ at $(s, s^\prime)$, we need to show that this direct sum maps isomorphically onto $V \cong \mathrm{Tan}_{s + s^\prime}(V)$. Certainly, $ks \oplus ks^\prime \hookrightarrow V$ because the classes of $s$ and $s^\prime$ in $\mathbb P(V)$ are distinct. To separate them and the other two summands in \eqref{eqn:Lievv} (each of which is mapped injectively into $V$ thanks to \eqref{eqn:sS-ratl}) we will use the action of the center $Z_M$. 

For this, we fix a maximal $k$-torus $T \subset P_s \cap P_{s^\prime}$ in a Borel $k$-subgroup $B \subset P_s$ of $H$ and use the resulting Dynkin diagram of $H$ with its simple roots. The $T$-weights that appear in $\mathrm{Lie}(\mathscr R_{\mathrm{u}}(P_{s}))$ are the opposites of the ones in $\mathrm{Lie}(\mathscr R_{\mathrm{u}}(P_{s'}))$, see \cite{SGA3IIInew}*{exposé~XXVI, théorème~4.3.2~(vi${}^\prime$)}. To analyze them, let $\alpha$ be the $T$-root of $H$ that corresponds to the fork node of $V$ in the $D_6$ Dynkin diagram of $H$, see \Cref{eg:half-spin}. By \S\ref{pp:pure}, removing $\alpha$ leaves us the Dynkin diagram of $M$ with respect to $T \cap M \subset B \cap M$. Its simple roots all appear in $\mathrm{Lie}(M)$, so they are trivial on $Z_M$. In particular, the $Z_M$-action on $\mathrm{Lie}(\mathscr R_{\mathrm{u}}(P_{s}))$ and $\mathrm{Lie}(\mathscr R_{\mathrm{u}}(P_{s'}))$ depends only on the coefficient of $\alpha$ in the simple root expression of each $T$-root in these Lie algebras. These coefficients are nonzero: else the corresponding root would contribute to $\mathrm{Lie}(M)$ instead. Since, by \cite{BouLGLA}*{chapitre VI, planche~IV,~(IV)}, the coefficient of $\alpha$ in the highest root of $H$ is $1$, we conclude from \cite{BouLGLA}*{chapitre~VI, section~1, no.~8, proposition 25 (i)} that it is $1$ for every root in $\mathrm{Lie}(\mathscr R_{\mathrm{u}}(P_{s}))$ and $-1$ for every root in $\mathrm{Lie}(\mathscr R_{\mathrm{u}}(P_{s'}))$. Thus, $Z_M$ acts via $\alpha|_{Z_M}$ on $\mathrm{Lie}(\mathscr R_{\mathrm{u}}(P_{s}))$ and via $\alpha^{-1}|_{Z_M}$~on~$\mathrm{Lie}(\mathscr R_{\mathrm{u}}(P_{s'}))$.

Since $M$ is simply connected, the $\beta^\vee$ for simple roots $\beta \neq \alpha$ of $H$ form a $\mathbb Z$-basis for $X_*(T \cap M)$. Thus, by the construction of the Dynkin diagram \cite{BouLGLA}*{chapitre VI, section 4, no. 2} and by \cite{BouLGLA}*{chapitre VI, section 1, numéro 10}, $\alpha^{-1}|_{T\cap M}$ is the fourth fundamental weight of $M$. By \cite{BouLGLA}*{chapitre VI, planche I, (VI), (VIII)}, we get that $\alpha|_{Z_M}$ has order $3$, so that $\alpha|_{Z_M} \neq \alpha^{-1}|_{Z_M}$.

In conclusion, the two middle summands in \eqref{eqn:Lievv} have distinct, nontrivial $Z_M$-weights, while $Z_M$ acts trivially on the two outer summands. Consequently, since the map $\mathscr O \rightarrow V$ is $Z_M$-equivariant and since we have already observed that $ks \oplus ks^\prime \subset V$, we find that
\be\label{eqn:Lie}
\mathrm{Tan}_{(s, s^\prime)}(\mathscr O) \subseteq \mathrm{Tan}_{s+ s^\prime}(V) \cong V.
\ee
This $k$-vector space inclusion is an isomorphism because, by \S\ref{pp:setup}, the $k$-dimension of $V$ is $32$, while, by \S\ref{pp:pure}, that of $\mathrm{Tan}_{(s, s^\prime)}(\mathscr O)$ is also $32$. Thus, by \cite{BLR90}*{Section 2.2, Corollary 10}, $\mathscr O \rightarrow V$ is étale.

For the claim about the fibral degree being $\le 2$, suppose that $s + s^\prime = t + t^\prime$ for another pair of opposite pure spinors $(t, t^\prime)$. The étale map $\mathscr O \rightarrow V$ is $M$-equivariant, so, since $M$ fixes $t + t^\prime$ and is connected, it must also fix $(t, t^\prime)$. The $Z_M$-action then gives $t, t^\prime \in ks \oplus ks^\prime \subset V$, so  we wish to show that the only pure spinors in $ks \oplus ks^\prime$ are the $k^\times s$ and the $k^\times s^\prime$. Since $T \subset P_s \cap P_{s^\prime}$, both $ks$ and $ks^\prime$ are $T$-stable lines, while $ks \subset V$ is the unique highest $T$-weight line thanks to our choice of $B$ (see \S\ref{pp:pure}~\ref{m:p-a}). Consequently, $T$ acts on $ks$ and on $ks^\prime$ via distinct characters. Thus, if some $as + bs^\prime$ with $ab \neq 0$ was a pure spinor, then, thanks to the $T$-action, the entire $(ks \oplus ks^\prime)\setminus \{ 0\}$ would consist of pure spinors. This would give the inclusion $ks \oplus ks^\prime \subseteq \mathrm{Tan}_s(\mathscr S)$, which contradicts \eqref{eqn:Lie-Lie}--\eqref{eqn:Lievv}. In conclusion, no $as + bs^\prime$ as above is a pure spinor, so the fibral degrees of $\mathscr O \rightarrow V$ are indeed $\le 2$. 

For the remaining transitivity of the $L$-action on $\mathscr O \subset \mathscr S \times_k \mathscr S$, we first use the transitivity of the $L$-action on $\mathscr S$ (see \S\ref{pp:pure}) to arrange the first coordinate to be the $s$ above. We then act by $\mathscr R_{\mathrm u}(P_s)$ and use \cite{SGA3IIInew}*{exposé XXVI, corollaire 4.3.5 (i)}, \S\ref{pp:pure}, and the fact that unipotent groups have no nontrivial characters to arrange the second coordinate to be $b s^\prime$ for some unit $b$. At this point, since $ks$ and $ks^\prime$ have distinct $T$-weights, we may act by $T$ to arrive at $(as, as^\prime)$ for some unit $a$. We then use the scaling action of the $\mathbb G_m$ that contributes to $L$ in \eqref{eqn:L} to make $a = 1$.
\epf

\bcor \label{cor:LM-ratl} Let $k$, $H$, and $V$ be as in \uS\uref{pp:setup}, and let $v \in \mathcal O(k)$.
\benum
\m\label{m:LM-a} The stabilizer $\mathrm{Stab}_L(v)$ is $k$-smooth, $\pi_0(\mathrm{Stab}_L(v)) \cong \mathbb Z/2\mathbb Z$, and $\mathrm{Stab}_L(v)^0$ is semisimple, simply connected of type $A_5$ and over $\overline k$ is the derived group of the intersection of the two opposite parabolic $H$-stabilizers of the lines of the two opposite pure spinor summands of $v$.

\m\label{m:LM-b} We have $L/\mathrm{Stab}_L(v) \cong \mathcal O$, while $L/\mathrm{Stab}_L(v)^0$ is either $\mathscr O$ or its quadratic twist over $\mathcal O$.

\m \label{m:LM-c}
Both $L/\mathrm{Stab}_L(v)$ and $L/\mathrm{Stab}_L(v)^0$ are rational varieties over $k$ of dimension $32$.
\eenum

\ecor

\bpf
Part \ref{m:LM-a} is geometric, so we may argue it over $\overline{k}$. By \Cref{thm:dt}, especially, by the transitivity of the $L$-action, we may then assume that $v = s + s^\prime$ with $(s, s^\prime)\in \mathscr O(k)$ as in the proof of \Cref{thm:dt}. In the notation of that proof, with the presentation $L \cong (H\times_k \mathbb G_m)/\mu_2 $ as in \eqref{eqn:L}, we have $\mathrm{Stab}_{H \times_k \mathbb G_m}(ks, ks^\prime) \cong (P_s \cap P_{s'}) \times_k \mathbb G_m$ and the common Levi $P_s \cap P_{s'}$ is an extension of $\mathbb G_m$ by the derived group $M$ that is semisimple, simply connected of type $A_5$. Thus, by \Cref{thm:dt}, to obtain \ref{m:LM-a} it suffices to argue that $M \subset H \subset H \times_k \mathbb G_m$ maps isomorphically onto $\mathrm{Stab}_L(s, s')$. 

We saw in the proof of \Cref{thm:dt} that $M \subset \mathrm{Stab}_H(s, s')$ and that, with a fixed $k$-Borel $T \subset B \subset P_s$, our $P_s \cap P_{s'}$ acts on $ks$ via the highest $T$-weight $\lambda$ in $V$ viewed as a character of $P_s \cap P_{s'}$. The opposite Borel $T \subset B^-$ is the semidirect product of the corresponding opposite Borel in $P_s \cap P_{s'}$ and $\mathscr R_{\mathrm u}(P_{s'})$, so $T \subset B^- \subset P_{s'}$. Thus, by symmetry, $P_s \cap P_{s'}$ acts on $ks'$ via the highest $T$-weight in $V$ relative to $B^-$, in other words, by the lowest $T$-weight relative to $B$. By \Cref{eg:half-spin}, as an $H$-representation $V$ is self-dual, so $P_{s} \cap P_{s'}$ acts on $ks^\prime$ via $\lambda\i$.

We conclude that $\mathrm{Stab}_{H \times_k \mathbb G_m}(s, s^\prime)$ parametrizes schematic points $(h, a)$ of $(P_s \cap P_{s'}) \times_k \mathbb G_m$ with 
\[
a \lambda(h) = 1 = a \lambda^{-1}(h), \qxq{in other words, with} \lambda(h)^2 = 1 \qxq{and} a = \lambda(h).
\]
Thus, since $\lambda$ is nontrivial on the $\mu_2 \subset H \times_k \bG_m$ in  \eqref{eqn:L} (see \S\ref{pp:setup}), we may adjust by the section $\lambda(h)$ of this $\mu_2$ to make $h$ lie in $\mathrm{Ker}(\lambda)$ and $a = 1$. Since the highest weight $\lambda$ is a fundamental weight (see \Cref{eg:half-spin}), its kernel as a character of $P_s \cap P_{s'}$ is the derived group $M$. Therefore, we get that $M \cong \mathrm{Stab}_H(s, s^\prime)$ and, since  $\mathrm{Stab}_L(s, s^\prime)$ is the image of $\mathrm{Stab}_{H \times_k \mathbb G_m}(s, s^\prime)$, that $M \xrightarrow{\sim}\mathrm{Stab}_L(s, s^\prime)$.

For \ref{m:LM-b} and \ref{m:LM-c}, we return to an arbitrary $k$. Their claims about $L/\mathrm{Stab}_L(v)$ follow from \Cref{thm:dt} (with \S\ref{pp:setup}). To study $L/\mathrm{Stab}_L(v)^0$, let $k_v$ be the quadratic étale $k$-algebra that is the coordinate ring of the $k$-fiber of $\mathscr O$ above $v$ and let $\sigma$ be the nontrivial $k$-automorphism of $k_v$ supplied by \Cref{thm:dt}, so the preimage of $v$ is a $k_v$-point $(s, \sigma(s)) \in \mathscr O(k_v) \subset \mathscr S(k_v) \times \mathscr S(k_v)$. Viewed differently, $s$ is a $k$-point of $\mathrm{Res}_{k_v/k}(\mathscr S_{k_v})$ that lies in the open $\overline {\mathscr O} \subset \mathrm{Res}_{k_v/k}(\mathscr S_{k_v})$ consisting of those $k$-scheme valued points $x$ for which $(x, \sigma(x))$ lies in the open $\mathscr O_{k_v} \subset \mathscr S_{k_v} \times_{k_v} \mathscr S_{k_v}$. This $\overline {\mathscr O}$ is simply the twist of $\mathscr O$ by the quadratic $k$-algebra $k_v$: indeed, we have $(\mathrm{Res}_{k_v/k}(\mathscr S_{k_v}))_{k_v} \cong \mathscr S_{k_v} \times_{k_v} \mathscr S_{k_v}$ given by $x \mapsto (x, \sigma(x))$, with $\sigma$ on the left amounting to the switch of the two factors on the right, so $\overline {\mathscr O}_{k_v} \simeq \mathscr O_{k_v}$ where the $\mathrm{Gal}(k_v/k)$-action on the left corresponds to the interchange of the two $\mathscr S_{k_v}$ factors on the right. That is, $\overline{\mathscr O}$ is the $k_v/k$-twist of $\mathscr O$ over $\mathcal O$, although, of course, when $k_v \simeq k \times k$ we also get $\overline{\mathscr O} \simeq \mathscr O$ by restricting their isomorphism after base change to $k_v$ to one of the factors $k$. By checking after base change to $k_v$ and using \Cref{thm:dt}, we find that the $L$-translation action based at $s \in \overline {\mathscr O}(k)$ gives a $k$-isomorphism $L/\mathrm{Stab}_L(v)^0 \cong \overline{\mathscr O}$. For \ref{m:LM-b} and \ref{m:LM-c}, it remains to note that, by \Cref{lem:S-ratl}, both $\mathrm{Res}_{k_v/k}(\mathscr S_{k_v})$ and $\overline {\mathscr O}$ are rational, geometrically integral varieties over $k$ of dimension $32$.
\epf

\brem
We will not need this, but one may also describe the $H$-stabilizer of $v$:
\be\label{eqn:d}
\mathrm{Stab}_H(v) \cong \left\{
\begin{array}{ll}
 \mathrm{Stab}_L(v)^0 & \text{if }\,  \mathrm{char}(k)\neq 2, \\[1ex]
  \mathrm{Stab}_L(v) & \text{if }\, \mathrm{char}(k)= 2,
\end{array}
\right.
\ee
see \cite{GG17}*{Table 1}. Seeing the dichotomy \eqref{eqn:d} in \emph{loc.~cit.}~was very helpful for the development of the ideas of this article because it suggested that $\mathrm{Stab}_L(v)^0$ rather than $\mathrm{Stab}_H(v)$ might be the right object to focus on, and so that it might be fruitful to investigate the rationality of $L/\mathrm{Stab}_L(v)^0$ (as opposed to the much more evident rationality of $L/\mathrm{Stab}_L(v) \cong \mathcal O$).
\erem

Combined with the method for lifting $R$-triviality developed in \S\ref{sec:R-homog} and with our work on the kernel of the Rost invariant in \S\ref{sec:Rost}, \Cref{cor:LM-ratl} gives the following consequence for the $R$-equivalence on $L $.

\bprop\label{prop:LkRt}
For $k$, $H$, and $V$ as in \uS\uref{pp:setup}, if $\mathrm{dim}_2^{\mathrm{sep}}(k) \le 2$, then 
\[
L(k)/R \cong 1.
\]
\eprop

\bpf
If $k$ is finite, then, since $H$ is inner (see \S\ref{pp:setup}),  Lang's theorem \cite{Lan56}*{Theorem~2} implies that $H$ is split. Then $L$ is also split, and \Cref{thm:RR} gives the claim. In the remaining case of an infinite $k$, there is a $v \in \mathcal O(k)$ in the open $L$-orbit $\mathcal O \subset V$ from \Cref{thm:dt}, and we set $M \colonequals \mathrm{Stab}_L(v)^0$. \Cref{cor:LM-ratl} ensures that $M$ is semisimple, simply connected, absolutely indecomposable of type $A_5$ and that $L/M$ is rational over $k$. Thus, by \Cref{thm:RR} and \Cref{prop:retract}~\ref{m:R-viii}, we have 
\[
(L/M)(k)/R = \{*\} \qxq{and} M(k)/R \cong 1.
\]
Moreover, by \Cref{rem:dimp-finite-extn} and \Cref{thm:Rost-A2}, for every finite extension $k'/k$, every $M$-torsor over $k'(\!(t)\!)$ with a trivial Rost invariant is trivial. Thus, once we verify that for every field extension $K/k$, every $M$-torsor over $K$ that induces a trivial $L$-torsor has a trivial Rost invariant, the conditions of \Cref{prop:nonabelian-R} (notably, its condition \ref{iiipr}) will hold and we will get the desired $L(k)/R \cong 1$.

Our argument for the remaining claim about the vanishing of the Rost invariant will not use the separable dimension assumption, so we may assume that $K = k$. The Rost multiplier of the inclusion $M \hookrightarrow H$ of semisimple, simply connected, absolutely indecomposable $k$-groups is $1$: since all the coroots of $H$ have the same length (the Dynkin diagram of $H$ is simply laced, so use \cite{BouLGLA}*{chapitre VI, section 1, numéro 5, proposition 15}) and the coroots of $M$ are also coroots of $H$ (see \Cref{cor:LM-ratl}~\ref{m:LM-a}), this follows from \cite{GMS03}*{page 122, Proposition 7.9 (2) (with page 121, paragraph after Remark 7.5)}. Therefore, it remains to show that every $H$-torsor that induces a trivial $L$-torsor has a trivial Rost invariant. By the construction of $L$ (see \eqref{eqn:L}), such $H$-torsors are precisely the ones that reduce to $\mu_2$-torsors for the central $\mu_2 \subset H$ whose base change to $\overline k$ is the kernel of the ``other'' (different from $V_{\overline k}$) half-spin representation of $H_{\overline{k}}$.

By \S\ref{pp:high-wt}, \Cref{eg:half-spin}, and \cite{BouLGLA}*{chapitre VI, planche IV, (VIII)}, the inclusion $i\colon \mu_2 \subset Z_H$ is the unique one that on character groups induces the map $X^*(Z_H) \rightarrow X^*(\mu_2) \cong \mathbb Z/2\mathbb Z$ that does not kill the restriction $\chi|_{Z_H}$ of the fundamental weight of $H_{\overline{k}}$ corresponding to $V$ and kills the restriction $\chi^\prime|_{Z_H}$ of the other half-spin fundamental weight of $H_{\overline{k}}$. By \S\ref{pp:R-center}~\ref{RCiii}, on duals the isomorphism $Z_H \xrightarrow{\widehat \rho} Z_H^\circ$ sends the classes of the fundamental coweights to those of the corresponding fundamental weights. Thus, on the character groups the composition $Z_H \xrightarrow{\widehat \rho} Z_H^\circ \xrightarrow{i^\circ} \mu_2^\circ \cong \mu_2$ is an inclusion of the span of $\chi|_{Z_H}$, that is, this composition is $\chi|_{Z_H}$. Therefore, by \S\ref{pp:R-center}, the Rost invariant of $H$ on the image of our $H^1(k, \mu_2)$ is given by cupping with $H^2_{\mathrm{fppf}}(k, \chi)(-t_H)$, where $t_H$ is the Tits class of $H$. It then remains to note that $H^2_{\mathrm{fppf}}(k, \chi)(-t_H)$ vanishes because $V$ is defined over $k$, see \S\ref{pp:TA}.
\epf

\brem
Under the assumption that the data $k$, $H$, and $V$ come from a semisimple, simply connected, absolutely indecomposable $k$-group $G$ of Tits index $E_{7, 1}^{66}$ as in \Cref{prop:E7166-Levi}, the second half of the proof of \Cref{prop:LkRt} becomes simpler. Indeed, for the same reason as for $M \hookrightarrow H$, the Rost multiplier of $M \hookrightarrow G$ is also $1$ and this inclusion factors through $L$.
\erem

Our study of $E_{7,1}^{66}$ in \Cref{prop:E7166-Levi} now leads us to the following crucial consequence.

\bcor\label{cor:GRtrivial}
For a field $k$ with $\mathrm{dim}_2^{\mathrm{sep}}(k) \le 2$ and a semisimple, simply connected, absolutely indecomposable $k$-group $G$ of Tits index $E_{7,1}^{66}$, we have $G(k)/R \cong 1$.
\ecor

\bpf
Let $L$ be a Levi of a $k$-parabolic of $G$ considered in \Cref{prop:E7166-Levi}. The latter supplies the $H$ and $V$ as in \S\ref{pp:setup} and ensures that $L$ is as in \eqref{eqn:L}. \Cref{prop:LkRt} then shows that $L(k)/R \cong 1$. As is well known (see \cite{GS25}*{Corollary 3.10}), this immediately implies that $G(k)/R \cong 1$.
\epf

\section{The triviality of $K_1^G(k)$} \label{sec:final}
 
 We are ready to settle the isotropic case of \cite{Gil19}*{conjecture 9.5.2}. The key case in its proof—that of $E_{7,1}^{66}$—is the fruit of the previous sections, namely, of \Cref{cor:GRtrivial}. Our treatment of other groups, on the other hand, is similar to the proof of \cite{Gil09}*{théorème 8.6}, although the sharpness of our hypotheses \ref{m:m-i}–\ref{m:m-iii} forces us to find various slightly more general alternative routes in certain cases.

\bthm \label{thm:M}
Let $G$ be a semisimple, simply connected group over a field $k$ all of whose factors in the canonical decomposition $G \cong \prod \mathrm{Res}_{k'/k}(G')$ of \eqref{eqn:can-dec} are isotropic. Assume that if some $G^\prime$ is
\benumr
\m
\label{m:m-i} an $A_n$ and $D$ is its division algebra of \uS\uref{pp:conv} \uref{m:conv-i} \uref{m:conv-a}–\uref{m:conv-b}, then $\dim^{\mathrm{sep}}_p(k) \le 2$ for every prime $p$ whose square divides the index of $D$  \up{for every such $p$, we have $p^2 \mid  n + 1$}\uscolon

\m
\label{m:m-ii} nontrialitarian $D_n$ and its central simple algebra of \uS\uref{pp:conv} \uref{m:conv-ii} has index $\ge 8$, then~$\mathrm{dim}^{\mathrm{sep}}_2(k) \le 2$\uscolon

\m
\label{m:m-iii} of Tits index either ${}^2E_{6,1}^{35}$, or $E_{7, 1}^{48}$, or $E_{7, 1}^{66}$, or $E_{8, 1}^{91}$, or $E_{8, 1}^{133}$, then $\mathrm{dim}^{\mathrm{sep}}_2(k) \le 2$.

\eenum 
Then $G(k) = G(k)^+$\uscolon equivalently, then $K_1^G(k) \cong 1$\uscolon equivalently, then $G(k)/R \cong 1$.
\ethm
 
 The assumptions \ref{m:m-i}–\ref{m:m-iii} hold if $\mathrm{dim}^{\mathrm{sep}}_p(k) \le 2$ for every torsion prime $p$ of $G$,  see \cite{Ser94}*{section~2.2.3}. 
 
\bpf
The different formulations are equivalent by \Cref{prop:W-R}. By \S\ref{pp:K1Gk}, the formation of $K_1^G(k)$ commutes with products and restrictions of scalars, so  \Cref{rem:dimp-finite-extn} lets us assume that our isotropic $G$ is absolutely indecomposable. \Cref{prop:retract} and \cite{Tit66}*{pages 58–60} then reduce us to when either
\bitem
\m
$G$ is an $A_n$; or

\m
$G$ is a nontrialitarian $D_n$ and its central simple algebra of \S\ref{pp:conv} \ref{m:conv-ii} has index $> 4$; or

\m
the Tits index of $G$ is either ${}^2E_{6,1}^{35}$, or $E_{7, 1}^{48}$, or $E_{7, 1}^{66}$, or $E_{8, 1}^{91}$, or $E_{8, 1}^{133}$;
\eitem
and, in addition, $G$ is not quasi-split, so also $k$ is infinite.

When $G$ is of inner type $A_n$, we have $G \simeq \mathrm{SL}_r(D)$ as in \S\ref{pp:conv} \ref{m:conv-i} \ref{m:conv-a}. By \cite{Die43}*{théorème 1} and \cite{BT73}*{Proposition 6.2 (ii), (v)}, then $K_1^G(k) \cong (D^\times)^{\mathrm{Nrd}=1}/[D^\times, D^\times]$.\footnote{Division algebras permit a uniform statement, but \cite{Die43}*{section 5} also gives $K_1^{\mathrm{SL}_1(A)} \cong  (A^\times)^{\mathrm{Nrd} = 1}/[A^\times,A^\times]  $ for every nondivision central simple algebra $A$ over a field except precisely for the case when $A \simeq \mathrm{Mat}_{2 \times 2,\,\mathbb F_2}$.} Decomposing the latter as in \cite{Dra83}*{Section 9, Corollary 11 and Section 23, Lemma 6} leaves us with $D$ of prime power index. Then $K_1^G(k) \cong 1$ by applying either \cite{GB07}*{Theorem 1.1 and Corollary 5.5 (1)} or \Cref{prop:retract}~\ref{m:R-viii}.

When $G$ is of outer type $A_n$, we have $G \simeq \mathrm{SU}_r(D, h)$ as in \S\ref{pp:conv} \ref{m:conv-i} \ref{m:conv-b}. Moreover, Kummer theory supplies an $a \in (k'^\times)^{\sigma = -\mathrm{id}}$ (when $\mathrm{char}(k) = 2$, we may take $a = 1$), and $\mathrm{SU}_r(D, h) \cong \mathrm{SU}_r(D, ah)$, so we may freely switch between $\sigma$-Hermitian forms (such as $h$) and $\sigma$-anti-Hermitian forms (such as $ah$). Since $G$ is isotropic, so is $h$: in fact, $h(v, v) = 0$ for every $v \in D^r$ in a weight $\neq 0$ summand of the $k$-linear weight decomposition of $D^r$ with respect to any $\mathbb G_{m,\, k} \subset G$. For any $v \in D^r$ with $h(v, v) = 0$, the transvection $T_v \colon x \mapsto x + vah(v, x)$ preserves $h$ and pointwise fixes both the $k$-span of $v$ and its orthogonal complement, so the nondegeneracy of $h$ ensures that $T_v \in G(k)^+$. The nondegeneracy also implies that $T_v$ is not central in $G(k)$ (equivalently, in $G(\overline{k})$) when $v \neq 0$. Thus, since $gT_v g^{-1} = T_{gv}$ for $g \in G(k)$, the Tits simplicity theorem \cite{Tit64}*{Main theorem} shows that the $T_v$ generate $G(k)^+$ (our $k$ is already infinite). Consequently, \cite{Wal59}*{Theorem 1}\footnote{Since our $\sigma$ is nontrivial on $k^\prime$ and since being $\sigma$-Hermitian amounts to the same as being $\sigma$-anti-Hermitian in characteristic $2$, by \cite{Die63}*{chapitre I, section 10}, the condition \cite{Wal59}*{(5)} about $h$ being trace valued is automatic also in characteristic $2$, and the same applies to \cite{Yan74}*{top of p. 373}.} and \cite{Yan74}*{proof of Lemma 2} give the identification $K_1^G(k) \cong (D^\times)^{\mathrm{Nrd} \in k^\times}/\langle(D^\times)^{\sigma = 1}\rangle$. Moreover, by decomposing the latter using \cite{Yan78}*{Sections 2.2 and 2.7}, we may assume that $D$ is of prime power index. Then $K_1^G(k) \cong 1$ by applying either \cite{GB07}*{Theorem 1.1 and Corollary 5.5 (1)} or \Cref{prop:retract}~\ref{m:R-viii}.

When $G$ is of nontrialitarian type $D_n$ and its central simple algebra $A$ of \S\ref{pp:conv} \ref{m:conv-ii} has index $> 4$, we have $G \simeq \mathrm{Spin}(A, \sigma, f)$ as there. Since $G$ is isotropic, the same argument with a nontrivial $\mathbb G_{m,\, k}$-weight decomposition as in the case of the outer type $A_n$ combines with \cite{KMRT98}*{Corollary~1.2, Theorem 3.1, Theorem 4.2 (1), and after the second display on page~73} to show that the involution $\sigma$ of $A$ is isotropic in the sense of \cite{KMRT98}*{Definition 6.3}. Thus, by applying the Merkurjev comparison \cite{Mer98}*{Lemma 8.2 and Proposition 8.3} and, when $\mathrm{char}(k) = 2$, by first using the same lifting to characteristic $0$ argument based on \Cref{RS} and \cite{CF15}*{corollaire~2.7.0.32, propositions~4.4.0.34 and 4.5.2.7} as in the proof of \Cref{prop:retract}, we get
\[
G(k)/R \cong H(k)/R, \quad{\text{where}} \quad H \colonequals \mathrm{GL}_1(A) \times_{\mathrm{Nrd},\, \mathbb G_m,\, 2} \mathbb G_m
\]
and the pullback maps are the reduced norm of $\mathrm{GL}_1(A)$ and the squaring of $\mathbb G_m$. Since the index of $A$ is a power of $2$  (see \S\ref{pp:conv}) and $\mathrm{dim}_2^{\mathrm{sep}}(k) \le 2$, the reduced norm $\mathrm{Nrd}$ is surjective on $k$-points (see \S\ref{pp:sep-dim}). Thus, every $k$-point of $H$ is a product of a $k$-point of $\mathrm{SL}_1(A)$ and an image of a $k$-point under the $k$-map $\mathrm{GL}_1(A) \rightarrow H$ given by $g \mapsto (g^2, \mathrm{Nrd}(g))$ (which need not be a group homomorphism). Then $H(k)/R \cong 1$ because both 
\bitem
\m
$\mathrm{SL}_1(A)(k)/R \cong 1$ by the inner type $A_n$ case treated above and \Cref{prop:W-R} (our $A$ is not a division algebra because $\sigma$ is isotropic); and

\m
$\mathrm{GL}_1(A)(k)/R \cong 1$ due to the rationality of $\mathrm{GL}_1(A)$ (see \Cref{thm:RR}).
\eitem
When the Tits index of $G$ is $E_{7, 1}^{66}$, we apply \Cref{cor:GRtrivial}. The other remaining exceptional Tits indices ${}^2E_{6,1}^{35}$, $E_{7, 1}^{48}$, $E_{8, 1}^{91}$, and $E_{8, 1}^{133}$ are all excluded by \Cref{prop:no}.
\epf

\begin{appendix}

\section{Separable $p$-dimension via syntomic cohomology}\label{sec:appendix}
\label{sec:appendix}

We wish to reinterpret the (resp., separable) $p$-dimension of a field $k$ reviewed in \S\ref{pp:sep-dim} in terms of the vanishing of syntomic cohomology, taking advantage of the syntomic theory from \cite{BM23}. This extends \cite{brauer-purity}*{Theorem A.3 and Corollary A.4}, where we have reinterpreted the notion of a field of (separable) dimension $\le 1$ in terms of the vanishing of its flat cohomology in degrees $\ge 2$.

\bthm
\label{thm:dimp} 
Let $k$ be a field, $p$ a prime, and $d \in \mathbb Z_{\ge 0}$, and recall the notation $\mathscr F(j)$ from \uS\uref{pp:conv}. 
\benum
\m \label{m:dimp-a}
With the separable $p$-dimension $\mathrm{dim}^{\mathrm{sep}}_p(k)$ defined as in \uS\uref{pp:sep-dim}, the following are equivalent\ucolon
\benumr
\m\label{m:dimp-a-i}
$\mathrm{dim}^{\mathrm{sep}}_p(k) \le d$\uscolon

\m\label{m:dimp-a-ii}
$H^{d +1}_{\mathrm{\text{\upshape{ét}}}}(k', (\mathbb Z/p\mathbb Z)(d))\! \cong\! 0$ \up{equivalently, $H^{d +1}_{\mathrm{\text{\upshape{ét}}}}(k', (\mathbb Q_p/ \mathbb Z_p)(d))\! \cong\! 0$} for every finite separable extension $k^\prime/k$\uscolon

\m\label{m:dimp-a-iii}
$H^i_{\mathrm{\text{\upshape{ét}}}}(k, \mathscr F(j))\! \cong\! 0$ for every $p^\infty$-torsion sheaf $\mathscr F$ on $k_{\mathrm{\text{\upshape{ét}}}}$, every $i > d$, and every $i > j \ge 0$.
\eenum
These conditions propagate to every algebraic extension $\ell/k$, and, when $p = \mathrm{char}(k)$, the vanishing in \ref{m:dimp-a-ii} for a fixed $k^\prime$ propagates to every $d^\prime\ge d$ \up{and the same $k^\prime$}.

\m\label{m:dimp-b}
With the $p$-dimension $\mathrm{dim}_p(k)$ defined as in \uS\uref{pp:sep-dim}, the following are equivalent\ucolon
\benumr \addtocounter{enumi}{-1}
\primeitem\label{m:dimp-b-i}
$\mathrm{dim}_p(k) \le d$\uscolon

\addtocounter{enumi}{1}
\primeitem\label{m:dimp-b-ii}
$H^{d +1}_{\mathrm{\text{\upshape{ét}}}}(k, \mathbb (\mathbb Z/p\mathbb Z)(d + 1)) \cong 0$ and $H^{d +1}_{\mathrm{\text{\upshape{ét}}}}(k', (\mathbb Z/p\mathbb Z)(d)) \cong 0$ for every finite extension $k^\prime/k$\uscolon

\addtocounter{enumi}{1}
\primeitem\label{m:dimp-b-iii}
$H^i_{\mathrm{\text{\upshape{ét}}}}(k, \mathscr F(j)) \cong 0$ for every $p^\infty$-torsion sheaf $\mathscr F$ on $k_{\mathrm{\text{\upshape{ét}}}}$, every $i > d$, and every $j \ge 0$.

\eenum

\eenum
\ethm

\bpf
By \cite{BM23}*{Example 1.3}, in the case when $p = \mathrm{char}(k)$, we have the interpretation $(\mathbb Z/p^n \mathbb Z)(d) \cong W_n \Omega^d_{k,\, \mathrm{log}}[-d]$ in terms of the étale sheaf of logarithmic de Rham–Witt forms. Thus, in all cases, the norm-residue isomorphisms \cite{Voe11}*{Theorem 6.16} and \cite{BK86}*{Corollary 2.8} give the following identification with Milnor $K$-theory:
\[
H^{d}_{\mathrm{\text{\upshape{ét}}}}(k, (\mathbb Q_p/ \mathbb Z_p)(d)) \cong K_d^M(k) \otimes_{\mathbb Z} \mathbb Q_p/\mathbb Z_p.
\]
Thus, $H^{d}_{\mathrm{\text{ét}}}(k, (\mathbb Q_p/ \mathbb Z_p)(d))$ is $p$-divisible, so $H^{d + 1}_{\mathrm{\text{ét}}}(k, \mathbb (\mathbb Z/p\mathbb Z)(d)) \cong H^{d + 1}_{\mathrm{\text{ét}}}(k, (\mathbb Q_p/ \mathbb Z_p)(d))[p]$. When applied over a variable finite separable $k^\prime/k$ and combined with a limit and dévissage argument for $H^{d + 1}_{\mathrm{\text{ét}}}(k, (\mathbb Q_p/ \mathbb Z_p)(d)) \cong \varinjlim_{n > 0}H^{d + 1}_{\mathrm{\text{ét}}}(k, (\mathbb Z/ p^n\mathbb Z)(d))$, this gives the parenthetical equivalence~in~\ref{m:dimp-a-ii}. 

For the rest, we start with the case $p \neq \mathrm{char}(k)$. We let $\mathscr F$ absorb the twist and, by restriction-corestriction, may therefore ignore the twists in \ref{m:dimp-a-ii}–\ref{m:dimp-a-iii} and \ref{m:dimp-b-ii}–\ref{m:dimp-b-iii}. The definition of  cohomological dimension \cite{Ser02}*{Section I.3.1} (see also \S\ref{pp:sep-dim}) then shows that \ref{m:dimp-a-i}, \ref{m:dimp-a-iii}, \ref{m:dimp-b-i}, and \ref{m:dimp-b-iii} are all equivalent. By \cite{Ser02}*{Section I.4.1, Proposition 21${}^\prime$} and pushforwards from finite extensions, they are also equivalent to \ref{m:dimp-a-ii} and \ref{m:dimp-b-ii}. Finally, a limit argument and \ref{m:dimp-b-i}–\ref{m:dimp-b-iii} give the passage to $\ell/k$. %, whereas, for a fixed $k^\prime$, by \cite{Ser02}*{Section I.3.3, Corollary 1 to Proposition 14; Section I.4.1, Proposition 21}, the vanishing in \ref{m:dimp-a-ii} for a fixed $k^\prime$ propagates to $d^\prime\ge d$.

Suppose from now on that $p = \mathrm{char}(k)$ and let us start with the propagation of the vanishing in \ref{m:dimp-a-ii} for a fixed $k^\prime$ to $d' \ge d$. By \cite{GO08}*{Proposition~2.2.1}, we have the concrete Cartier-style description 
\be\label{eqn:Cartier}
H^{d +1}_{\mathrm{\text{\upshape{ét}}}}(k', \mathbb (\mathbb Z/p\mathbb Z)(d))  \cong \mathrm{Coker}\left(\Omega^d_{k'} \xrightarrow{a \frac{  \mathrm{d}x_1}{x_1}\wedge \dotsc\wedge\frac{\mathrm{d}x_d}{x_d} \mapsto (a - a^p) \frac{\mathrm{d}x_1}{x_1}\wedge \dotsc\wedge\frac{\mathrm{d}x_d}{x_d} }\Omega^d_{k'}/\mathrm{d}\Omega^{d - 1}_{k^\prime}\right).
\ee
Thus, for the claimed propagation, it suffices to observe that $\mathrm{d}\eta \wedge \frac{\mathrm{d}x_{d + 1}}{x_{d + 1}} = \mathrm{d}(\eta \wedge \frac{\mathrm{d}x_{d + 1}}{x_{d + 1}})$ for $\eta \in \Omega_{k^\prime}^{d - 1}$.

Since $\Omega^{d + 1}_{k,\, \mathrm{log}}$ vanishes precisely when the $p$-rank of $k$ is $\le d$ (see \cite{SP}*{\href{https://stacks.math.columbia.edu/tag/07P2}{Lemma 07P2}}), both \ref{m:dimp-a-i}$\Leftrightarrow$\ref{m:dimp-a-ii} and \ref{m:dimp-b-i}$\Leftrightarrow$\ref{m:dimp-b-ii} amount to the definitions from \S\ref{pp:sep-dim}. Moreover, pushforwards from finite separable extensions and the projection formula give \ref{m:dimp-a-iii}$\Rightarrow$\ref{m:dimp-a-ii}. For the converse \ref{m:dimp-a-ii}$\Rightarrow$\ref{m:dimp-a-iii}, by a limit argument, we may first assume that $\mathscr F$ is a finite commutative étale $k$-group scheme killed by some $p^n$ and then, by dévissage, that $\mathscr F$ is killed by $p$. At this point, since $\mathscr F(j) \cong \mathscr F \otimes_{\mathbb F_p}\underline{\Omega}^j_{k,\, \mathrm{log}}[-j]$ lives in at most one degree (namely, $j$) and $\mathrm{cd}_p(k) \le 1$ (see \cite{Ser02}*{Section~II.2.2, Proposition 3}), in \ref{m:dimp-a-iii} we may restrict to $i = j + 1$ and get that  $H^{j + 1}_{\mathrm{\text{ét}}}(k, \mathscr F(j))$ is a quotient of $H^{j + 1}_{\mathrm{\text{ét}}}(k, \widetilde {\mathscr F}(j))$ whenever $\mathscr F$ is a quotient of an $\mathbb F_p$-local system $\widetilde {\mathscr F}$ on $k_{\mathrm{\text{ét}}}$. By choosing $\widetilde {\mathscr F}$ to correspond to a direct sum of permutation modules given by the Galois orbits of the elements of a finite generating set for $\mathscr F(k^s)$ and combining the Shapiro lemma with the settled upward propagation in \ref{m:dimp-a-ii}, we get that \ref{m:dimp-a-ii}$\Rightarrow$\ref{m:dimp-a-iii}.

For the claim about $\ell/k$, propagation to separable algebraic extensions follows from \ref{m:dimp-a-ii} and a limit argument, so we may assume that $\ell/k$ is purely inseparable. By first replacing our $k$ by a $k^\prime$ from \ref{m:dimp-a-ii} if necessary (and $\ell$ by $\ell k^\prime$), we need to argue that $H^{d +1}_{\mathrm{\text{ét}}}(\ell, (\mathbb Q_p/ \mathbb Z_p)(d)) \cong 0$ (equivalently, that $H^{d +1}_{\mathrm{\text{ét}}}(\ell, \mathbb F_p(d)) \cong 0$) granted the analogous vanishing over $k$. For this, by a limit argument based on \cite{BM23}*{footnote 1}, we may assume that $\ell/k$ is finite and then, by dévissage, even that $\ell^p \subset k$. By the aforementioned identification with the logarithmic de Rham--Witt forms and by \cite{Ill79}*{after I.5.7.1} (alternatively, by \cite{Ill79}*{I.2.18.7 and I.5.7.2.1}), the absolute Frobenius of $\ell$ induces multiplication by $p^d$ on each $(\mathbb Q_p/ \mathbb Z_p)(d)$. Thus, by the $\mathrm{cd}_p(\ell)\le 1$ argument as above, multiplication by $p^d$ is surjective on $H^{d +1}_{\mathrm{\text{ét}}}(\ell, (\mathbb Q_p/ \mathbb Z_p)(d))$, while, by the factorization of the $p$-power map of $\ell$ through the inclusion $k \subset \ell$, this multiplication by $p^d$ factors through $H^{d +1}_{\mathrm{\text{ét}}}(k, (\mathbb Q_p/ \mathbb Z_p)(d))$. By \ref{m:dimp-a-ii}, the latter vanishes, so it follows that $H^{d +1}_{\mathrm{\text{ét}}}(\ell, (\mathbb Q_p/ \mathbb Z_p)(d))$ vanishes as well, as desired.

The implication \ref{m:dimp-b-iii}$\Rightarrow$\ref{m:dimp-b-ii} now follows from \ref{m:dimp-a}, especially, from its aspect about passing to $\ell/k$. For the remaining \ref{m:dimp-b-ii}$\Rightarrow$\ref{m:dimp-b-iii}, it suffices to note that \ref{m:dimp-b-ii}$\Rightarrow$\ref{m:dimp-a-ii}$\Rightarrow$\ref{m:dimp-a-iii} and that, as observed above, \ref{m:dimp-b-ii} implies that $\underline{\Omega}^{j}_{k,\, \mathrm{log}} \cong 0$ for $j  > d$, so that also  $\mathscr F(j) \cong 0$ for such $j$. 
\epf

\brem
It seems to be an open problem whether the vanishing in \ref{m:dimp-a-ii} for a fixed $k^\prime$ also propagates to $d' \ge d$ when $p \neq \mathrm{char}(k)$, see \cite{MW07}*{before Section 1} already for the case $d = 0$. 
\erem

\section{On the Kneser--Tits problem for $E_{7,1}^{66}$} \label{app:B}

\renewcommand{\thefootnote}{\fnsymbol{footnote}}

\center{Alexis Bouthier\footnote[1]{Sorbonne Universit\'e, Universit\'e Paris Cit\'e, CNRS, IMJ-PRG, F-75005 Paris, France. E-mail address: \url{abouthier@imj-prg.fr}} and Federico Scavia\footnote[2]{CNRS, Institut Galil\'ee, Universit\'e Sorbonne Paris Nord, 99 avenue Jean-Baptiste Cl\'ement, 93430 Villetaneuse, France. E-mail address: \url{scavia@math.univ-paris13.fr}}}

\justify
\renewcommand{\addtocontents}[2]{}

\renewcommand{\thenumberingbase}{\thesection.0.\arabic{numberingbase}}

The purpose of this appendix is to give an independent alternative proof of the last remaining case of the Kneser--Tits problem for fields of cohomological dimension $2$ and to settle the isotropic part of a question \cite{Gil09}*{Question 8.5} asked by Gille.

	\begin{theorem}\label{thm:main}
		Let $k$ be a field such that $\mathrm{cd}_2(k)\leq 2$, and let $G$ be an absolutely almost simple simply connected $k$-group of Tits index $E_{7,1}^{66}$. If $\mathrm{char}(k)=2$, assume that $k$ is perfect. Then~$W(k,G)=1$.
	\end{theorem}
	
	The assumption that $k$ is perfect in characteristic $2$ is made only for simplicity, and it can be dropped, provided that the conclusion of the next lemma holds for $k$, with only minor modifications of the arguments.
	
	\begin{lemma}\label{lem:norm}
		Let $k$ be a field with $\mathrm{char}(k) \neq 2$ such that $\mathrm{cd}_2(k)\leq 2$. For every finite \'etale $k$-algebra $F$ and every quaternion $F$-algebra $D$, the reduced norm $\mathrm{Nrd}_{D/F}\colon D^\times\to F^\times$ is surjective, that is, the pointed set $H^1(F,\mathrm{SL}_1(D))$ is trivial.
	\end{lemma}
	
	\begin{proof}
		See \S\ref{pp:sep-dim}.
		%It suffices to consider the case of a finite field extension $F/k$. Write $D=(a,b)_F$, and let $c\in F^\times$. Recall that the reduced norm form of the quaternion algebra $(a,b)_F$ is the two-fold Pfister form $\langle\!\langle a,b\rangle\!\rangle$. Therefore, $c$ is a reduced norm from $D$ if and only if the $3$-fold Pfister form $\langle\!\langle a,b,c\rangle\!\rangle$ is hyperbolic. By the norm-residue isomorphism, this is equivalent to the vanishing of $(a,b,c)\in H^3(F,\mu_2^{\otimes3})$. Since $\mathrm{cd}_2(F)\leq 2$, the latter group is zero, as desired.
	\end{proof}
	
\addtocounter{numberingbase}{-2}

\csub[Reduction to \texorpdfstring{$H$}{H}]

\addtocounter{numberingbase}{-1}

\renewcommand{\thenumberingbase}{\thesection.1.\arabic{numberingbase}}

	Let $k$ be a field with $\mathrm{cd}_2(k)\leq 2$, and let $G$ be an absolutely almost simple simply connected $k$-group of Tits index $E_{7,1}^{66}$. By \cite{BouLGLA}*{chapitre VI, planche~VI}, its Dynkin diagram is
	\[
	\begin{tikzpicture}[
		x=1cm, y=1cm,
		every node/.style={font=\small}
		]

		\draw (0,0) -- (5,0);
		\draw (2,0) -- (2,1);
		
		\foreach \x/\i in {0/1,1/3,2/4,3/5,4/6,5/7} {
			\fill (\x,0) circle (2pt);
			\node[below=5pt] at (\x,0) {$\alpha_{\i}$};
		}
		\fill (2,1) circle (2pt);
		\node[above=5pt] at (2,1) {$\alpha_2$};

		\draw (0,0) circle (5pt);
	\end{tikzpicture}
	\]
	Let $S\simeq\mathbb G_m$ be a maximal split torus of $G$,
	let $L\ce Z_G(S)$, and let $H\ce[L,L]$. Since $G$ is simply connected, the group $H$ is simply connected, anisotropic, and of inner type $D_6$; see \cite{PS11}*{Theorem~3, type $E_7$} and \cite{MT11}*{Proposition~12.14}.
	
	\begin{remark}\label{rem:char2}
		If $\mathrm{char}(k)=2$ and $k$ is perfect, no group $G$ as above exists. Indeed, the anisotropic kernel $H$ is isomorphic to $\mathrm{Spin}(A,\sigma,f)$, where $(A,\sigma,f)$ is a degree-$12$ algebra with quadratic pair, and since $\sigma$ is an involution of the first kind we have $2[A]=0$ in $\mathrm{Br}(k)$; see \cite{KMRT98}*{Theorem~3.1(1)}. Since $k$ is perfect of characteristic $2$, by \cite{KMRT98}*{Proposition~(5.11)} we have $\mathrm{Br}(k)[2]=0$, and hence $A$ is split, so that $H\simeq\mathrm{Spin}(q)$ for a nonsingular quadratic form $q$ of dimension $12$. Over a perfect field of characteristic $2$, every nonsingular binary quadratic form represents $1$, hence every nonsingular quadratic form of dimension at least $4$ is isotropic, and so $H$ is isotropic, a contradiction.
	\end{remark}
	
	In view of Remark~\ref{rem:char2}, we assume from now on that $\mathrm{char}(k)\ne2$. Since $G$ is not quasi-split (its anisotropic kernel has positive semisimple rank) and every connected reductive group over a finite field is quasi-split, we may also assume that $k$ is infinite.

	\begin{lemma}\label{lem:Levi-isogeny}
		The multiplication map $H\times S\to L$ induces a short exact sequence
		\[1\longrightarrow \mu_2\longrightarrow H\times S\longrightarrow L\longrightarrow 1.\]
	\end{lemma}
	
	\begin{proof}
		See \Cref{prop:E7166-Levi}~\ref{m:E7-b}.
	\end{proof}

	\begin{proposition}\label{prop:reduction}
		In order to prove that $W(k,G)=1$, it suffices to prove $H(k)/R=1$.
	\end{proposition}

	\begin{proof}
		By \cite{Gil09}*{Theorem~7.2}, we have $W(k,G)\simeq G(k)/R$. Using the big cell associated with $L$, we see that $G(k)/R\simeq L(k)/R$. Since the only
		torsion prime of type $D_6$ is $2$, Gille's theorem
		\cite{Gil01}*{Theorem~6(b)} applied to the sequence of Lemma~\ref{lem:Levi-isogeny} gives an
		exact sequence
		\[(H\times S)(k)/R\longrightarrow L(k)/R\longrightarrow H^1(k,\mu_2)/R\longrightarrow 1.\]
		By Kummer theory, $H^1(k,\mu_2)\simeq k^\times/k^{\times2}$, and every class $[a]\in k^\times/k^{\times2}$ is $R$-equivalent to $1$: indeed,
		the class of $1-t+ta\in k(t)^\times$
		specializes to $[1]$ at $t=0$ and to
		$[a]$ at $t=1$. Thus $H^1(k,\mu_2)/R=1$, and hence the induced map $(H\times S)(k)/R\to L(k)/R$
		is surjective. Finally, $S\simeq\mathbb G_m$ is $R$-trivial, so the map
		$H(k)/R\to L(k)/R$ is surjective. Combining this with $G(k)/R\simeq L(k)/R$ proves the assertion.
	\end{proof}

	\begin{lemma}[Garibaldi--Qu\'eguiner-Mathieu]\label{lem:quadratic}
		There exists a central simple $k$-algebra $A$ of degree $12$ with an orthogonal involution $\sigma$ of trivial discriminant and trivial Clifford invariant such that $H\simeq\mathrm{Spin}(A,\sigma)$. Moreover, there is a quadratic subfield $K=k(i)\subset A$ with $d\ce i^2\in k^\times$ and $\sigma(i)=-i$. Letting $B\ce Z_A(K)=Z_A(i)$ and $\tau\ce\sigma|_{B}$, the involution $\tau$ is unitary for $K/k$, and
		\[\deg_K B=6,\qquad \mathrm{ind}(B)\leq 2,\qquad (A,\sigma)_K\text{ is hyperbolic}.\]
	\end{lemma}
	
	\begin{proof}
		By \cite{GQM09}*{Remark~3.3}, the anisotropic kernel $H$ of $G$ is isogenous to $\mathrm{SO}(A,\sigma)$ for an algebra with involution satisfying the hypotheses of \cite{GQM09}*{Theorem~3.1}. Since $H$ is simply connected, it follows that $H\simeq\mathrm{Spin}(A,\sigma)$. Moreover, \cite{GQM09}*{Theorem~3.1} and its proof yield a quadratic \'etale $k$-subalgebra $K=k(i)\subset A$ satisfying $i^2\in k^\times$ and $\sigma(i)=-i$, such that $\tau\ce\sigma|_{B}$ is a unitary involution of $B\ce Z_A(K)$ for $K/k$, and $(A,\sigma)_K$ is hyperbolic.
		
		Suppose that $K\simeq k\times k$ is split, and let $e_1=(1,0)$ and $e_2=(0,1)$ be its primitive idempotents. Since $\sigma(i)=-i$, the restriction of $\sigma$ to $K$ is the nontrivial automorphism of $K/k$, and hence exchanges the two factors. Thus $\sigma(e_1)=e_2=1-e_1$ and hence by \cite{KMRT98}*{Proposition~6.7} the involution $\sigma$ is hyperbolic. Consequently, $H\simeq\mathrm{Spin}(A,\sigma)$ is isotropic, a contradiction. Hence $K$ is a field. By the double-centralizer theorem \cite{KMRT98}*{Theorem~1.5},
		\[\deg_K B=\frac{\deg_k A}{[K:k]}=6,
		\qquad
		[B]=[A\otimes_k K]\quad\text{in $\mathrm{Br}(K)$.}\]
		Since $\sigma$ is an involution of the first kind, we have $2[A]=0$, and hence $2[B]=0$. Thus $\mathrm{ind}(B)$ is a power of $2$. As it also divides $\deg_K B=6$, we obtain $\mathrm{ind}(B)\leq 2$.
	\end{proof}

\addtocounter{subsection}{-3}

\csub[The group \texorpdfstring{$M$}{M}]

\addtocounter{numberingbase}{-2}

\renewcommand{\thenumberingbase}{\thesection.2.\arabic{numberingbase}}

	Let $M_0\ce Z_{\mathrm{SO}(A,\sigma)}(i)$; by \cite[Remark~1.15]{GQM09}, we have $M_0=\mathrm{U}(B,\tau)$. Let $\pi\colon H\to\mathrm{SO}(A,\sigma)$ be the canonical central isogeny, let $M\ce\pi^{-1}(M_0)$, and let $T_0\ce R^1_{K/k}\mathbb G_m$.
	
	\begin{lemma}\label{lem:centralizer}
		We have an isomorphism
		\[M\simeq\{(x,t)\in M_0\times T_0:\mathrm{Nrd}_{B/K}(x)=t^2\}\]
		which identifies $\pi|_M$ with the first projection $M\to M_0$.
	\end{lemma}
	
	\begin{proof}
		Since the derived subgroup $(M_0)^{\mathrm{der}}$ of $M_0$ is the simply connected $k$-group $\mathrm{SU}(B,\tau)$, the map $\pi$ identifies the derived subgroup $M^{\mathrm{der}}$ of $M$ with $(M_0)^{\mathrm{der}}$. Let $T=M/M^{\mathrm{der}}$. The reduced norm
		gives an isomorphism $M_0/(M_0)^{\mathrm{der}}\simeq T_0$, so the composite
		$\mathrm{Nrd}_{B/K}\circ\pi$ induces a degree-two
		isogeny of rank-$1$ tori $\psi\colon T\to T_0$, and hence there is an isomorphism	$\varphi\colon T\to T_0$ such that $\psi$ factors as the composite $T\xrightarrow{\varphi}T_0\xrightarrow{\times 2}T_0$.
		Let $\chi\colon M\to T_0$ be the composite
		of the quotient map $M\to T$ with $\varphi$. 
		Then $\chi^2=\mathrm{Nrd}_{B/K}\circ\pi$, and hence we obtain a homomorphism
		\[M\longrightarrow\{(x,t)\in M_0\times T_0:\mathrm{Nrd}_{B/K}(x)=t^2\},\qquad m\longmapsto(\pi(m),\chi(m)),\]
		which induces isomorphisms on the derived subgroups and on the abelianization tori, and hence it is an isomorphism.
	\end{proof}
	
	\begin{lemma}\label{lem:M}
		The group $M$ is $R$-trivial.%\footnote{In fact, $M$ is rational, but the proof is more difficult (use \cite[Theorem 4.2]{KMRT} to construct a section of $M\to T_0$).}
	\end{lemma}
	
	\begin{proof}
		We use the description of $M$ given in Lemma~\ref{lem:centralizer}. By \cite[Proposition~2.4]{CM98},
		the group $M_0$ is rational. The morphism
		\[M_0\longrightarrow M,\qquad x\longmapsto(x^2,\mathrm{Nrd}_{B/K}(x))\]
		sends the identity to the identity, and its image contains every square in $M(k)$, proving that $M(k)/R$ is killed by $2$. On the other hand, multiplication induces a
		central isogeny
		\[(M_0)^{\mathrm{der}}\times T_0\longrightarrow M,\qquad (v,t)\longmapsto(tv,t^3).\]
		The torus $T_0$ is rational and, since $\mathrm{ind}(B)\leq 2$, the group
		$M_0^{\mathrm{der}}\simeq \mathrm{SU}(B,\tau)$ is rational by
		\cite[Corollary~6.5]{CM98}. Thus  $M(k)/R$ is killed by $3$, and hence it is trivial.
	\end{proof}

\csub[The \'etale algebras \texorpdfstring{$E_h$}{E_h} and the groups \texorpdfstring{$J_h$}{J_h}]

\addtocounter{numberingbase}{-2}

\renewcommand{\thenumberingbase}{\thesection.3.\arabic{numberingbase}}

\renewcommand{\theequation}{\thesection.3.\arabic{numberingbase}}

	Let $h\in H(k)$, set $a\ce \pi(h)$ and $j\ce aia^{-1}$, and define
	\[z\ce ij+ji,\qquad w\ce ij-ji.\]
	The identities $i^2=j^2=d$ and $\sigma(i)= -i$, $\sigma(j)=-j$ give
	\begin{equation}\label{eq:quaternion-relations}
		\begin{aligned}
			zi=iz\quad, zj=jz,\quad wi=-iw,\quad wj=-jw,\\
			w^2=z^2-4d^2,\quad \sigma(z)=z,\quad \sigma(w)=-w.
		\end{aligned}
	\end{equation}
	
	\addtocounter{numberingbase}{-1}

	\begin{lemma}\label{lem:generic-geometry}
	
		Let $(W,q)$ be a $12$-dimensional quadratic space over $k^s$ such that 
		\[
		(A,\sigma)_{k^s}\simeq (\mathrm{End}_{k^s}(W),\mathrm{ad}_q),
		\]
		and view $i,j,z$ as endomorphisms of $W$. There exists a dense open subscheme $\Omega\subset H$ such that, for every $h\in\Omega(k)$, the following assertions hold. 
		\benumr
			\item The endomorphism $z$ has three distinct eigenvalues
			$\lambda_1,\lambda_2,\lambda_3\in k^s$, each of
			multiplicity $4$, such that $\lambda_l\ne\pm2d$ for $l=1,2,3$.
			\item The corresponding eigenspace decomposition
			\[
			W=W_1\perp W_2\perp W_3
			\]
			is orthogonal, with each $W_s$ nondegenerate of
			dimension $4$.
			\item For $l\in\{1,2,3\}$, the eigenspace $W_l$ is stable under $i$ and $j$, and the two eigenspaces of $i$ on $W_l$ have dimension $2$.
		\end{enumerate}
	\end{lemma}

	\begin{proof}
		Over $k^s$, choose bases such that
		\[
		i=\sqrt d\begin{pmatrix}I_6&0\\0&-I_6\end{pmatrix},
		\qquad
		j=\sqrt d\begin{pmatrix}P&Y\\Z&-P^t\end{pmatrix},
		\qquad
		z=2d\begin{pmatrix}P&0\\0&P^t\end{pmatrix},
		\]
		where $Y,Z$ are alternating, $PY=YP^t$, and
		$P^2+YZ=I_6$.
		
		Let $\Omega'\subset H$ be the open subset where $z^2-4d^2\in A^\times$, or equivalently where both $Y$ and $Z$ are invertible. Over $\Omega'$, the matrix $P$ is self-adjoint with respect to the alternating form
		$Y^{-1}$, so that $\det(tI_6-P)=f(t)^2$, for some cubic polynomial $f$, and the characteristic polynomial $\chi_z(t)$ of $z$ satisfies $\chi_z(t)=p_z(t)^4$, where $p_z(t)=(2d)^3f(t/(2d))$. Since $\mathrm{char}(k)\neq 2$, comparison of coefficients shows that $p_z$ has coefficients in $k[\Omega']$. 
		
		Let $\Omega\subset \Omega'$ be the open subset where $p_z$ is separable. In order to show that $\Omega$ is non-empty, we choose distinct
		$p_1,p_2,p_3\in k^s\setminus\{\pm1\}$ and set
		\[
		P=\mathrm{diag}(p_1I_2,p_2I_2,p_3I_2),\quad
		Y=\mathrm{diag}(J,J,J),\quad
		Z=(I_6-P^2)Y^{-1},
		\]
		where
		\[
		J=\begin{pmatrix}0&1\\-1&0\end{pmatrix}.
		\]
		Then $j^2=d$ and $\sigma(j)=-j$. For each $l=1,2,3$, letting $x_l$ be the restriction of $i+j$ to $W_l$, we have
		\[
		x_li=jx_l,\quad
		x_l^2=2d(1+p_l),\quad
		\sigma(x_l)x_l=-2d(1+p_l),\quad \det(x_l)=4d^2(1+p_l)^2,
		\]
		so that the element
		\[
		g_l\ce \frac{x_l}{\sqrt{-2d(1+p_l)}}
		\in\mathrm{SO}(W_l)
		\]
		satisfies $g_l(i|_{W_l})g_l^{-1}=j|_{W_l}$. Hence $g\ce \bigoplus_l g_l\in\mathrm{SO}(W)$ satisfies \[
		gig^{-1}=j,\qquad p_z(t)=\prod_{l=1}^3(t-2dp_l),
		\qquad z^2-4d^2\in A^\times,
		\]
		and picking a lift of $g$ to $\mathrm{Spin}(W)(k^s)$ we deduce that $\Omega$ is non-empty.
		\benumr
		\m For every $h\in\Omega(k)$, the generalized eigenspaces of $P$ are $2$-dimensional nondegenerate subspaces for the alternating form $Y^{-1}$. Since every self-adjoint endomorphism of a $2$-dimensional space is scalar, we have $p_z(z)=0$. Thus $z$ is similar to $\mathrm{diag}
		(\lambda_1I_4,\lambda_2I_4,\lambda_3I_4)$ over $k^s$, where the $\lambda_l$ are pairwise distinct and
		$\lambda_l^2\neq 4d^2$ for $l=1,2,3$.
		
		\m This follows from the fact that, since $\sigma(z)=z$, eigenspaces corresponding
		to distinct eigenvalues are orthogonal.
		
		\m Let $l\in\{1,2,3\}$. Since $[z,i]=[z,j]=0$, the subspace $W_l$ is stable
		under $i$ and $j$. Writing $W_l^\pm=\ker(i|_{W_l}\mp\sqrt d\,\mathrm{id}_{W_l})$,
		the fact that $\sigma(i)=-i$ implies
		\[
		W_l=W_l^+\oplus W_l^-,
		\qquad
		W_l^\pm\subset(W_l^\pm)^\perp.
		\]
		As $W_l$ is nondegenerate of dimension $4$, we have $\dim W_l^+=\dim W_l^-=2$.\qedhere
	\eenum
	\end{proof}

	\begin{lemma}\label{lem:genericE}
		Let $\Omega$ be as in Lemma~\ref{lem:generic-geometry}. For every $h\in\Omega(k)$, the algebra $E_h\ce k[z]$ is \'etale of degree $3$ and $z^2-4d^2\in E_h^\times$.
	\end{lemma}
	
	\begin{proof}
		By Lemma~\ref{lem:generic-geometry}, the element $z$ has three distinct eigenvalues over $k^s$, each of multiplicity $4$, and so its minimal polynomial is separable of degree $3$, that is, $E_h=k[z]$ is a degree-$3$ \'etale $k$-algebra. Moreover, none of the eigenvalues of $z$ equals $\pm 2d$, so $z^2-4d^2$ is invertible.
	\end{proof}
	
	For the rest of this section, let $\Omega\subset H$ be as in Lemma~\ref{lem:generic-geometry}, let $h\in\Omega(k)$, and set $E\ce E_h=k[z]$.
	Since $\sigma(z)=z$, the degree-$3$ \'etale $k$-algebra $E$ is fixed by $\sigma$. Define 
	\[D\ce E[i,j]=E[i,w]\subset Z_A(E),\qquad C=Z_{Z_A(E)}(D).
	\]
	From \eqref{eq:quaternion-relations} we see that $D$ is the quaternion $E$-algebra
	\[	D\simeq(d,z^2-4d^2)_E	\]
	and that $\sigma|_D$ is the canonical quaternion involution.
	
	\begin{lemma}\label{lem:block}
		The centralizer $Z_A(E)$ is an Azumaya $E$-algebra of degree $4$,
		and multiplication identifies it with $D\otimes_E C$, where $C$ is
		a quaternion $E$-algebra. The restrictions of $\sigma$ to $D$ and $C$
		are their canonical quaternion involutions.
	\end{lemma}
	\begin{proof}
		Consider the orthogonal decomposition $W=W_1\perp W_2\perp W_3$ given by  Lemma~\ref{lem:generic-geometry}. It follows that
		\[Z_A(E)\otimes_k k^s\simeq \prod_{s=1}^3\mathrm{End}_{k^s}(W_s)\simeq M_4(k^s)^3,\]
		so $Z_A(E)$ is an Azumaya $E$-algebra of degree $4$. By the double-centralizer theorem, multiplication identifies $Z_A(E)$ with $D\otimes_E C$, where $C$ is a quaternion $E$-algebra. Finally, since the spaces $W_s$ are nondegenerate, the restriction of $\sigma$ to $Z_A(E)$ is orthogonal. Since $\sigma|_D$ is canonical, hence symplectic, the restriction $\sigma|_C$ is symplectic as well, and therefore canonical; see \cite[Propositions~2.21 and~2.23(1)]{KMRT98}.
	\end{proof}
	
	We define
	\[\widetilde J_h\ce R_{E/k}(\mathrm{SL}_1(D)\times\mathrm{SL}_1(C)).\]
	By Lemma~\ref{lem:block}, multiplication gives an isomorphism $D\otimes_E C\simeq Z_A(E)\subset A$, and hence induces a homomorphism $\widetilde J_h\to\mathrm{SO}(A,\sigma)$.
	Since $\widetilde J_h$ is semisimple and simply connected, this homomorphism lifts to a homomorphism 
	\[\widetilde J_h\longrightarrow H,\]
	and we let $J_h\subset H$ and $J_{h,0}\subset\mathrm{SO}(A,\sigma)$	be the corresponding images.
	
	\begin{lemma}\label{lem:Jstructure} \hfill
		\benumr
			\item There is an exact sequence
			\[1\longrightarrow R^1_{E/k}\mu_2\longrightarrow\widetilde J_h\longrightarrow J_h\longrightarrow1.\]
			\item We have $J_h=\pi^{-1}(J_{h,0})$.
			\item We have \[J_{h,0}(k)=\{xy:\ x\in D^\times,\ y\in C^\times,\ \mathrm{Nrd}_{D/E}(x)\mathrm{Nrd}_{C/E}(y)=1\}.\]
		\end{enumerate}
	\end{lemma}

\vspace{5pt}	
	
	\begin{proof} \hfill
		\benumr
		\m By Lemma~\ref{lem:generic-geometry}, we have an isomorphism $E\otimes_kk^s\xrightarrow{\sim} (k^s)^3$ given by evaluation at $\lambda_1$, $\lambda_2$ and $\lambda_3$, and there is an orthogonal decomposition
		\[W=W_1\perp W_2\perp W_3,\qquad \dim_{k^s}W_l=4,\]
		where $W_l$ is the eigenspace of $z$ corresponding to $\lambda_l$, and we define $D_l$ and $C_l$ as the factors of $D\otimes_k k^s$ and $C\otimes_k k^s$ corresponding to $\lambda_l$, for $l=1,2,3$. Since $E=k[z]$, every element of $Z_A(E)$ commutes
		with $z$ and hence preserves each eigenspace $W_l$. Thus, by Lemma~\ref{lem:block}, we obtain a commutative diagram
		\[
		\begin{tikzcd}
			(D\otimes_E C)\otimes_k k^s
			\arrow[r,"\sim"]
			\arrow[d,hook]
			&
			\prod_{l=1}^3(D_l\otimes_{k^s}C_l)
			\arrow[r,"\sim"]
			\arrow[d,hook]
			&
			\prod_{l=1}^3\mathrm{End}_{k^s}(W_l)
			\arrow[d,hook]
			\\
			A\otimes_k k^s
			\arrow[r,equal]
			&
			A\otimes_k k^s
			\arrow[r,"\sim"]
			&
			\mathrm{End}_{k^s}(W),
		\end{tikzcd}
		\]
		where the vertical map on the left is given by multiplication, and the one on the right is the block-diagonal embedding, and hence a commutative diagram of $k^s$-groups \addtocounter{numberingbase}{-2}
		\begin{equation}\label{eq:rectangle}
		\begin{tikzcd}
			(\widetilde J_h)_{k^s}
			\arrow[r,"\sim"]
			\arrow[d,"\varphi_{k^s}"]
			&
			\prod_{l=1}^3
			(\mathrm{SL}_1(D_l)\times\mathrm{SL}_1(C_l))
			\arrow[r,"\prod_l\rho_l"]
			\arrow[d]
			&
			\prod_{l=1}^3\mathrm{SO}(W_l)
			\arrow[d,hook]
			\\
			H_{k^s}
			\arrow[r,"\sim"]
			&
			\mathrm{Spin}(W)
			\arrow[r,"\pi_{k^s}"]
			&
			\mathrm{SO}(W).
		\end{tikzcd}
		\end{equation} \addtocounter{numberingbase}{1}
		Recall that $\mathrm{Spin}_4\simeq \mathrm{SL}_2\times\mathrm{SL}_2$ and that under this identification the central isogeny $\mathrm{SL}_2\times \mathrm{SL}_2\to \mathrm{SO}_4$ has kernel the diagonal $\mu_2$. Thus $\rho_l$ is the spin covering $\mathrm{Spin}(W_l)\to\mathrm{SO}(W_l)$, for $l=1,2,3$. By the universal property of simply connected coverings, we obtain a commutative diagram
		\begin{equation}\label{eq:square}
		\begin{tikzcd}
			(\widetilde J_h)_{k^s}\arrow[r,"\sim"]\arrow[d,"\varphi_{k^s}"] & \prod_{l=1}^3\mathrm{Spin}(W_l) \arrow[d]\\
			H_{k^s} \arrow[r,"\sim"] & \mathrm{Spin}(W),
		\end{tikzcd}
		\end{equation}
		where the right vertical map is induced by the orthogonal decomposition $W=W_1\perp W_2\perp W_3$. It follows from \eqref{eq:rectangle} that $\mathrm{ker}(\varphi_{k^s})$ is contained in the product of the three central subgroups $\mu_2\subset\mathrm{Spin}(W_l)$. The right vertical map in \eqref{eq:square} is induced by multiplication in the Clifford
		algebra of $W$ (see \cite[\S 13.B]{KMRT98}) and so, since the central subgroup
		$\mu_2\subset\mathrm{Spin}(W_l)$ consists of
		scalar elements, the restriction $\mu_2^3\to \mu_2$ is given by multiplication, that is, it is the norm map. Since the Galois action permutes the three factors according to its action on $E\otimes_kk^s\simeq (k^s)^3$, we conclude that $\ker\varphi\simeq R^1_{E/k}\mu_2$, which proves (i).	
		
		\m 
		By (i), for every $l=1,2,3$, the central subgroup $\mu_2\subset\mathrm{Spin}(W_l)$ maps isomorphically onto
		$\ker(\mathrm{Spin}(W)\to\mathrm{SO}(W))$, and hence $\ker\pi\subset J_h$. Since $\pi(J_h)=J_{h,0}$ by definition, we obtain $J_h=\pi^{-1}(J_{h,0})$, proving (ii).
		
	\m Define
		\[
		P\ce\{(x,y)\in\mathrm{GL}_1(D)\times
		\mathrm{GL}_1(C):
		\mathrm{Nrd}_{D/E}(x)\mathrm{Nrd}_{C/E}(y)=1\},
		\]
		and let $\eta\colon R_{E/k}P\to J_{h,0}$ be the homomorphism induced by multiplication.
		With the notation of (i), we have a commutative square 	
		\[
		\begin{tikzcd}
			(R_{E/k}P)_{k^s}
			\arrow[r,"\sim"]
			\arrow[d,"\eta_{k^s}"]
			&
			\prod_{l=1}^3 P_l
			\arrow[d,"\prod_l\eta_l"]
			\\
			(J_{h,0})_{k^s}
			\arrow[r,"\sim"]
			&
			\prod_{l=1}^3\mathrm{SO}(W_l),
		\end{tikzcd}
		\]
		where for $l=1,2,3$ we set
		\[P_l\ce\{(x,y)\in\mathrm{GL}_1(D_l)\times\mathrm{GL}_1(C_l):\mathrm{Nrd}(x)\mathrm{Nrd}(y)=1\}.\]
		The homomorphism $\eta_l\colon P_l\to\mathrm{SO}(W_l)$ is induced by multiplication. It follows that $\eta_l$ restricts to $\rho_l$ on $\mathrm{SL}_1(D_l)\times\mathrm{SL}_1(C_l)$, so that in particular $\eta_l$ is surjective, and that $\mathrm{ker}(\eta_l)$ consists of the pairs $(\lambda,\lambda^{-1})$ with $\lambda\in\mathbb G_m$. Thus $\eta$ fits into a short exact sequence
		\[1\longrightarrow R_{E/k}\mathbb G_m \longrightarrow R_{E/k}P \xrightarrow{\phantom{x}\eta\phantom{x}} J_{h,0}\longrightarrow 1.\]
		By Hilbert's Theorem~90, we have $H^1(k,R_{E/k}\mathbb G_m)=1$, and the conclusion follows.  \qedhere
	\eenum
	\end{proof}
	
	\addtocounter{numberingbase}{1}
	\begin{proposition}\label{prop:J}
		For every $h\in\Omega(k)$, one has $J_h(k)/R=1$.
	\end{proposition}
	
	\begin{proof}
		Gille's theorem \cite[Theorem~6(b)]{Gil01} applied to the sequence of Lemma~\ref{lem:Jstructure} gives an exact sequence
		\[
		\widetilde J_h(k)/R\longrightarrow J_h(k)/R
		\longrightarrow H^1(k,R^1_{E/k}\mu_2)/R\longrightarrow1.
		\]
		The group $\widetilde J_h$ is rational because quaternion norm-one groups are smooth affine quadrics with a rational point and rationality is preserved by Weil restriction. In particular, $\widetilde J_h(k)/R=1$.
		
		It remains to prove that $H^1(k,R^1_{E/k}\mu_2)/R$ is trivial. Since $[E:k]=3$, we have $R_{E/k}\mu_2\simeq R^1_{E/k}\mu_2\times\mu_2$, and it suffices to show $H^1(k, R_{E/k}\mu_2)/R=1$. By Kummer theory and Shapiro's lemma, $H^1(k,R_{E/k}\mu_2)
		\simeq E^\times/E^{\times2}$. 	Every class $[a]\in E^\times/E^{\times2}$ is $R$-equivalent to $1$ via the class of $1-t+ta$, as desired.
	\end{proof}
	\addtocounter{numberingbase}{-2}

	\csub[Proof of Theorem~\ref{thm:main}]
	\addtocounter{numberingbase}{-4}
\renewcommand{\thenumberingbase}{\thesection.4.\arabic{numberingbase}}

\vspace{-20pt}
	\begin{lemma}\label{lem:common-centralizer}
		 Let $\Omega\subset H$ be as in Lemma~\uref{lem:generic-geometry}. For every $h\in\Omega(k)$, we have 
		 \[
		 M_0(k)=(J_{h,0}(k)\cap M_0(k))\cdot\pi(M(k)).
		 \]
	\end{lemma}
	
	\begin{proof}
		Define $Q_{h,0}\ce J_{h,0}\cap M_0$, and let $m\in M_0(k)$ and
		$t\ce \mathrm{Nrd}_{B/K}(m)\in T_0(k)$. By Hilbert's Theorem~90, there exists
		$u\in K^\times$ such that $t=u/\bar u$, where the bar denotes
		the nontrivial automorphism of $K/k$, and we set $n=u\bar u\in k^\times$.
		By Lemma~\ref{lem:norm}, we may find $c\in C^\times$ with
		$\mathrm{Nrd}_{C/E}(c)=n^{-1}$. Since $\mathrm{Nrd}_{D/E}(u)=n$,
		Lemma~\ref{lem:Jstructure}(iii) gives $uc\in J_{h,0}(k)$.
		Since $u$ and $c$ commute with $i$, we have $uc\in Q_{h,0}(k)$. 
		
		On the other hand, since $Z_D(i)=E[i]\simeq E\otimes_kK$, by Lemma~\ref{lem:block} and the double-centralizer theorem we have \[Z_{B}(E)=Z_{Z_A(E)}(i)=Z_D(i)\otimes_EC=E[i]\otimes_E C,\]
		and therefore
		\[
		\mathrm{Nrd}_{B/K}(uc)
		=N_{E[i]/K}(u^2\mathrm{Nrd}_{C/E}(c))
		=N_{E[i]/K}(t)=t^3.
		\]
		Thus $\mathrm{Nrd}_{B/K}((uc)^{-1}m)=t^{-2}$ is a square in $T_0(k)$, and so
		Lemma~\ref{lem:centralizer} gives $(uc)^{-1}m\in\pi(M(k))$, that is, $m\in Q_{h,0}(k)\cdot\pi(M(k))$, as desired.
	\end{proof}
	
	\begin{proposition}\label{prop:open}
		Let $\Omega\subset H$ be as in Lemma~\uref{lem:generic-geometry}. There exists a nonempty open subscheme $\Omega\subset H$ such that $\Omega(k)\subset RH(k)$.
	\end{proposition}
	
	\begin{proof}
		Fix $h\in\Omega(k)$, set $a=\pi(h)$ and $j=aia^{-1}$, and consider the corresponding algebras $E,D,C$.
		Then $x\ce i+j$ satisfies $x^2=z+2d\in E^\times$ and $xi=jx$. Therefore $x\in D^\times$ and $xix^{-1}=j$. By Lemma~\ref{lem:norm}, there exists $c\in C^\times$ with
		$\mathrm{Nrd}_{C/E}(c)=\mathrm{Nrd}_{D/E}(x)^{-1}$, and we define $r\ce xc$.
		Lemma~\ref{lem:Jstructure}(iii) then gives $r\in J_{h,0}(k)$ and $rir^{-1}=j$, so that $r^{-1}a\in M_0(k)$. Lemma~\ref{lem:common-centralizer} now implies
		\[
		a\in J_{h,0}(k)\cdot M_0(k)=J_{h,0}(k)\cdot \pi(M(k)).
		\]
		By Lemma~\ref{lem:Jstructure}(ii), we deduce that $h\in J_h(k)\cdot M(k)\subset RH(k)$,
		where the inclusion follows from Proposition~\ref{prop:J} and Lemma~\ref{lem:M}.
	\end{proof}
	
	\begin{proof}[Proof of Theorem~\uref{thm:main}]
		By Remark~\ref{rem:char2}, we may assume that $\mathrm{char}(k)\neq 2$, and by Proposition~\ref{prop:reduction}, it suffices to show $H(k)/R=1$. Proposition~\ref{prop:open} gives a dense open subset $\Omega\subset H$ such that $\Omega(k)\subset RH(k)$. Since $k$ is infinite and $H$ is reductive, hence unirational, $H(k)$ is Zariski-dense in $H$, so that for each $h\in H(k)$ we may choose $x\in\Omega(k)\cap h^{-1}\Omega(k)$. Then both $x$ and $hx$ belong to $RH(k)$, and hence so does $h=(hx)x^{-1}$.
	\end{proof}

	\begin{remark}
	A recent preprint \cite{Ope26} claims to establish the Margulis--Platonov conjecture over number fields. This would imply that $G(k)/R=1$ for every semisimple simply connected group $G$ over a number field $k$ (this was already known when $G$ has no almost simple factor anisotropic of type $E_6$, see \cite[Theorem~8.3]{Gil09}), and hence that $G(k)/R$ is commutative for every reductive group $G$ over a number field $k$, answering \cite[Question~7.5(1)]{Gil09} in the case of number fields.
	\end{remark}

\end{appendix}

%%%%%%%%%%%%%%%%%%%%%%%%%%%%%%%%%%

%\subfile{}

\end{document}